\documentclass[11pt,reqno]{elsarticle}

\makeatletter
\def\ps@pprintTitle{%
 \let\@oddhead\@empty
 \let\@evenhead\@empty
 \def\@oddfoot{}%
 \let\@evenfoot\@oddfoot}
\makeatother

\usepackage{amsmath,amsthm,amscd,wrapfig,amsfonts,amssymb,mathtools}
\usepackage{enumitem,comment,mathrsfs}
\usepackage{booktabs,multirow,lscape,cool,bm}
\usepackage{url,parskip,caption}
\usepackage[top=1.0in, bottom=1.0in, left=1.0in, right=1.0in]{geometry}

\newtheorem{remark}{Remark}

\makeatletter
\g@addto@macro\normalsize{%
  \setlength\abovedisplayskip{.4em}
  \setlength\belowdisplayskip{.4em}
  \setlength\abovedisplayshortskip{.4em}
  \setlength\belowdisplayshortskip{.4em}
}

\begin{document}

\begin{frontmatter}


\title
{
Phase function methods for second order linear ordinary differential equations with turning points
}

\begin{abstract}
It is well known that second order linear ordinary differential equations with slowly varying coefficients admit slowly varying 
phase functions. This observation is the basis of the Liouville-Green method and many other techniques for the asymptotic approximation of the solutions of such equations.
More recently,  it was exploited by the author
to develop  a highly efficient solver  for second order linear ordinary differential equations whose solutions
are oscillatory.    In many cases of interest, that algorithm
 achieves near optimal accuracy  in time independent of the frequency of oscillation
of the solutions.    Here we show that, after minor modifications,  it
also allows for the efficient solution of second order differential equations which have turning points.
That is, it is  effective in the case
of equations whose solutions are oscillatory in some regions and behave
like linear combinations of increasing and decreasing exponential functions in others.
We present the results of numerical  experiments demonstrating the properties of our method,
including some which show that it can used to evaluate many classical special functions
in time independent of the parameters on which they depend.

\end{abstract}


\author{James Bremer}
\ead{bremer@math.toronto.edu}
\address{Department of Mathematics, University of Toronto}

\end{frontmatter}

\begin{abstract}
{It is well known that second order linear ordinary differential equations with slowly varying coefficients admit slowly varying 
phase functions. This observation is the basis of the Liouville-Green method and many other techniques for the asymptotic approximation of the solutions of such equations.
More recently,  it was exploited by the author
to develop  a highly efficient solver  for second order linear ordinary differential equations whose solutions
are oscillatory.    In many cases of interest, that algorithm
 achieves near optimal accuracy  in time independent of the frequency of oscillation
of the solutions.    Here we show that, after minor modifications,  it
also allows for the efficient solution of second order differential equations which have turning points.
That is, it is  effective in the case
of equations whose solutions are oscillatory in some regions and behave
like linear combinations of increasing and decreasing exponential functions in others.
We present the results of numerical  experiments demonstrating the properties of our method,
including some which show that it can used to evaluate many classical special functions
in time independent of the parameters on which they depend.
}
{
ordinary differential equations; fast algorithms; special functions.
}
\end{abstract}



\begin{section}{Introduction}

It has long been known that the solutions of 
second order linear ordinary differential equations
with slowly varying coefficients can be approximated via slowly varying
phase functions.    This principle is the basis of many asymptotic
techniques, including the venerable Liouville-Green method.  
If $q$ is smooth and strictly positive on the interval $[a,b]$, then
\begin{equation}
u_0(t) = \frac{\cos\left(\alpha_0(t)\right)}{\sqrt{\alpha_0'(t)}}
\ \ \mbox{and}\ \ \
v_0(t) = \frac{\sin\left(\alpha_0(t)\right)}{\sqrt{\alpha_0'(t)}},
\label{introduction:lgbasis}
\end{equation}
where $\alpha_0$ is defined via
\begin{equation}
\alpha_0(t) = \int_a^t \sqrt{q(s)}\ ds,
\label{introduction:lgphase}
\end{equation}
are a pair of Liouville-Green approximates for the second order linear
ordinary differential equation
\begin{equation}
y''(t) + q(t) y(t) = 0, \ \ a < t < b.
\label{introduction:ode}
\end{equation}
When $q$ is slowly varying, so too is 
the function defined via (\ref{introduction:lgphase}),
and this is the case regardless of the magnitude of $q$.   By contrast, the solutions 
of (\ref{introduction:ode}) become increasing oscillatory as $q$ 
grows in magnitude.  
For a careful discussion of the Liouville-Green method,
including rigorous error bounds for the approximates  (\ref{introduction:lgbasis}),
we refer the reader to  Chapter~6 of \cite{Olver}.

If $\alpha$ is a phase function for (\ref{introduction:ode})
--- so that
\begin{equation}
u(t) = \frac{\cos\left(\alpha(t)\right)}{\sqrt{\alpha'(t)}}
\ \ \mbox{and}\ \ \
v(t) = \frac{\sin\left(\alpha(t)\right)}{\sqrt{\alpha'(t)}}
\label{introduction:basis}
\end{equation}
are solutions of (\ref{introduction:ode}) --- then $\alpha'$ satisfies Kummer's differential equation
\begin{equation}
q(t) - (\alpha'(t))^2 + \frac{3}{4} \left(\frac{\alpha''(t)}
{\alpha'(t)}\right)^2
- \frac{\alpha'''(t)}{2\alpha'(t)} = 0,\ \ a <t <b.
\label{introduction:kummer}
\end{equation}
The Liouville-Green phase (\ref{introduction:lgphase}) is the (crude) approximation
obtained by deleting the expression
\begin{equation}
\frac{3}{4} \left(\frac{\alpha''(t)}
{\alpha'(t)}\right)^2
- \frac{\alpha'''(t)}{2\alpha'(t)} = 0
\end{equation}
from (\ref{introduction:kummer}) and solving the resulting equation.  
In \cite{SpiglerZeros},  an iterative scheme for refining the 
Liouville-Green phase is introduced.  In order to avoid the complications 
which square roots bring, it constructs a sequence $\phi_0, \phi_1, \phi_2, \ldots$
of asymptotic approximations to a solution of 
the differential equation
\begin{equation}
\phi(t) = q(t) -\frac{1}{4}\frac{\phi''(t)}{\phi(t)} + \frac{5}{16} \left(\frac{\phi'(t)}{\phi(t)}\right)^2
\label{introduction:spigler}
\end{equation}
satisfied by functions of the form $\phi(t) = \left(\alpha'(t)\right)^2$,
where $\alpha'$ is the derivative of a phase function for (\ref{introduction:ode}).
The formulas
\begin{equation}
\left\{
\begin{aligned}
\phi_0(t) &= q(t) \\
\phi_{n+1}(t) &= q(t) -\frac{1}{4}\frac{\phi_n''(t)}{\phi_n(t)} + \frac{5}{16} \left(\frac{\phi_n'(t)}{\phi_n(t)}\right)^2
\end{aligned}
\right.
\label{introduction:spigleriter}
\end{equation}
defining the scheme of \cite{SpiglerZeros} make clear that  when $q$ and its derivatives are slowly varying, 
the asymptotic approximations it produces will be as well.
Error bounds  which hold under various assumptions on the form of $q$ 
are given in \cite{SpiglerPhase1,SpiglerPhase2,SpiglerZeros}.  In the important case in which 
\begin{equation}
q(t) = \lambda^2 q_0(t)
\label{introduction:q}
\end{equation}
with $q_0$ smooth and positive,  there exist solutions $u$, $v$ 
of (\ref{introduction:ode}) such that if $\alpha_n$ is the approximate phase
function obtained from the iterate $\phi_n$ of the scheme (\ref{introduction:spigleriter}), then
\begin{equation}
\begin{aligned}
u(t) &= \frac{\cos\left(\alpha_n(t)\right)}{\sqrt{\alpha_n'(t)}} + \mathcal{O}\left(\frac{1}{\lambda^{n+1}}\right)
\ \ \mbox{as} \ \ \lambda\to\infty\ \ \mbox{and}\\[1.2em]
v(t) &= \frac{\sin\left(\alpha_n(t)\right)}{\sqrt{\alpha_n'(t)}} + \mathcal{O}\left(\frac{1}{\lambda^{n+1}}\right)
\ \ \mbox{as} \ \ \lambda\to\infty.
\end{aligned}
\label{introduction:spigleruv}
\end{equation}

The approximations generated by the scheme of \cite{SpiglerZeros} are closely connected to the 
(perhaps more familiar) classical ones obtained from the Riccati equation
\begin{equation}
r'(t) + (r(t))^2 + q(t) = 0
\label{introduction:riccati}
\end{equation}
satisfied by the logarithmic derivatives of the solutions of (\ref{introduction:ode}).
We note that the solutions of (\ref{introduction:riccati}) are related to those
of (\ref{introduction:kummer}) via the formula
\begin{equation}
r(t) = i \alpha'(t) - \frac{1}{2} \frac{\alpha''(t)}{\alpha'(t)}.
\label{introduction:r}
\end{equation}
In the event that $q$ takes the form (\ref{introduction:q}), inserting the ansatz
\begin{equation}
r(t) = \sum_{k=0}^n \lambda^{1-k} r_k(t) 
\label{introduction:ricciter1}
\end{equation}
into (\ref{introduction:riccati}) and solving the equations which result from collecting
like powers of $\lambda$ leads to the formulas
\begin{equation}
\begin{aligned}
r_0(t) &= i \sqrt{q_0(t)}, \\
r_1(t) &= -\frac{q_0'(t)}{4q_0(t)},\\
\vdots\\
r_{n}(t) &= -\frac{i}{2\sqrt{q_0(t)}} \left(r_{n-1}'(t) +\sum_{j=0}^{n-1} r_j(t) r_{n-j}(t) \right).
\end{aligned}
\label{introduction:ricciter2}
\end{equation}
It is a classical result (a proof of which can be found in Chapter~7 of \cite{Miller}) 
that there exist solutions $u$, $v$ of (\ref{introduction:ode}) such that
\begin{equation}
\begin{aligned}
u(t) &= \exp\left( i \sum_{n=0}^n  \lambda^{1-n} \int_a^t r_n(s)\ ds \right) + \mathcal{O}\left(\frac{1}{\lambda^{n+1}}\right) \ \ \mbox{as} \ \ \lambda\to\infty
\ \ \mbox{and}\\[1.2em]
v(t) &= \exp\left( -i\sum_{n=0}^n  \lambda^{1-n} \int_a^t r_n(s)\ ds \right) + \mathcal{O}\left(\frac{1}{\lambda^{n+1}}\right) \ \ \mbox{as} \ \ \lambda\to\infty.
\end{aligned}
\label{introduction:wkbapprox}
\end{equation}
Since the formula defining $r_{n+1}$ depends on all of the previous iterates
$r_0,r_1,\ldots,r_n$  whereas  the $(n+1)^{st}$ iterate $\phi_{n+1}$ 
 depends only on the  $n^{th}$ iterate $\phi_n$, 
the approximations produced by the scheme (\ref{introduction:spigleriter})
are typically much simpler than those which result from the classical scheme.

While asymptotic methods such as (\ref{introduction:spigleriter}) and
(\ref{introduction:ricciter1}), (\ref{introduction:ricciter2})
are extremely useful for generating {\it symbolic} expressions which represent
the solutions of second order linear ordinary differential equations,
they leave much to be desired as numerical methods.   
Such schemes suffer from at least  two serious problems:

\begin{enumerate}
\item

Like all asymptotic methods,  there is a limit to the accuracy which they can obtain,
and this limit often falls short of the level of accuracy  indicated by the condition number
of the problem.  Moreover, in the case of  asymptotic methods for  (\ref{introduction:ode}),
the obtainable accuracy generally depends in a complicated way on behavior of the function $q$ 
and its derivatives, thus making it difficult to control numerical errors.

\item
The higher order approximations constructed by 
these techniques depend on higher order derivatives of $q$, which cannot be 
calculated numerically without significant loss of accuracy.


\end{enumerate}

In \cite{BremerKummer}, an algorithm for constructing nonoscillatory
phase functions which represent the solutions of equations
of the form (\ref{introduction:ode}) when $q$ is smooth and strictly positive
is described.   It operates by solving Kummer's equation (\ref{introduction:kummer})
numerically.  Most of the solutions of Kummer's equation are oscillatory,
and the principal challenge addressed by the algorithm
of \cite{BremerKummer} is the identification of the values of the first two
derivatives of a nonoscillatory phase function at a point on the interval $(a,b)$.  
Once this has been done, (\ref{introduction:kummer}) 
can be solved numerically using any method which  applies to stiff ordinary differential equations.
The algorithm of \cite{BremerKummer} does not require knowledge of  the
derivatives of $q$ and calculates the phase function $\alpha$ to  near machine precision
regardless of the magnitude of $q$.

A theorem which shows the existence of a nonoscillatory
phase function for (\ref{introduction:ode}) in the case in which
$q$ is of the form (\ref{introduction:q})  appears in a companion paper \cite{BremerRokhlin} to \cite{BremerKummer}.  
It applies when the function $p(x) = \widetilde{p}(t(x))$, where
$p(t)$ is defined via
\begin{equation}
\widetilde{p}(t) = 
\frac{1}{q_0(t)} \left(
\frac{5}{4}\left(\frac{q_0'(t)}{q_0(t)}\right)^2 - \frac{q_0''(t)}{q_0(t)}
\right)
=
4 \left(q_0(t)\right)^{\frac{1}{4}}
\frac{d}{dt} \left(
\frac{1}{\left(q_0(t)\right)^{\frac{1}{4}}}
\right)
\end{equation}
and $t(x)$ is the inverse function of 
\begin{equation}
x(t) = \int_a^t \sqrt{q_0(s)}\ ds,
\label{introduction:xt}
\end{equation}
has a rapidly decaying Fourier transform.  More explicitly, the theorem
asserts that if the Fourier transform of $p$ satisfies a bound of the form
\begin{equation}
\left|\widehat{p}(\xi)\right|\leq \Gamma \exp\left(-\mu\left|\xi\right|\right),
\end{equation}
then there exist functions $\nu$ and $\delta$ such that 
\begin{equation}
\left|\nu(t)\right| \leq  \frac{\Gamma}{2\mu} \left(1 + \frac{4\Gamma}{\lambda}\right) \exp(-\mu \lambda),
\end{equation}
\begin{equation}
\left|\widehat{\delta}(\xi)\right| \leq \frac{\Gamma}{\lambda^2}\left(1+\frac{2\Gamma}{\lambda}\right)\exp(-\mu|\xi|)
\end{equation}
and
\begin{equation}
\alpha(t) = \lambda \sqrt{q_0(t)} \int_a^t \exp\left(\frac{\delta(u)}{2}\right)\ du
\end{equation}
is a phase function for 
\begin{equation}
y''(t) + \lambda^2 \left(q_0(t) + \frac{\nu(t)}{4\lambda^2} \right) y(t) = 0.
\end{equation}
The definition of the function $p(x)$ is ostensibly quite complicated;
however, $p(x)$ is, in fact, simply equal to twice the Schwarzian derivative of the inverse function
$t(x)$ of (\ref{introduction:xt}).
This theorem ensures that, even for relatively modest values of $\lambda$, the phase
function constructed by \cite{BremerKummer} is nonoscillatory. 
When $\lambda$ is small, the algorithm of \cite{BremerKummer} is not guaranteed
to produce a phase function which is nonoscillatory, but, in stark contrast to asymptotic methods,
in all cases it results in one which approximates a solution of (\ref{introduction:kummer}) with near machine 
precision accuracy.

Here, we describe a variant of the algorithm of \cite{BremerKummer} 
which allows for the numerical solution 
of second order linear ordinary differential equations with turning points.
We focus on the case in which the coefficient $q$ in (\ref{introduction:ode}) is a smooth function
on the interval $(a,b)$ with a single zero at the point $c \in (a,b)$,
and further assume that 
\begin{equation}
q(t) \sim C (t-c)^k \ \ \mbox{as} \ \  t\to c
\label{introduction:qasym}
\end{equation}
with $k$ a positive integer.
Most equations with multiple turning points  can be handled by the repeated
use of the algorithm described here (and we consider two such examples 
in the numerical experiments discussed in this article).  When $q$ is nonpositive
on the entire interval $(a,b)$, alternate methods are indicated
(for example, the technique used in \cite{BremerRadial}
to solve second order differential equations in the nonoscillatory regime).

Two problems arise when the algorithm of \cite{BremerKummer} is applied to second
order differential equations of this type.
First,  Kummer's equation (\ref{introduction:kummer}) and the closely-related
Riccati equation (\ref{introduction:riccati}) encounter numerical
difficulties in the nonoscillatory regime, where the values of $\alpha'$ 
are small.   More explicitly, when these equations are solved numerically,
the obtained values of $\alpha'$ only approximate it with \emph{absolute} accuracy
on the order of machine precision in the nonoscillatory region.  Since $\alpha'$ is small
there, this means  that the \emph{relative} accuracy of these approximations
is extremely poor.  That such difficulties are encountered is 
unsurprising given that $\alpha'$ appears in the denominator of two terms in Kummer's equation, 
and considering the form (\ref{introduction:q}) of the  solutions of Riccati's equation.
We address this difficulty by solving Appell's equation rather than Kummer's equation.
Appell's equation is a certain third order linear ordinary differential
equation satisfied by the product of any two solutions of (\ref{introduction:ode}),
including the reciprocal of $\alpha'$, which is the sum of the squares
of the two functions appearing in (\ref{introduction:basis}).
As the experiments discussed in  this paper show, there is no difficulty in calculating 
$\alpha'$ with high relative accuracy in both the oscillatory and nonoscillatory regimes by solving Appell's 
equation numerically.

The second difficulty addressed by the algorithm of this paper is that, 
in the case of turning points of even orders,
there need not exist a slowly varying phase function which extends across
the turning point.  We include a simple analysis of the second order differential equation
\begin{equation}
y''(t) + t^{k} y(t) = 0
\label{introduction:evenode}
\end{equation}
in the case in which $k$ is an even integer to demonstrate this.
To overcome this problem, we simply  construct two nonoscillatory 
phase functions, each defined only on one side of the turning point.
For turning points of odd order, a single phase function suffices.

The remainder of this article is structured as follows.
Section~\ref{section:phase} reviews certain basic facts regarding
phase functions for second order linear ordinary differential equations.
In Section~\ref{section:turning}, we analyze the second
order differential equation (\ref{introduction:evenode})
in order to show that we cannot expect 
  nonoscillatory phase functions to extend
 across turning points of even order.
Section~\ref{section:algorithm} details our numerical algorithm.
Numerical experiments conducted to demonstrate its properties are
discussed in Section~\ref{section:experiments}.  We close
with a few brief remarks in Section~\ref{section:conclusions}.

\end{section}

\begin{section}{Phase functions for second order linear ordinary differential equations}
\label{section:phase}

In this section, we review several basic facts regarding phase functions
for second order linear ordinary differential equations.  We first discuss
the case of a general second order linear ordinary differential equation
in Subsection~\ref{section:phase:general}, and
 then we consider the effect of applying the standard transform
which reduces  such equations  to their  so-called ``normal forms''
in Subsection~\ref{section:phase:normal}.

\begin{subsection}{The general case}
\label{section:phase:general}

If $y(t) = \exp(r(t))$ satisfies the second order differential equation
\begin{equation}
y''(t) + p(t) y'(t) + q(t) y(t) = 0, \ \ a < t < b,
\label{phase:ode}
\end{equation}
then it can be easily seen that r(t) solves the Riccati equation
\begin{equation}
r''(t) + (r'(t))^2 + p(t) r'(t) + q(t) = 0,\ \ a < t <  b.
\label{phase:riccati}
\end{equation}
Letting  $r(t) = i \alpha(t) + \beta(t)$  yields the system
of ordinary differential equations
\begin{equation}
\left\{
\begin{aligned}
q(t) - (\alpha'(t))^2 + p(t) \beta'(t) + (\beta'(t))^2+ \beta''(t)  &=0 \\
p(t) \alpha'(t) + 2 \alpha'(t) \beta'(t) + \alpha''(t) &= 0.
\end{aligned}
\right.
\label{phase:system}
\end{equation}
The second of these can be rearranged as 
\begin{equation}
\beta'(t) = -\frac{1}{2} \frac{\alpha''(t)}{\alpha'(t)} - \frac{p(t)}{2},
\label{phase:beta}
\end{equation}
and when this expression is inserted into the first equation in (\ref{phase:system}) 
we arrive at 
\begin{equation}
q(t)- \frac{(p(t))^2}{4} - \frac{p'(t)}{2} - (\alpha'(t))^2  + \frac{3}{4} \left(\frac{\alpha''(t)}{\alpha'(t)}\right)^2
- \frac{\alpha'''(t)}{2\alpha'(t)} = 0,\ \ a <t <b.
\label{phase:kummer}
\end{equation}
It is clear  that if  $\alpha$ does not vanish on $(a,b)$ and satisfies (\ref{phase:kummer})
there, then
\begin{equation}
r(t) = i \alpha(t) -\frac{1}{2} \log\left(\alpha'(t)\right) 
- \int p(t)\ dt,
\label{phase:rsol}
\end{equation}
where the implicit constant of integration and the choice of branch cut for
the logarithm are irrelevant,
is a solution of the Riccati equation (\ref{phase:riccati}) and 
\begin{equation}
\sqrt{\frac{\omega(t)}{\alpha'(t)}} \exp(i \alpha(t))
=
 \sqrt{\frac{\omega(t)}{\alpha'(t)}} \cos(\alpha(t))
+ i  \sqrt{\frac{\omega(t)}{\alpha'(t)}} \sin(\alpha(t)),
\end{equation}
where
\begin{equation}
\omega(t) = \exp\left(-\int p(t)\ dt\right),
\label{phase:wron}
\end{equation}
is a solution of the original ordinary differential equation (\ref{phase:ode}).
We note that Abel's identity implies that the Wronskian of any pair of solutions of (\ref{phase:ode})
is a constant multiple of (\ref{phase:wron}).

Equation~(\ref{phase:kummer}) is known as Kummer's equation, after E.~E.~Kummer who
considered it in \cite{Kummer} and we  refer to its solutions 
 as phase functions for the ordinary differential equation (\ref{phase:ode}).
When the coefficients $p$ and $q$ are real-valued,
 $\alpha(t)$ is also real-valued, and the functions
\begin{equation}
u(t) = \sqrt{\frac{\omega(t)}{\alpha'(t)}} \cos(\alpha(t))
\ \ \mbox{and} \ \ \ 
v(t) = \sqrt{\frac{\omega(t)}{\alpha'(t)}} \sin(\alpha(t))
\label{phase:basis}
\end{equation}
are linearly independent real-valued solutions of (\ref{phase:ode}) which form a basis in its 
space of solutions.

If $\alpha$ is a phase function and $u$ and $v$
are as in (\ref{phase:basis}), then a straightforward 
computation shows that the modulus function
\begin{equation}
m(t) = (u(t))^2 + (v(t))^2 = \frac{\omega(t)}{\alpha'(t)}
\end{equation}
satisfies the differential equation
\begin{equation}
m'''(t) + 3p(t) m''(t) + (2 (p(t))^2 + p'(t) + 4 q(t))m'(t) +
(4p(t)q(t) + 2q'(t) ) m(t) = 0.
\label{phase:appell}
\end{equation}
%
We refer to (\ref{phase:appell}) as Appell's equation in light of the article
\cite{Appell}.  

Given any pair $u, v$  of solutions of (\ref{phase:ode}) whose Wronskian $\omega(t)$
and the associated modulus function $m(t) = (u(t))^2 + (v(t))^2$ are
nonzero on $(a,b)$,  it can be shown by a straightforward calculation that
\begin{equation}
\alpha'(t) = \frac{\omega(t)}{(u(t))^2+(v(t))^2}
\label{phase:alphap}
\end{equation}
satisfies (\ref{phase:kummer}) on that interval.  It follows that any antiderivative of $\alpha'$ 
is a phase function for (\ref{phase:ode}) on $(a,b)$.  
Requiring that (\ref{phase:basis}) holds
determines $\alpha$  up to an additive constant which is an integral multiple of $2\pi$, but further restrictions
are required to determine it uniquely.
We will nonetheless, by a slight abuse of terminology, refer to  ``the phase function generated
by the pair of solutions $u,v$''  with the understanding that the choice of constant
is of no consequence.

\end{subsection}

\begin{subsection}{Reduction to normal form}
\label{section:phase:normal}

If $y$ solves (\ref{phase:ode}), then
\begin{equation}
\tilde{y}(t) = \exp\left(-\frac{1}{2} \int p(t)\ dt \right) y(t) 
\label{phase:transform}
\end{equation}
satisfies
\begin{equation}
\tilde{y}''(t) +  \tilde{q}(t) \tilde{y}(t) = 0, \ \ a < t < b,
\label{phase:normal}
\end{equation}
where
\begin{equation}
\tilde{q}(t) = q(t) - \frac{p(t)}{2} - \frac{(p(t))^2}{4}.
\label{phase:qtilde}
\end{equation}
Equation~(\ref{phase:normal}) is often called the normal form of (\ref{phase:ode}).
The effects of this transformation on the various formulas discussed above
are easy to discern.  In particular, Riccati's equation becomes
\begin{equation}
\tilde{r}''(t) + (\tilde{r}'(t))^2 + \tilde{q}(t) = 0,
\label{phase:norm_riccati}
\end{equation}
Kummer's equation now takes the form
\begin{equation}
\tilde{q}(t) - (\tilde{\alpha}'(t))^2 + \frac{3}{4} \left(\frac{\tilde{\alpha}''(t)}
{\tilde{\alpha}'(t)}\right)^2
- \frac{\tilde{\alpha}'''(t)}{2\tilde{\alpha}'(t)} = 0,\ \ a <t <b,
\label{phase:norm_kummer}
\end{equation}

Appell's equation becomes
\begin{equation}
\tilde{m}'''(t)  + 4 \tilde{q}(t) \tilde{m}'(t) + 2 \tilde{q}'(t) \tilde{m}(t) = 0,
\end{equation}
and it is the functions
\begin{equation}
\tilde{u}(t) = \frac{\cos(\tilde{\alpha}(t))}{\sqrt{\tilde{\alpha}'(t)}}
\ \ \mbox{and} \ \ \ 
\tilde{v}(t) = \frac{\sin(\tilde{\alpha}(t))}{\sqrt{\tilde{\alpha}'(t)}}
\label{phase:norm_sols}
\end{equation}
which form a  basis in the space of solutions of (\ref{phase:normal}).

We observe that when (\ref{phase:qtilde}) is inserted into (\ref{phase:norm_kummer})
it simply becomes (\ref{phase:kummer}), so that Kummer's equation is unchanged
by this transformation.  Similarly, (\ref{phase:alphap}) is invariant 
under  (\ref{phase:transform}), which implies that
\begin{equation}
\alpha'(t) = \tilde{\alpha}'(t).
\end{equation}
So while the solution of Riccati's equation and the basis functions are modified
by the transformation (\ref{phase:transform}), the phase function $\alpha$ is invariant.  
It does not matter, then, whether the form (\ref{phase:ode}) or its normal form (\ref{phase:normal})
is considered.   Either set of solutions (\ref{phase:basis}) or (\ref{phase:norm_sols}) 
can be recovered once $\alpha$ has been calculated.

\end{subsection}

\end{section}

\begin{section}{Turning points of even order}
\label{section:turning}

We now consider  the equation
\begin{equation}
y''(t) + t^k y(t) = 0
\label{turning:ode}
\end{equation}
in the case in which  $k$ is an even integer.  Since the solutions of 
(\ref{turning:ode}) can be expressed as linear combinations of  Bessel functions,
we begin with a short discussions of Bessel's differential equation
in Subsection~\ref{section:turning:bessel}.  We then
exhibit a basis in the space of solutions of (\ref{turning:ode}) in
Subsection~\ref{section:turning:basis}.  Finally, in Subsection~\ref{section:turning:phase},  
we show that (\ref{turning:ode}) does not admit a phase function which is nonoscillatory
on all of the real line.


\begin{subsection}{Bessel functions}
\label{section:turning:bessel}

Standard solutions of Bessel's differential equation 
\begin{equation}
z^2 y''(z) + z y'(z) + (z^2-\nu^2)y(z) = 0
\label{turning:bessel}
\end{equation}
include the Bessel functions
\begin{equation}
J_\nu(z) = \left(\frac{z}{2}\right)^\nu \sum_{k=0}^\infty \frac{1}{\Gamma(k+1)\Gamma(k+\nu+1)}\left(\frac{z}{2}\right)^k
\label{turning:jnu}
\end{equation}
and
\begin{equation}
Y_\nu(z) = \cot(\pi\nu) J_\nu(z) - \csc(\pi\nu)J_{-\nu}(z)
\label{turning:ynu}
\end{equation}
of the first and second kinds of order $\nu$.  In general, these functions are multi-valued,
and the usual convention is to regard them as defined on the cut plane
$\mathbb{C} \setminus \left(-\infty,0\right]$.  We note, however, 
that in the case of integer  values of $\nu$, $J_\nu(z)$ is entire.

We use $\alpha_\nu^{\mbox{\tiny bes}}$ to denote the phase function for Bessel's differential
equation generated by the pair $J_\nu,\ Y_\nu$.  The following formula, which
appears in   \cite{Hartman73}, expresses the corresponding modulus function
as the Laplace transform of a Legendre function:
\begin{equation}
J_\nu^2(z) + Y_\nu^2(z)
=
 \frac{2}{\pi} \int_0^\infty \exp(-zt) P_{\nu-1/2}\left(1+\frac{t^2}{2}\right)\ dt.
\label{turning:beslaplace}
\end{equation}
Among other things, it implies that this modulus function is  
completely monotone on $(0,\infty)$.   A function $f$ is completely monotone on the 
interval $(a,b)$ provided
\begin{equation}
(-1)^n f^{(n)}(t) \geq 0 \ \ \mbox{for all} \ \ a < t <b.
\end{equation}
It is well known that $f$ is completely monotone on $(0,\infty)$ if and only if it is the 
Laplace transform of a nonnegative Borel measure  (see, for instance, Chapter~4 of \cite{Widder}).
The derivation of (\ref{turning:beslaplace}) in \cite{Hartman73} relies on a considerable
amount of machinery.  However, it can be verified simply by observing that the integral 
in (\ref{turning:beslaplace}) satisfies Appell's differential equation and then matching
the first few terms of the asymptotic expansions 
of the  left- and right-hand sides of  (\ref{turning:beslaplace}) at infinity.

Since the Wronskian of the pair $J_\nu(z),\ Y_\nu(z)$ is $\frac{2}{\pi z}$
(this fact can be found in Chapter~7 of \cite{HTFII}),
\begin{equation}
\frac{d}{dz}\alpha_\nu^{\mbox{\tiny bes}}(z)= \frac{2}{\pi z} \frac{1}{J_\nu^2(z) + Y_\nu^2(z)}.
\end{equation}
In particular, Bessel's differential equation admits a phase function which is slowly varying
in the extremely strong sense that its derivative is the  reciprocal of the product of $z$ and a completely monotone function.

\end{subsection}

\begin{subsection}{A basis in the space of solutions}
\label{section:turning:basis}

It follows from the expansions (\ref{turning:jnu}) and (\ref{turning:ynu}) that 
when $0 < \nu < 1$ the following formulas hold:
\begin{equation}
\begin{aligned}
\lim_{t\to 0^+}  J_\nu\left(2\nu t^{\frac{1}{2\nu}}\right)\sqrt{t} &= 0,\\
\lim_{t\to 0^+}  \frac{d}{dt} \left(J_\nu\left(2\nu t^{\frac{1}{2\nu}}\right)\sqrt{t}\right) &= \frac{\nu^\nu}{\Gamma(\nu+1)},\\
\lim_{t\to 0^+}  Y_\nu\left(2\nu t^{\frac{1}{2\nu}}\right)\sqrt{t} &= -\frac{\nu^{-\nu}\Gamma(\nu)}{\pi}\ \ \mbox{and} \\
\lim_{t\to 0^+}  \frac{d}{dt} \left(Y_\nu\left(2\nu t^{\frac{1}{2\nu}}\right)\sqrt{t}\right) &= -\frac{\nu^{\nu}\Gamma(-\nu)\cos\left(\pi\nu\right)}{\pi}.
\end{aligned}
\label{turning:limits}
\end{equation}
The limits in (\ref{turning:limits}) imply that functions
$u_k^{\mbox{\tiny even}}$ and  $v_k^{\mbox{\tiny even}}$ defined via
\begin{equation}
u_k^{\mbox{\tiny even}}(t) = 
\begin{cases}
\sqrt{\frac{\pi t}{2+k}}\  J_{\frac{1}{2+k}}\left(\frac{ t^{1+\frac{k}{2}}}{1+\frac{k}{2}}\right) & \mbox{if}\ \ t > 0\\
-\sqrt{\frac{-\pi t}{2+k}}\  J_{\frac{1}{2+k}}\left(\frac{(-t)^{1+\frac{k}{2}}}{1+\frac{k}{2}} \right) &  \mbox{if}\ \  t < 0
\end{cases}
\label{turning:ueven}
\end{equation}
and
\begin{equation}
v_k^{\mbox{\tiny even}}(t) = 
\begin{cases}
\sqrt{\frac{\pi t}{2+k}}\  Y_{\frac{1}{2+k}}\left(\frac{ t^{1+\frac{k}{2}}}{1+\frac{k}{2}}\right) &\mbox{if}\ \ t > 0\\
-2\cot(\pi\nu) \sqrt{\frac{-\pi t}{2+k}}\  J_{\frac{1}{2+k}}\left(\frac{ (-t)^{1+\frac{k}{2}}}{1+\frac{k}{2}}\right)
+\sqrt{\frac{-\pi t}{2+k}}\  Y_{\frac{1}{2+k}}\left(\frac{ (-t)^{1+\frac{k}{2}}}{1+\frac{k}{2}}\right)
  &\mbox{if}\ \ t < 0
\end{cases}
\label{turning:veven}
\end{equation}
are continuously differentiable at $0$.
They are obviously smooth on the set $U=\mathbb{R}\setminus\{0\}$
and it follows from direct substitution that they satisfy (\ref{turning:ode}) on $U$.
So they are, in fact, smooth on $\mathbb{R}$ and they satisfy
(\ref{turning:ode}) there.
It is also a consequence of the formulas appearing in (\ref{turning:limits})  that the Wronskian
of this pair of solutions is $1$.  

\end{subsection}

\begin{subsection}{There is no phase function which is nonoscillatory on the whole real line}
\label{section:turning:phase}

We now use the well-known asymptotic approximations
\begin{equation}
\begin{aligned}
\sqrt{\frac{\pi z}{2}} J_\nu(z) &\sim \cos\left(z-\frac{\pi\nu}{2}-\frac{\pi}{4}\right)\ \ \mbox{as} \ \ z \to\infty \ \ \mbox{and} \\[1.3em]
\sqrt{\frac{\pi z}{2}}Y_\nu(z) &\sim \sin\left(z-\frac{\pi\nu}{2}-\frac{\pi}{4}\right)\ \ \mbox{as} \ \ z \to\infty
\end{aligned}
\label{turning:asym}
\end{equation}
for the Bessel functions (which can be found in Chapter~7 of \cite{HTFII}, among many other sources) to show that there do not exist
constants $c_1$ and $c_2$ such that the modulus function
\begin{equation}
m_k^{\mbox{\tiny bes}}(t) = c_1^2 \left(u_k^{\mbox{\tiny even}}(t)\right)^2
+
c_2^2 \left(v_k^{\mbox{\tiny even}}(t)\right)^2
\end{equation}
for the pair of solutions 
\begin{equation}
c_1 u_k^{\mbox{\tiny even}}(t),\  c_2 v_k^{\mbox{\tiny even}}(t)
\end{equation}
is nonoscillatory on both of the intervals $(0,\infty)$ and $(-\infty,0)$.  

From (\ref{turning:asym}) and the definitions of 
$u_k^{\mbox{\tiny even}}$ and  $v_k^{\mbox{\tiny even}}$ given in the preceding subsection,
we see that 
\begin{equation}
m_k^{\mbox{\tiny bes}}(t) \sim
\frac{k+2}{\pi} t^{-\frac{k}{2}}
\left(
c_1^2 \sin^2\left(\frac{k\pi + 8 t^{1+\frac{k}{2}}}{8+4k}\right)
+
c_2^2 \cos^2\left(\frac{k\pi + 8 t^{1+\frac{k}{2}}}{8+4k}\right)
\right)
\ \ \mbox{as} \ \ t \to \infty.
\end{equation}
In order to prevent oscillations on the interval $(0,\infty)$, we must have $c_1=c_2$.   
Without loss of generality we may assume that $c_1=c_2=1$.  But, it then follows
from (\ref{turning:ueven}), (\ref{turning:veven}) and (\ref{turning:asym})  that
\begin{equation}
\begin{aligned}
m_k^{\mbox{\tiny bes}}(t) \sim
\frac{(k+2) (-t)^{-k/2} \csc ^2\left(\frac{\pi }{k+2}\right) \left(-4 \cos \left(\frac{\pi
   }{k+2}\right) \sin \left(\frac{4 t (-t)^{k/2}}{k+2}\right)+\cos \left(\frac{2 \pi
   }{k+2}\right)+3\right)}{2 \pi }
\end{aligned}
\label{turning:osc}
\end{equation}
as $t \to -\infty$.
Clearly, the function appearing on the right-hand side of (\ref{turning:osc}) oscillates
on $(-\infty,0)$.
\end{subsection}

\end{section}

\begin{section}{Numerical Algorithm}
\label{section:algorithm}

In this section, we describe our algorithm for solving the second order differential equation
\begin{equation}
y''(t) + q(t) y(t) = 0, \ \ a < t <b,
\label{algorithm:ode}
\end{equation}
in the case in which $q$ is smooth on $(a,b)$ with a single turning point $c$.
We will further assume that 
\begin{equation}
q(t) \sim C (t-c)^k \ \mbox{as} \ \ t\to c
\end{equation}
with $k$ a positive integer and $C>0$.  Obvious modifications to the algorithm apply in the case in
which $k$ is an odd integer and $C <0$.  When $k$ is even and $C<0$, the solutions
of (\ref{algorithm:ode}) are nonoscillatory on the whole interval $(a,b)$ and 
other methods are indicated.

The algorithm takes as input  the interval $(a,b)$; the location of the turning point $c$;
an external subroutine for evaluating the coefficient $q(t)$ and, optionally, its derivative
$q'(t)$;  a positive integer $l$  specifying the 
order of the piecewise Chebyshev expansions which are used to represent phase functions;
and a double precision number $\epsilon$ which controls the accuracy of the obtained solution.   
In the case in which the derivative of the coefficient is not specified, $q'(t)$
is evaluated using spectral differentiation.  This usually results
in a small loss of precision;    several experiments in which this effect is measured are discussed
in Section~\ref{section:experiments}.

The algorithm outputs one or more phase functions, which are represented via
$l^{th}$ order piecewise Chebyshev expansions.  By an  $l^{th}$ order piecewise Chebyshev 
expansions  on the interval $[a_0,b_0]$, we mean  a sum of the form
\begin{equation}
\begin{aligned}
&\sum_{i=1}^{m-1} \chi_{\left[x_{i-1},x_{i}\right)} (t) 
\sum_{j=0}^{l} \lambda_{ij}\ T_j\left(\frac{2}{x_{i}-x_{i-1}} t + \frac{x_{i}+x_{i-1}}{x_{i}-x_{i-1}}\right)\\
+
&\chi_{\left[x_{m-1},x_{m}\right]} (t) 
\sum_{j=0}^{l} \lambda_{mj}\ T_j\left(\frac{2}{x_{m}-x_{m-1}} t + \frac{x_{m}+x_{m-1}}{x_{m}-x_{m-1}}\right),
\end{aligned}
\end{equation}
where $a_0 = x_0 < x_1 < \cdots < x_m = b_0$ is a partition of $[a_0,b_0]$,
$\chi_I$ is the characteristic function on the interval $I$ and 
$T_j$ denotes the Chebyshev polynomial of degree $j$.
The phase functions define a basis in the space of solutions of (\ref{algorithm:ode}),
and using these basis functions it is straightforward to construct a solution of 
(\ref{algorithm:ode}) satisfying any reasonable set of boundary conditions.

When the turning point is of odd order, our algorithm
constructs a single phase function $\alpha$  such that
\begin{equation}
\begin{aligned}
u(t) = \frac{\cos\left(\alpha(t)\right)}{\sqrt{\alpha'(t)}}
\ \ \mbox{and} \ \ 
v(t) = \frac{\sin\left(\alpha(t)\right)}{\sqrt{\alpha'(t)}}
\end{aligned}
\label{algorithm:basis_odd}
\end{equation}
form a basis in the space of solutions of (\ref{algorithm:ode}).  
In this case, the algorithm returns as output  three piecewise $l^{th}$ order Chebyshev expansions, 
one representing the phase  function $\alpha(t)$, one 
representing its derivatives $\alpha'(t)$ and a third expansion representing $\alpha''(t)$.
Using these expansions, the values of the functions $u$ and $v$, as well
the values of their first derivatives, can be easily evaluated.
The phase function $\alpha$ is not always given over the whole interval $[a,b]$.  
The value of $\alpha'(t)$ decays rapidly in the nonoscillatory regime,
as $t$ decreases from $c$ to $a$, and it can become smaller
in magnitude than the smallest positive IEEE double precision number.  
In such cases, our algorithm constructs the phase function
on a truncated domain $[\tilde{a},b]$, where $\tilde{a}$ is chosen so that
the value of $\alpha'$ is close to the smallest representable IEEE double precision number
at $\tilde{a}$.  As the experiments of Subsection~\ref{section:experiments:airy} and
\ref{section:experiments:high} show, despite this potential limitation,
phase functions can represent the solutions 
of (\ref{algorithm:ode}) with high relative accuracy deep into the nonoscillatory region.
Indeed, in our experience, the results are often better than standard solvers
in this regime.

If the turning point is of even order,
then two phase functions are constructed: the ``left'' phase function $\alpha_{\mbox{\tiny left}}$ 
given on  the interval $[a,c]$ and the ``right'' phase function $\alpha_{\mbox{\tiny right}}$
given on  $[c,b]$.     The output in this case comprises six piecewise $l^{th}$ order Chebyshev expansions:
one for each of the functions $\alpha_{\mbox{\tiny left}}$, $\alpha'_{\mbox{\tiny left}}$,  $\alpha''_{\mbox{\tiny left}}$,
$\alpha_{\mbox{\tiny right}}$, $\alpha'_{\mbox{\tiny right}}$ and $\alpha''_{\mbox{\tiny right}}$,
as well as four real-valued coefficients  $c_{11}, c_{12}, c_{21}, c_{22}$ 
such that the functions
\begin{equation}
u(t) =  
\begin{cases}
 \frac{\cos\left(\alpha_{\mbox{\tiny left}}(t)\right)}{\sqrt{\alpha_{\mbox{\tiny left}}'(t)}}
& \ \ t \leq c
\\[1.3em]
c_{11} \frac{\cos\left(\alpha_{\mbox{\tiny right}}(t)\right)}{\sqrt{\alpha_{\mbox{\tiny right}}'(t)}}
+
c_{12}\frac{\sin\left(\alpha_{\mbox{\tiny right}}(t)\right)}{\sqrt{\alpha_{\mbox{\tiny right}}'(t)}}
& \ \ t > c
\end{cases}
\label{algorithm:ueven}
\end{equation}
and
\begin{equation}
v(t) =  
\begin{cases}
 \frac{\sin\left(\alpha_{\mbox{\tiny left}}(t)\right)}{\sqrt{\alpha_{\mbox{\tiny left}}'(t)}}
& \ \ t \leq c
\\[1.3em]
c_{21} \frac{\cos\left(\alpha_{\mbox{\tiny right}}(t)\right)}{\sqrt{\alpha_{\mbox{\tiny right}}'(t)}}
+
c_{22} \frac{\sin\left(\alpha_{\mbox{\tiny right}}(t)\right)}{\sqrt{\alpha_{\mbox{\tiny right}}'(t)}}
& \ \ t > c
\end{cases}
\label{algorithm:veven}
\end{equation}
form a basis in the space of solutions of (\ref{algorithm:ode}).


We begin the description of our algorithm by detailing a
fairly standard adaptive Chebyshev solver for ordinary differential equations
which it utilizes in Subsection~\ref{section:algorithm:odesolver}.  
We then describe the ``windowing'' procedure of \cite{BremerKummer},
which is also a component of the scheme of this paper, in Subsection~\ref{section:algorithm:windowing}.  
The method proper is described in Subsection~\ref{section:algorithm:main}.

\begin{subsection}{An adaptive solver for ordinary differential equations}
\label{section:algorithm:odesolver}

We now describe a fairly standard adaptive Chebyshev solver 
which is repeatedly used by the algorithm of this paper. It solves the problem
\begin{equation}
\left\{
\begin{aligned}
\bm{y}'(t) &= F(t,\bm{y}(t)), \ \ \ a_0 < t < b_0,\\
\bm{y}\left(t_0\right) &= \bm{v},
\end{aligned}
\right.
\label{algorithm:system}
\end{equation}
where $F:\mathbb{R}^{n+1} \to \mathbb{R}^n$ is smooth, $t_0 \in [a_0,b_0]$ and $\bm{v} \in \mathbb{R}^n$,
It takes as input the interval $(a_0,b_0)$, the point $t_0$, the vector $\bm{v}$,
and an external subroutine for evaluating the function $F$.  It uses the parameters
$\epsilon$ and $l$ which are supplied as inputs to the algorithm of this paper.

It outputs $n$ piecewise $l^{th}$ order Chebyshev expansions,
one for each of the components $y_i(t)$ of the solution $\bm{y}$ of (\ref{algorithm:system}).
Throughout, the solver maintains a list of ``accepted intervals'' of $[a_0,b_0]$.
An interval is accepted if the solution is deemed to be adequately represented by an $l^{th}$ 
order Chebyshev expansion on that interval.  

The solver operates in two phases.    In the first, it constructs piecewise Chebyshev expansions
representing the  solution on $[a_0,t_0]$.  During this phase,
a list of subintervals of $[a_0,t_0]$ to process is maintained.
Assuming $a_0 \neq t_0$, this list initially contains $[a_0,t_0]$.  Otherwise, it is 
empty.  
The following steps are repeated as long as the list of 
subintervals of $[a_0,t_0]$ to process is not empty:

\begin{enumerate}

\item
Extract the subinterval $[c_0,d_0]$ from the list of intervals
to process such that $d_0$ is as large as possible.

\item
Solve the terminal value problem
\begin{equation}
\left\{
\begin{aligned}
\bm{u}'(t) &= F(t,\bm{u}(t)), \ \ \ c_0< t < d_0,\\
\bm{u}(d_0) &= \bm{w}.
\end{aligned}
\right.
\label{algorithm:tvp2}
\end{equation}
If $d_0 = t_0$, then $\bm{w}=\bm{v}$.    Otherwise,  the value of the solution at the 
point $d_0$ has already been approximated, and we use that estimate
for $\bm{w}$.

If the problem is linear, a Chebyshev integral equation method (see, for instance, \cite{GreengardSolver})
is used to solve (\ref{algorithm:tvp2}).  Otherwise, 
the trapezoidal method (see, for instance, \cite{Ascher}) is first used to produce an initial
approximation $\mathbf{u_0}$ of the solution and then Newton's method is applied to refine it.
The linearized problems are solved using a Chebyshev integral equation method.

In any event, the result is a set of $l^{th}$ order Chebyshev expansions 
\begin{equation}
u_{i}(t) =  \sum_{j=0}^{l} \lambda_{ij}\ T_j\left(\frac{2}{d_0-c_0} t + \frac{d_0+c_0}{d_0-c_0}\right),
\ \ \ i=1,\ldots,n,
\label{algorithm:exps}
\end{equation}
approximating  the components $u_1,\ldots,u_n$ of the solution of (\ref{algorithm:tvp2}).

\item
Compute the quantities
\begin{equation}
\xi_i = \frac{\sqrt{\sum_{j=\frac{l}{2}+1}^l \lambda_{ij}^2}}{\sqrt{\sum_{j=0}^l \lambda_{ij}^2}}, \ \ \ i=1,\ldots,n,
\end{equation}
where the $\lambda_{ij}$ are the coefficients in the expansions (\ref{algorithm:exps}).
If any of the resulting values is  larger than $\epsilon$,
then we split the interval into two halves $\left[c_0,\frac{c_0+d_0}{2}\right]$ and 
$\left[\frac{c_0+d_0}{2},d_0\right]$ and place them on the list of subintervals 
of $[a_0,t_0]$ to process.  Otherwise, we place the interval
$[c_0,d_0]$ on the list of accepted subintervals.

\end{enumerate} 

In the second phase, the solver  constructs piecewise Chebyshev expansions
representing the  solution on $[t_0,b_0]$.  During this phase,
a list of subintervals of $[t_0,b_0]$ to process is maintained.
Assuming $b_0 \neq t_0$, this list initially contains $[t_0,b_0]$.  Otherwise, it is 
empty.    
The following steps are repeated as long as the list of 
subintervals of $[t_0,b_0]$ to process is not empty:

\begin{enumerate}

\item
Extract the subinterval $[c_0,d_0]$ from the list of intervals
to process such that $c_0$ is as small as possible.

\item
Solve the initial value problem
\begin{equation}
\left\{
\begin{aligned}
\bm{u}'(t) &= F(t,\bm{u}(t)), \ \ \ c_0< t < d_0,\\
\bm{u}(c_0) &= \bm{w}.
\end{aligned}
\right.
\label{algorithm:ivp2}
\end{equation}
If $c_0 = t_0$, then $\bm{w}=\bm{v}$.    Otherwise,  the value of the solution at the 
point $c_0$ has already been approximated, and we use that estimate
for $\bm{w}$.

If the problem is linear, a straightforward Chebyshev integral equation method is used.  Otherwise,
the trapezoidal method is first used to produce an initial
approximation $\mathbf{u_0}$ of the solution and then Newton's method is applied to refine it.
The linearized problems are solved using a straightforward Chebyshev integral equation method.

In any event, the result is a set of $l^{th}$ order Chebyshev expansions 
\begin{equation}
u_i(x)  \approx \sum_{j=0}^l \lambda_{ij}\ T_j\left(\frac{2}{d_0-c_0} t + \frac{c_0+d_0}{c_0-d_0}\right),\ \ \ i=1,\ldots,n,
\label{algorithm:exps3}
\end{equation}
approximating  the components $u_1,\ldots,u_n$ of the solution of (\ref{algorithm:tvp2}).

\item
Compute the quantities
\begin{equation}
\frac{\sqrt{\sum_{j=\frac{l}{2}+1}^l \lambda_{ij}^2}}{\sqrt{\sum_{j=0}^l \lambda_{ij}^2}}, \ \ \ i=1,\ldots,n,
\end{equation}
where the $\lambda_{ij}$ are the coefficients in the expansions (\ref{algorithm:exps3}).
If any of the resulting values is  larger than $\epsilon$,
then we split the interval into two halves $\left[c_0,\frac{c_0+d_0}{2}\right]$ and 
$\left[\frac{c_0+d_0}{2},d_0\right]$ and place them on the list of subintervals 
of $[t_0,b_0]$ to process.  Otherwise, we place the interval
$[c_0,d_0]$ on the list of accepted intervals.

\end{enumerate} 

At the conclusion of this procedure,  we have $l^{th}$ order piecewise Chebyshev expansions
representing each component of the solution.
\end{subsection}

\begin{subsection}{The ``windowing'' procedure of \cite{BremerKummer}}
\label{section:algorithm:windowing}

We now describe a variant of the ``windowing'' procedure of \cite{BremerKummer}, which is used as a step 
in the algorithm of this paper.   It takes as input an interval $[a_0,b_0]$
on which the coefficient $q$ in (\ref{algorithm:ode}) is nonnegative
and  makes use of the external subroutine for evaluating the coefficient $q$ 
which is specified as input to the algorithm of this paper. 
The derivative of $q$ is not needed by this procedure.
It  outputs the values of the first and second derivatives
of a nonoscillatory phase function for  (\ref{algorithm:ode}) on the interval
$[a_0,b_0]$ at the point $a_0$.  The algorithm can easily be modified to provide
the values of these functions at $b_0$ instead.

We first construct a ``windowed'' version $\tilde{q}$ of $q$ which closely
approximates $q$ on the leftmost quarter of the interval $[a_0,b_0]$ and
is approximately equal to 
\begin{equation}
 q\left(\frac{a_0+b_0}{2}\right)
\end{equation}
on the rightmost quarter of the interval $[a_0,b_0]$.  More explicitly,
we set
\begin{equation}
\nu = \sqrt{ q\left(\frac{a_0+b_0}{2}\right)}
\end{equation}
and define  $\widetilde{q}$ via the formula
\begin{equation}
\widetilde{q}(t) = \phi(t) \nu^2 + (1-\phi(t)) q(t),
\end{equation}
where $\phi$ is given by 
\begin{equation}
\phi(t) = \frac{
1+\mbox{erf}\left(\frac{12}{b_0-a_0} \left(t-\frac{a_0+b_0}{2}\right)\right)
}{2}.
\label{algorithm:errfun}
\end{equation}
The constant in (\ref{algorithm:errfun}) is chosen so that 
\begin{equation}
\left|\phi(a_0) \right|,\ \left|\phi(b_0) - 1\right| < \epsilon_0,
\end{equation}
where $\epsilon_0$ denotes IEEE double precision machine zero.  
We next solve the terminal value problem
\begin{equation}
\left\{
\begin{aligned}
\widetilde{q}(t) - (\widetilde{\alpha}'(t))^2 + \frac{3}{4} \left(\frac{\widetilde{\alpha}''(t)}
{\widetilde{\alpha}'(t)}\right)^2
- \frac{\widetilde{\alpha}'''(t)}{2\widetilde{\alpha}'(t)} = 0 \\
\widetilde{\alpha'}(b_0) = \nu\\
\widetilde{\alpha''}(b_0) = 0
\end{aligned}
\right.
\label{algorithm:windowkummer}
\end{equation}
using the solver described in Subsection~\ref{section:algorithm:odesolver}.
Although it outputs $l^{th}$ order piecewise Chebyshev expansions representing
 $\widetilde{\alpha}'$ and  $\widetilde{\alpha}''$, it is only the 
values of these functions at the point $a_0$ which concern us.  These
are the outputs of the windowing procedure, and they closely approximate
the values of  $\alpha'(a_0)$ and $\alpha''(a_0)$, where $\alpha$ is the 
the desired nonoscillatory phase function for the original equation
(an error estimate under mild conditions on $q$ is proven in \cite{BremerKummer}).

\end{subsection}

\begin{subsection}{A phase function method for differential equations with a turning point}
\label{section:algorithm:main}

Our algorithm operates differently depending on the order of the turning point $c$.
If $c$ is of odd order (i.e., $k$ is an odd integer), then we first apply the windowing
algorithm of Subsection~\ref{section:algorithm:windowing} on the interval $[c,b]$,
where $q$ is nonnegative and the solutions of (\ref{algorithm:ode}) 
oscillate.  By doing so, we obtain the values of the first two
derivatives of a  nonoscillatory phase function $\alpha$ for (\ref{algorithm:ode}) 
at the point $c$.  We next compute the value of $\alpha'''(c)$ using Kummer's equation:
\begin{equation}
\alpha'''(c) = 2\alpha'(c)q(c)-  2(\alpha'(c))^3  + \frac{3}{2} \frac{\left(\alpha''(c)\right)^2}{\alpha'(c)}.
\end{equation}
The first three derivatives of the function
\begin{equation}
w(t) = \frac{1}{\alpha'(t)},
\end{equation}
which satisfies Appell's equation
\begin{equation}
w'''(t) + 4q(t) w'(t)  +2 q'(t) w(t) = 0,
\label{algorithm:main:appell}
\end{equation}
at the point $c$ are given by
\begin{equation}
\begin{aligned}
w(c)   = \frac{1}{\alpha'(c)},\ \ \
w'(c)  = - \frac{\alpha''(c)}{\left(\alpha'(c)\right)^2} \ \ \ \mbox{and} \ \ \
w''(c) = 2 \frac{\left(\alpha''(c)\right)^2}{\left(\alpha'(c)\right)^3}
-\frac{\alpha'''(c)}{\left(\alpha'(c)\right)^2}.
\end{aligned}
\label{algorithm:main:wvals}
\end{equation}
We now apply the algorithm of Subsection~\ref{section:algorithm:odesolver}
find the  solution of (\ref{algorithm:main:appell}) which satisfies
(\ref{algorithm:main:wvals}).  Appell's equation includes
a term involving the derivative $q'(t)$ of the coefficient $q(t)$.
If the values of $q'(t)$ are not specified, then they are computed
using spectral differentiation.  In particular, on each interval
 considered by the solver of Subsection~\ref{section:algorithm:odesolver},
a Chebyshev expansion representing $q(t)$ is formed and used
to evaluate the derivatives of $q'(t)$.  
The solution of Appell's equation can become extremely large 
in the  nonoscillatory interval  $[a,c]$.  The adaptive solver
used to solve Appell's equation terminates its first 
phase early if the value of $w(t)$ on the interval under consideration
grows above $10^{300}$. This has the effect of truncating
the interval over which the phase function is given so that
the value of $\alpha'$ is close to the smallest representable
IEEE double precision number at the left endpoint.
The output of the solver
is three $l^{th}$ order Chebyshev expansions,
one for each of the functions $w(t)$, $w'(t)$ and $w''(t)$.  
At this stage, we construct $l^{th}$ order piecewise  order Chebyshev  expansions
representing each of the functions
\begin{equation}
\alpha'(t) = \frac{1}{w(t)}
\end{equation}
and
\begin{equation}
\alpha''(t) = -\frac{w'(t)}{(w(t))^2}
\end{equation}
over an  interval  of the form $[a',b]$ with $a \leq a' < b$.  An $l^{th}$ order piecewise Chebyshev expansion
representing the phase function $\alpha$ itself is constructed via spectral integration;
the particular choice of antiderivative is determined by the requirement that 
$\alpha(c) = 0$.

In the case of a turning point of even order, we first apply the windowing
algorithm of Subsection~\ref{section:algorithm:windowing} to the interval
$[a,c]$ in order to obtain the values of the first two derivatives of the phase
function $\alpha_{\mbox{\tiny left}}$ at the point $a$.  One this has been
done, Appell's equation is solved over the interval $[a,c]$
using the algorithm of  Subsection~\ref{section:algorithm:odesolver}
and  $l^{th}$ order piecewise Chebyshev expansions
of $\alpha_{\mbox{\tiny left}}'$ and $\alpha_{\mbox{\tiny left}}''$ are constructed
using the obtained solution.
An  $l^{th}$ order piecewise Chebyshev expansion of $\alpha_{\mbox{\tiny left}}$
is constructed via spectral integration; the value of $\alpha_{\mbox{\tiny left}}$
at $c$ is taken to be $0$.  This procedure is repeated on the interval $[c,b]$
in order to construct 
$l^{th}$ order piecewise Chebyshev expansions of $\alpha_{\mbox{\tiny right}}$
and its first two derivatives.  The value of $\alpha_{\mbox{\tiny right}}$ at
$c$ is taken to be $0$.  It is easy to see that the coefficients $c_{11}$ and  $c_{12}$
in (\ref{algorithm:ueven}) are given by the formulas
\begin{equation}
\begin{aligned}
c_{11} =  \sqrt{\frac{\alpha_{\mbox{\tiny right}}'(c)}{\alpha_{\mbox{\tiny left}}'(c)}} \ \ \ \mbox{and} \ \ \
c_{12} = - \frac{\alpha_{\mbox{\tiny right}}'(c) \alpha_{\mbox{\tiny left}}''(c) + \alpha_{\mbox{\tiny left}}'(c)\alpha_{\mbox{\tiny right}}''(c) }
{2\left(\alpha_{\mbox{\tiny left}}'(c)\right)^\frac{3}{2}\left(\alpha_{\mbox{\tiny right}}'(c)\right)^\frac{3}{2}},
\end{aligned}
\end{equation}
while the coefficients $c_{21}$ and $c_{22}$ in  (\ref{algorithm:veven}) are
\begin{equation}
\begin{aligned}
c_{21} = 0 \ \ \ \mbox{and} \ \ \ 
c_{22} =\sqrt{\frac{\alpha_{\mbox{\tiny right}}'(c)}{\alpha_{\mbox{\tiny left}}'(c)}}.
\end{aligned}
\end{equation}

\end{subsection}

\end{section}

\begin{figure}[!h]
\hfil
\includegraphics[width=.43\textwidth]{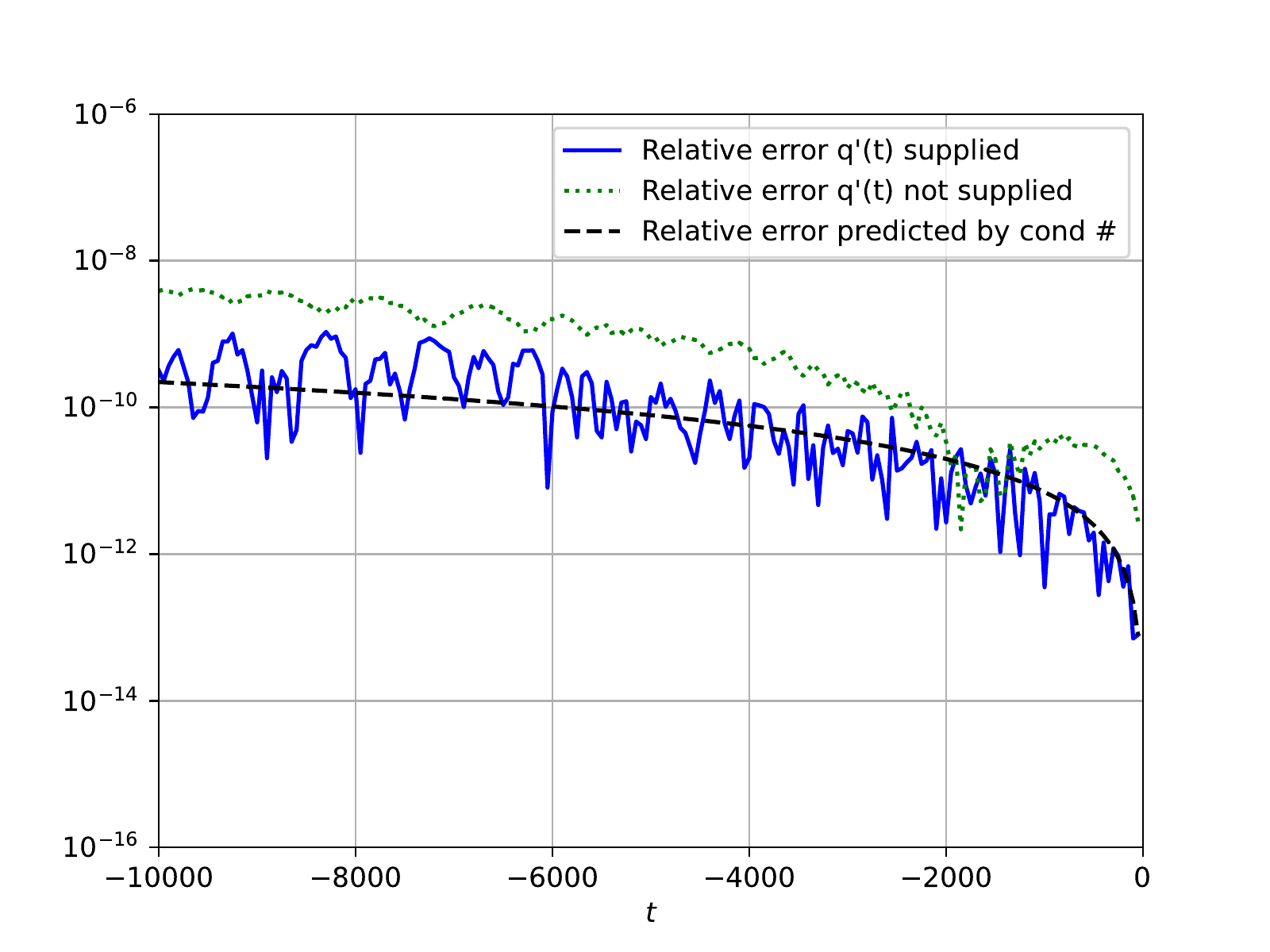}
\hfil
\includegraphics[width=.43\textwidth]{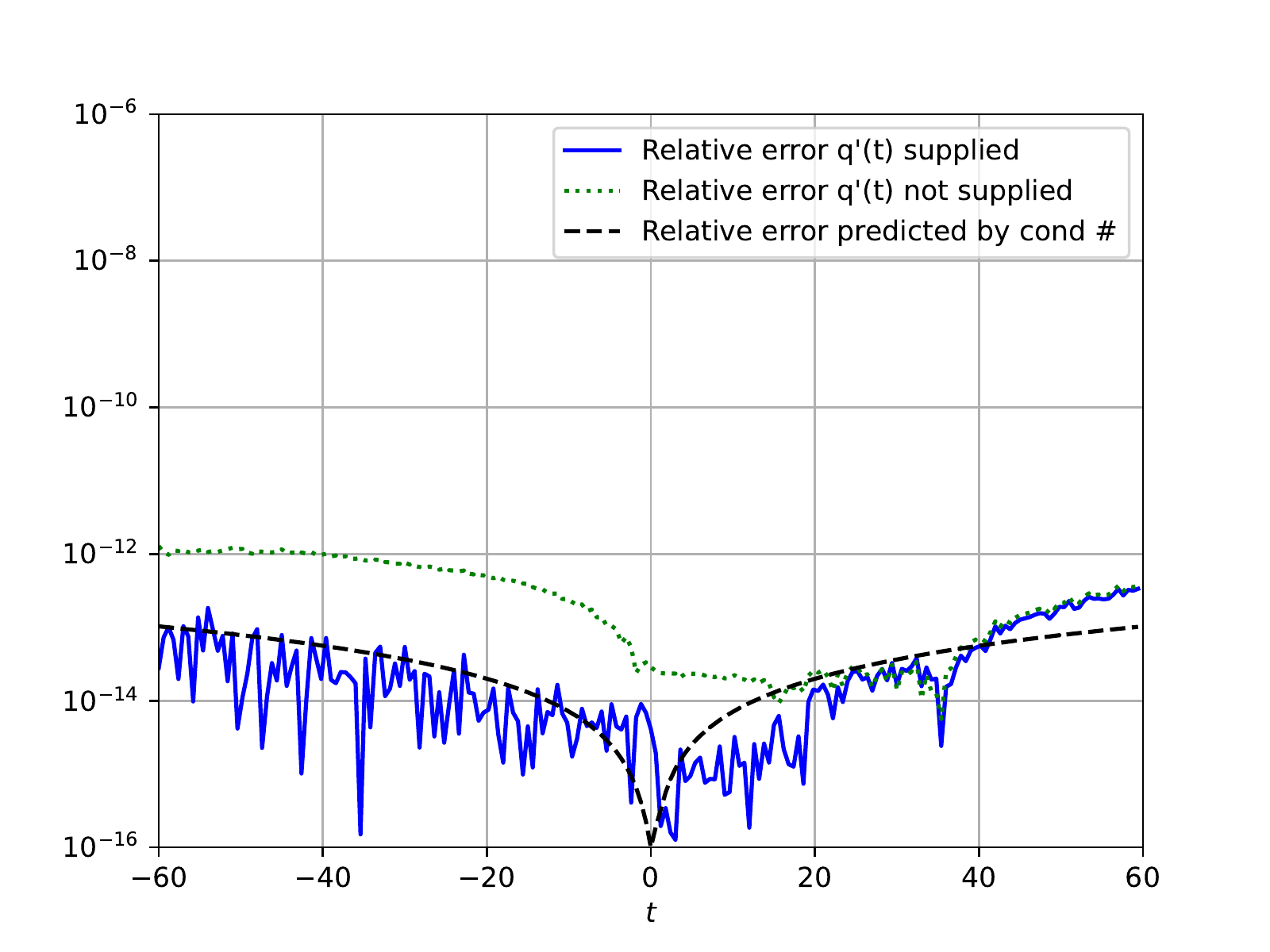}
\hfil

\hfil
\includegraphics[width=.43\textwidth]{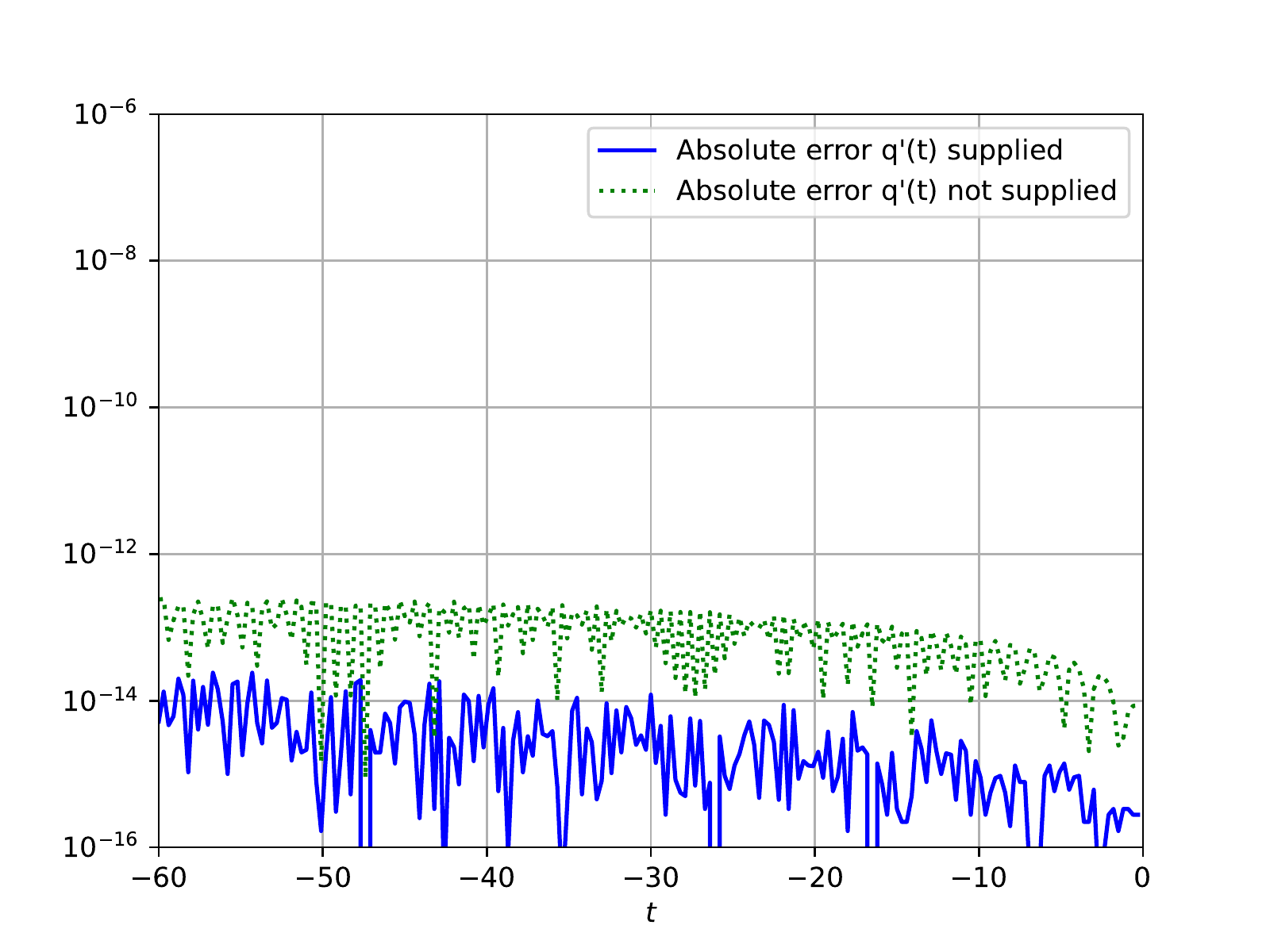}
\hfil
\includegraphics[width=.43\textwidth]{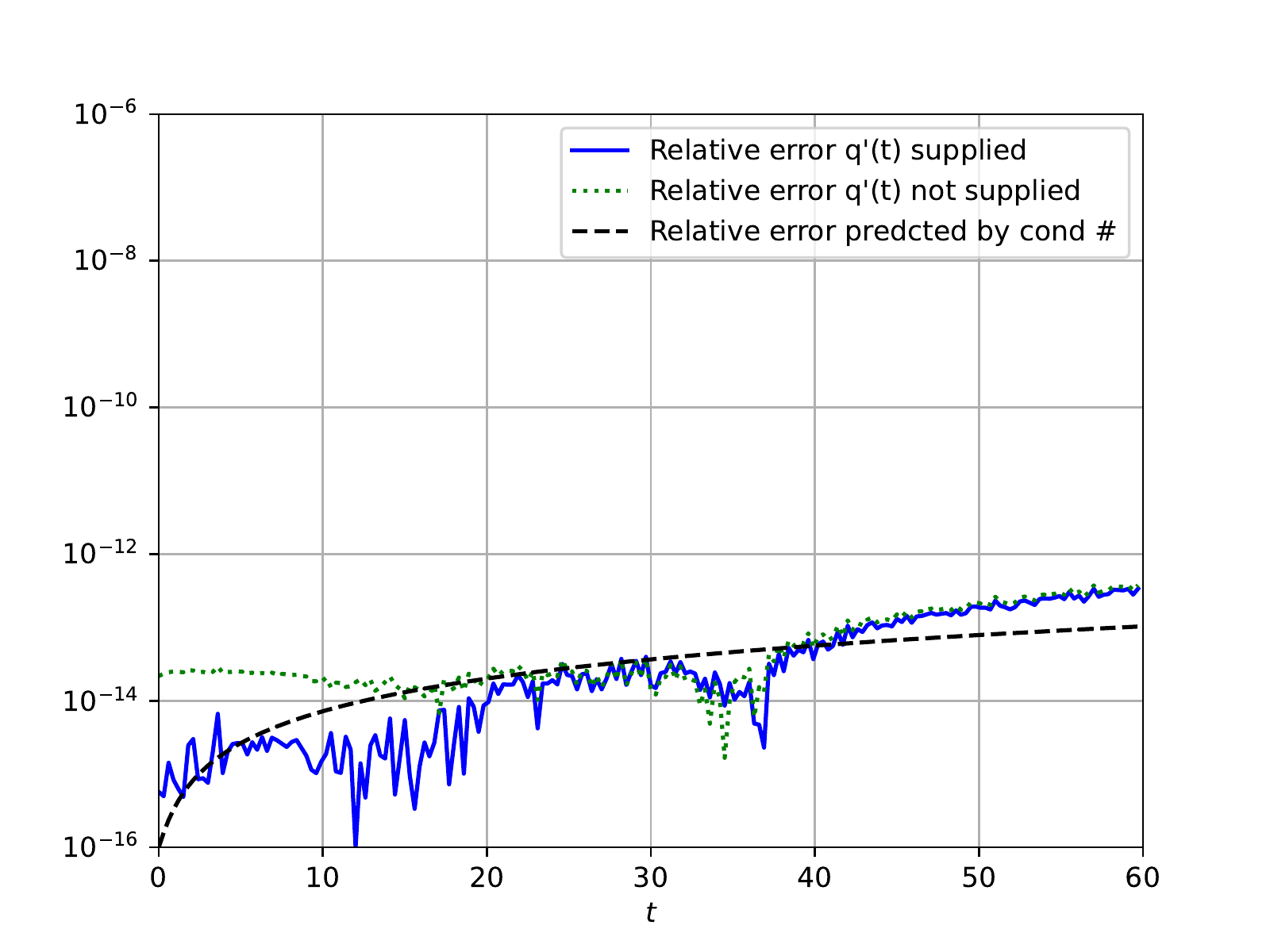}
\hfil

\caption{The results of the experiments discussed in Subsection~\ref{section:experiments:airy},
which concerned the Airy functions.   The plots in the first row compare
the relative errors  incurred when evaluating $f^{\mbox{\tiny airy}}(t) = \mbox{Ai}(t) + i \mbox{Bi}(t)$ 
using the phase functions $\alpha^{\mbox{\tiny airy}}$ 
and $\widetilde{\alpha}^{\mbox{\tiny airy}}$  with the relative error predicted by its condition number of 
evaluation.   The plot in the upper left is given over the interval $(-10000,0)$,
and that on the upper right is given over the interval $(-60,60)$.
The plot on the lower left gives the absolute error in the calculated
values of $\mbox{Ai}(t)$  in the interval $(-60,0)$ as a function of $t$,
while that on the lower right gives the relative error in the calculated
value of $\mbox{Ai}(t)$ in the interval $(0,60)$ as a function of $t$  
and compares it with the relative error predicted by the condition number
of evaluation of $\mbox{Ai}(t)$.
}
\label{figure:airyplots1}
\end{figure}

\begin{section}{Numerical Experiments}

In this section, we present the results of numerical experiments
which were conducted to illustrate the properties of the
algorithm of this paper.  The code for these experiments
was written in Fortran and compiled with version 12.10
of the GNU Fortran compiler.     The experiments were performed on 
a desktop computer equipped  with an AMD Ryzen 3900X processor and 32GB 
of memory.  No attempt  was made to parallelize our code.    

In the first four experiments described here, we considered problems 
with known solutions and so were able to measure the error in our
obtained solutions by direct comparison with them.
Moreover, in these experiments,
we were able to  compare the relative accuracy of our method with that predicted by 
the condition number of evaluation of the solution.  
By the relative accuracy predicted by the condition number
of evaluation of a function $f$, we mean the quantity
$\kappa_f \times \epsilon_0$, where
\begin{equation}
\kappa_f(t)  = \left|\frac
{f'(t)}
{f(t)}
t 
\right|
\end{equation}
is the condition number of evaluation of $f$ and 
and $\epsilon_0 \approx 2.220446049250313 \times 10^{-16}$ is machine zero
for IEEE double precision arithmetic.
The product of machine zero and the condition number of a function
is a reasonable heuristic for  the relative accuracy which can be 
expected when evaluating $f$ using IEEE double precision arithmetic.
A thorough discussion of the notion of ``condition number of evaluation
of a function'' can be found, for instance, in \cite{higham}.

Explicit formulas are not available for the solutions of the problems
discussed in  Subsections~\ref{section:experiments:bumps}, \ref{section:experiments:two} and
\ref{section:experiments:many}.  Consequently, we ran the conventional solver
described in Subsection~\ref{section:algorithm:odesolver}
using quadruple precision arithmetic (i.e., using Fortran REAL*16 numbers) to construct
reference solutions in these experiments.
We used extended precision arithmetic because we found in the course
of conducting the experiments for this article that,
in most cases, our  method obtains higher accuracy than the 
conventional solver described in Subsection~\ref{section:algorithm:odesolver}.


\label{section:experiments}

\begin{subsection}{Airy functions}
\label{section:experiments:airy}

The experiments of this subsection concern the Airy functions $\mbox{Ai}$ and $\mbox{Bi}$,
which  are standard solutions of the differential equation
\begin{equation}
y''(t) - t y(t) = 0, \ \ -\infty < t <\infty.
\label{experiments:airy:ode}
\end{equation}
Their definitions can be found in many sources (for instance, \cite{Olver}).
Obviously, (\ref{experiments:airy:ode}) has a turning point of order $1$ at $0$.
The Airy functions oscillate on the interval $(-\infty,0)$, while $\mbox{Ai}$ is decreasing
on $(0,\infty)$ and $\mbox{Bi}$ is increasing on  $(0,\infty)$.

We first used the algorithm of this
section to construct two  phase functions
$\alpha^{\mbox{\tiny airy}}$ and $\widetilde{\alpha}^{\mbox{\tiny airy}}$ for  (\ref{experiments:airy:ode})
using the algorithm of this paper.  In the case of 
$\alpha^{\mbox{\tiny airy}}$, the values of the derivative of the coefficient $q(t)=-t$ were
supplied as input and in the case of $\widetilde{\alpha}^{\mbox{\tiny airy}}$, the derivative
of $q(t)$ was calculated via spectral differentiation.
The phase functions were constructed on the interval $(-10000,b)$, where $b = 64.43359375$.
The right endpoint was chosen by our algorithm so that the value of the derivative of the phase function
$\alpha^{\mbox{\tiny airy}}$ was extremely small there.  In fact:
\begin{equation}
\begin{aligned}
\frac{d \alpha^{\mbox{\tiny airy}}}{dt} (b) &\approx 2.5585823472966988\times 10^{-299},\\
\mbox{Ai}(b) &\approx 1.7776196565817 \times 10^{-151} \ \ \mbox{and} \\
\mbox{Bi}(b) &\approx 1.1153864494678 \times 10^{149}.
\end{aligned}
\end{equation}

In a first experiment, we measured the relative error incurred when evaluating
the function $f^{\mbox{\tiny airy}}(t) = \mbox{Ai}(t) + i \mbox{Bi}(t)$
using each of the two phase functions for (\ref{experiments:airy:ode}).
We choose to consider $f^{\mbox{\tiny airy}}(t)$ rather than the Airy functions $\mbox{Ai}$
and $\mbox{Bi}$ individually because it does not vanish and 
its absolute value  is nonoscillatory on the real axis (see Figure~\ref{experiments:airy:fplot}
for a plot of the absolute value of this function).
It took approximately $2$ milliseconds to construct each of the phase functions.

\begin{figure}[h!]
\begin{center}
\includegraphics[width=.50\textwidth]{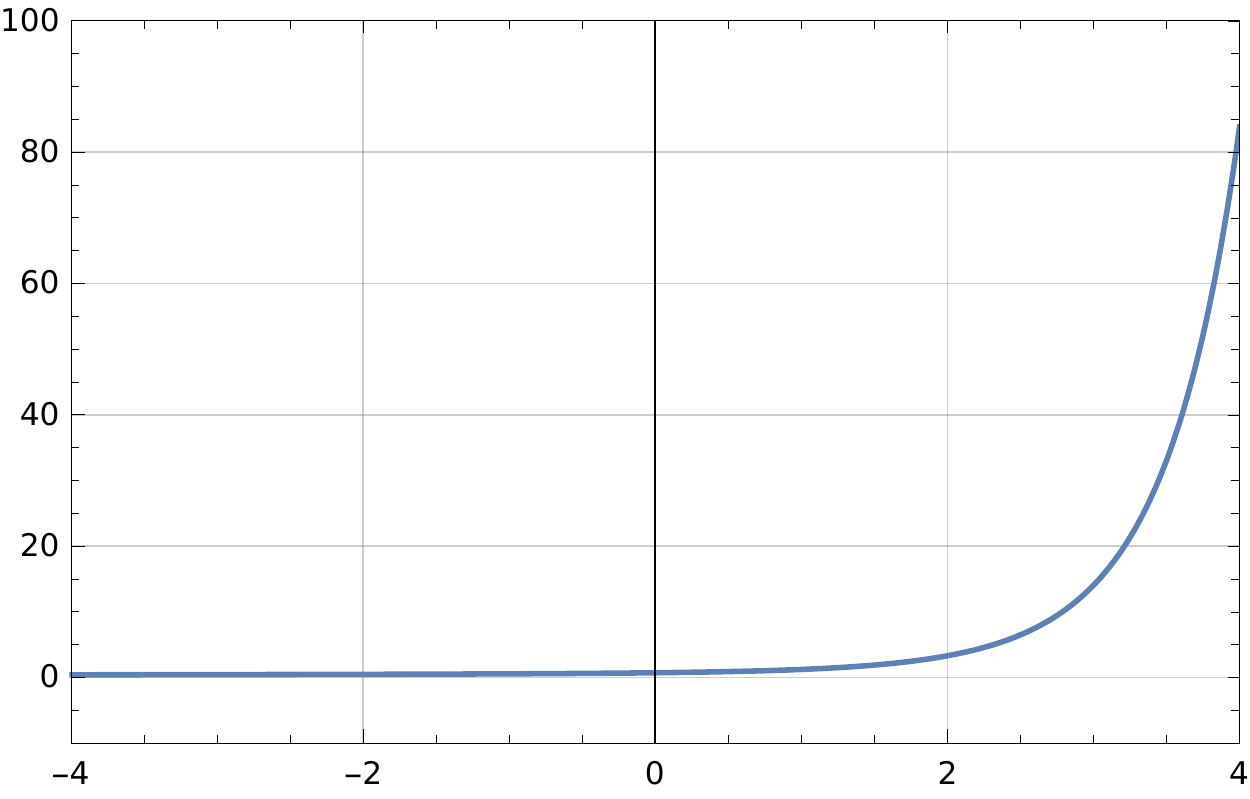}
\end{center}
\caption{A plot of  $\left|f^{\mbox{\tiny airy}}(t)\right|$,
where $f^{\mbox{\tiny  airy}}(t) = \mbox{Ai}(t) + i\mbox{Bi}(t)$.
In contrast to $f^{\mbox{\tiny  airy}}$, the absolute values of most solutions of Airy's differential equation 
are oscillatory on the negative real axis.
}
\label{experiments:airy:fplot}
\end{figure}

In our first experiment, we used each of the two phase functions to evaluate $f^{\mbox{\tiny airy}}$ at
$200$ equispaced points on the interval $(-10,000,0)$, where the Airy functions
are oscillatory.  The plot on the left-hand side of the top row of 
Figure~\ref{figure:airyplots1} gives the results.
There, the relative errors in the values of  $f^{\mbox{\tiny airy}}$
calculated using each of the two phase functions are plotted, as are the relative
errors  predicted by the condition number
of evaluation of this function.   

In our next experiment, we used the each of the two phase functions
to evaluate $f^{\mbox{\tiny airy}}$ at $200$ equispaced points on
the interval $(-60,60)$.  The plot on the right-hand side of the top
row of Figure~\ref{figure:airyplots1} compares the relative errors
which were incurred when doing so with those predicted by the condition number of 
evaluation of  $f^{\mbox{\tiny airy}}$ .

In a third experiment related to the Airy functions, we used the phase function for (\ref{experiments:airy:ode})
to evaluate $\mbox{Ai}$ at $200$ equispaced
points in the interval $(-60,0)$ using each of the two phase functions
for (\ref{experiments:airy:ode}).  Since $\mbox{Ai}$ is oscillatory on this
interval and has  many zeros there, it is not sensible to measure the relative errors in 
these calculations. Instead, we measured the absolute errors incurred when
$\mbox{Ai}$ was evaluated using
$\alpha^{\mbox{\tiny airy}}$ and $\widetilde{\alpha}^{\mbox{\tiny airy}}$.
The plot on the left-hand side of the second
row of  Figure~\ref{figure:airyplots1} gives the results.

In  a last experiment concerning the Airy functions, we used the two phase
functions for (\ref{experiments:airy:ode}) to evaluate $\mbox{Ai}$
at $200$ equispaced points on the  interval $(0,60)$. 
The plot on the right-hand side of the
 second row of Figure~\ref{figure:airyplots1} gives the results.
There, the relative errors in the values of  $\mbox{Ai}$
calculated using each of the two phase functions are plotted, as are the relative
errors  predicted by the condition number
of evaluation of $\mbox{Ai}$.

\end{subsection}


\begin{subsection}{Bessel functions}
\label{section:experiments:bessel}

In this set of experiments, we constructed phase functions for the normal form
\begin{equation}
y''(t) + \left(1 + \frac{\frac{1}{4}-\nu^2}{t^2}\right)y(t) = 0, \ \  0 <  t < \infty,
\label{experiments:bessel:ode}
\end{equation}
of Bessel's differential equation.   When $\nu > \frac{1}{2}$,
\begin{equation}
c = \frac{1}{2} \sqrt{4\nu-1}
\end{equation}
is a turning point of order $1$ for  (\ref{experiments:bessel:ode}).
The coefficient in (\ref{experiments:bessel:ode}) is negative on the interval $(0,c)$ and  positive on $(c,\infty)$.
Standard solutions include the Bessel functions of the first and second kinds $J_\nu$
and $Y_\nu$ (definitions of which were given in Subsection~\ref{section:turning:bessel}).


In the first experiment, 
we sampled  $m=200$ equispaced points $x_1,x_2,\ldots,x_m$ in the interval $[0,6]$
and constructed phase functions $\alpha_\nu^{\mbox{\tiny bes}}$ 
and $\widetilde{\alpha}_\nu^{\mbox{\tiny bes}}$ for each 
$\nu = 10^{x_1}, 10^{x_2}, \ldots,10^{x_m}$.  
In the case of 
$\alpha^{\mbox{\tiny bes}}_\nu$, the derivative of the coefficient was
supplied as input and in the case of $\widetilde{\alpha}^{\mbox{\tiny bes}}_\nu$, the derivative
was calculated via spectral differentiation.
We then used each of these phase functions to evaluate the function
$f_\nu^{\mbox{\tiny bes}}(t) = J_\nu(t) + i Y_\nu(t)$ at $5,000$ equispaced points
on the interval $(0,100\nu)$.  We considered the function
$f_\nu^{\mbox{\tiny bes}}$ rather than the Bessel functions individually
because its absolute value is nonoscillatory; indeed,
it  was shown in Section~\ref{section:turning:bessel} that 
$\left|f_\nu^{\mbox{\tiny bes}}\right|^2$ is completely monotone on $(0,\infty)$.
The plot on the left-hand side of Figure~\ref{figure:besselplots1}
gives  the maximum observed relative errors as functions of  $\nu$, as well
as the maximum relative error predicted by the condition number of
evaluation of $f^{\mbox{\tiny bes}}_\nu$.  The plot on the right-hand side
of Figure~\ref{figure:besselplots1} gives the time require to construct
the phase function $\alpha_\nu^{\mbox{\tiny bes}}$ as a function of $\nu$.

In a second set of experiments, for each $\nu=1,10,10^2,\ldots,10^5$, we 
measured the relative errors incurred
when the phase functions $\alpha_\nu^{\mbox{\tiny bes}}$ 
and $\widetilde{\alpha}_\nu^{\mbox{\tiny bes}}$ were used to evaluate
$f^{\mbox{\tiny bes}}_\nu(t)$ at $200$ equispaced spaced points in the interval $(0,100 \nu)$.
Figure~\ref{figure:besselplots2} gives the results.  Each plot there 
corresponds to a different value of $\nu$, and gives the relative errors
incurred when the phase functions were used to
evaluate  $f^{\mbox{\tiny bes}}_\nu(t)$, 
as well as the relative error predicted by its condition number of evaluation,
as functions of $t$.

\end{subsection}

\begin{subsection}{Associated Legendre functions}
\label{section:experiments:alf}

The experiments described in this subsection concern the associated Legendre
differential equation
\begin{equation}
(1-t^2) y''(t) - t y(t) + \left(\nu(\nu+1) - \frac{\mu^2}{1-t^2}\right)y(t) = 0, \ \ -1 <t <1.
\label{experiments:alf:ode}
\end{equation}
Standard solutions of (\ref{experiments:alf:ode}) defined on the interval
$(-1,1)$ include the Ferrer's functions of the first and second kinds
$P_\nu^{\mu}$ and $Q_\nu^{\mu}$ (see, for instance, \cite{HTFI} or \cite{Olver} for definitions).
The Ferrer's functions of negative orders are generally better behaved than those
of positive orders --- for instance, when $\mu \geq 0$ is not an integer,
$P_\nu^{\mu}(t)$ is singular at $t=1$ whereas $P_\nu^{-\mu}(t)$ is not ---
and they coincide in the important case of integer values of $\nu$ and $\mu$.   
When $\nu \geq \left|\mu\right|$ (this is a standard requirement), the normalizations
\begin{equation}
\begin{aligned}
\widetilde{P}_\nu^{\mu}(t) &= \sqrt{\left(\nu+\frac{1}{2}\right)\frac{\Gamma\left(\nu+\mu+1\right)}{\Gamma\left(\nu-\mu+1\right)}}
\ P_\nu^{\mu}(t),\ \ \mbox{and} \\[1.3em]
\widetilde{Q}_\nu^{\mu}(t) &= \sqrt{\left(\nu+\frac{1}{2}\right)\frac{\Gamma\left(\nu+\mu+1\right)}{\Gamma\left(\nu-\mu+1\right)}}
\ Q_\nu^{\mu}(t)
\end{aligned}
\end{equation}
are well-defined, and we prefer them to the standard Ferrer's functions because
the latter can take on extremely large and extremely small values,
even when $\nu$ and $\mu$ are of modest sizes.    We note that 
the $L^2(-1,1)$ norm of $\widetilde{P}_\nu^{-\mu}$ is $1$ when both $\nu$ and $\mu$ are integers.

We claim that the phase function generated by the pair 
\begin{equation}
\widetilde{P}_\nu^{-\mu}(t), \ \ \frac{2}{\pi}\widetilde{Q}_\nu^{-\mu}(t)
\end{equation}
is nonoscillatory in a rather strong sense.  To see this, we first observe that
\begin{equation}
\begin{aligned}
\widetilde{P}_\nu^{-\mu}(t) + i \frac{2}{\pi}\widetilde{Q}_\nu^{-\mu}(t)
=
&\frac{2}{\pi}\sqrt{\frac{\nu+\frac{1}{2}}{\Gamma\left(\nu+\mu+1\right)\Gamma\left(\nu-\mu+1\right)}}
\exp\left(i\frac{\pi}{2} \left(\mu-\nu\right)\right)
\left(\frac{1}{1-t^2}\right)^{\frac{\nu+1}{2}}\\
&\int_0^\infty \exp\left(ix \frac{t}{\sqrt{1-t^2}}\right) K_\mu(x) x^\nu\ dx,
\end{aligned}
\label{phase:alf:int}
\end{equation}
where $K_\nu$ is the modified Bessel function of the third kind of order $\mu$.
This formula appears in a slightly different form in Section~7.8 of \cite{HTFII},
and it can be verified quite easily --- for instance, by showing that the expression
appearing on the right-hand side satisfies Appell's differential equation
and then verifying that the functions appearing on the left- and right-hand
sides and their first two derivatives agree at $0$.
Next, we make the change of variables
\begin{equation}
t = \frac{p}{\sqrt{1+p^2}}
\end{equation}
to obtain
\begin{equation}
\begin{aligned}
\widetilde{P}_\nu^{-\mu}\left(\frac{p}{\sqrt{1+p^2}}\right) + i \frac{2}{\pi}\widetilde{Q}_\nu^{-\mu}
\left(\frac{p}{\sqrt{1+p^2}}\right)=
&\frac{2}{\pi}\sqrt{\frac{\nu+\frac{1}{2}}{\Gamma\left(\nu+\mu+1\right)\Gamma\left(\nu-\mu+1\right)}}
\exp\left(i\frac{\pi}{2} \left(\mu-\nu\right)\right)\\
&\left(1+p^2\right)^{\frac{\nu+1}{2}}
\int_0^\infty \exp\left(itp\right) K_\mu(t) t^\nu\ dt.
\end{aligned}
\label{phase:alf:int2}
\end{equation}
We now observe that
\begin{equation}
\begin{aligned}
&\left(\widetilde{P}_\nu^{-\mu}\left(\frac{p}{\sqrt{1+p^2}}\right)\right)^2
+
\left(\frac{2}{\pi}\widetilde{Q}_\nu^{-\mu}\left(\frac{p}{\sqrt{1+p^2}}\right)\right)^2
\\ 
&\ \ =
\frac{4}{\pi^2}\ \frac{(1+p^2)^{\nu+1}\left(\nu+\frac{1}{2}\right)}{\Gamma\left(\nu-\mu+1\right) \Gamma\left(\nu+\mu+1\right) } 
\int_{-\infty}^\infty \exp(itp) G(t)\ dt,
\end{aligned}
\label{phase:alf:mod}
\end{equation}
where
\begin{equation}
G(t) =   \int_{\max(0,-t)}^\infty K_\mu(t+s) (t+s)^\nu K_\mu(s) s^\nu\ ds.
\label{phase:alf:convo}
\end{equation}
Since
\begin{equation}
K_\mu(z) \sim \sqrt{\frac{\pi}{2z}}\exp(-z)\ \ \mbox{as} \ \ z \to \infty,
\end{equation}
the function $K_\mu(t) t^\nu$ behaves similarly to a bump function centered
at the point $\nu-1/2$ when $\nu$ is of large magnitude.  It follows that the function $G$ defined in (\ref{phase:alf:convo})
resembles a bump function centered around the point $0$, again assuming $\nu$ is of large magnitude.
In particular, the function (\ref{phase:alf:mod}) is nonoscillatory in the sense that its
Fourier transform resembles a rapidly decaying bump function centered at $0$.
Figure~\ref{phase:figure:1} contains plots of the functions $1/\Gamma(\nu) K_\mu(t) t^\nu$
and $1/\Gamma(\nu)^2 G(t)$ when $\nu=100$ and $\mu=20$.    The scaling by the reciprocal of $\Gamma(\nu)$ was introduced
because the magnitudes of these functions are quite large otherwise.

\begin{figure}[!h]
\hfil
\includegraphics[width=.45\textwidth]{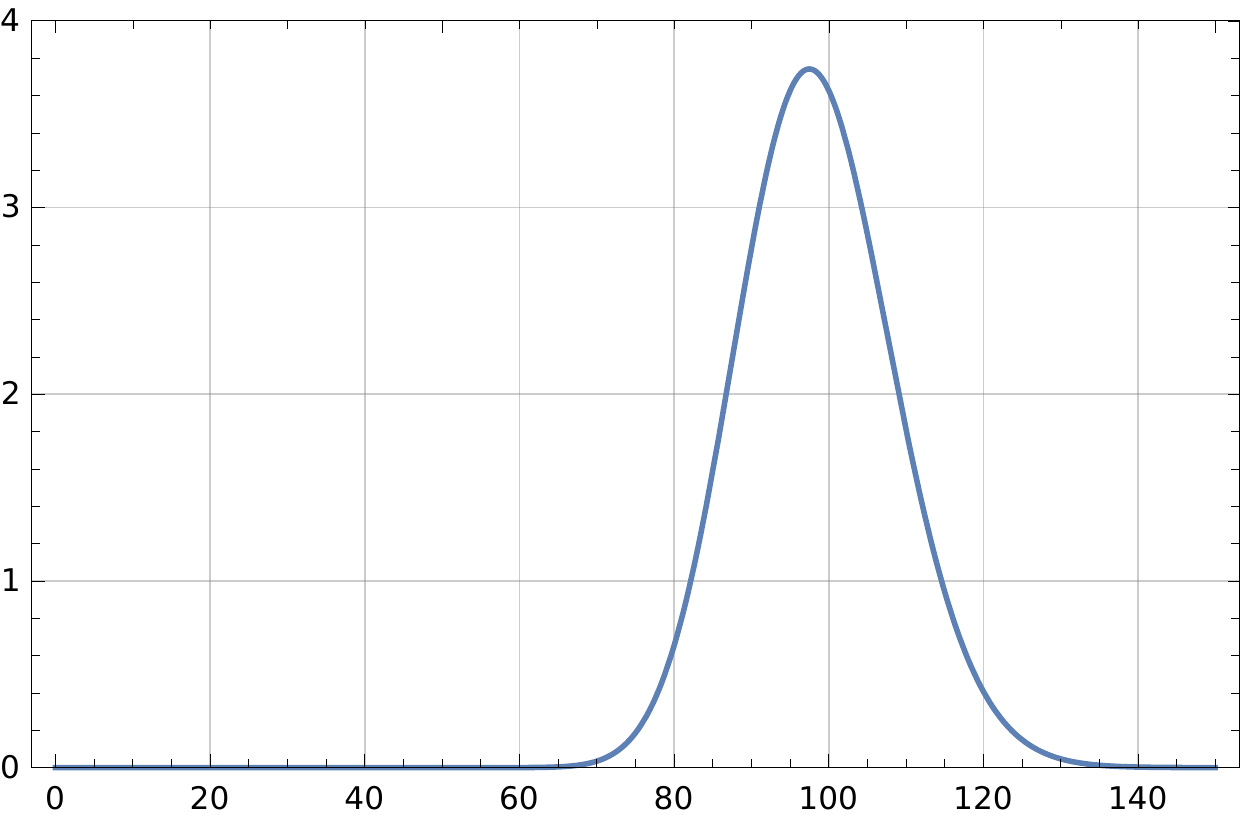}
\hfil
\includegraphics[width=.45\textwidth]{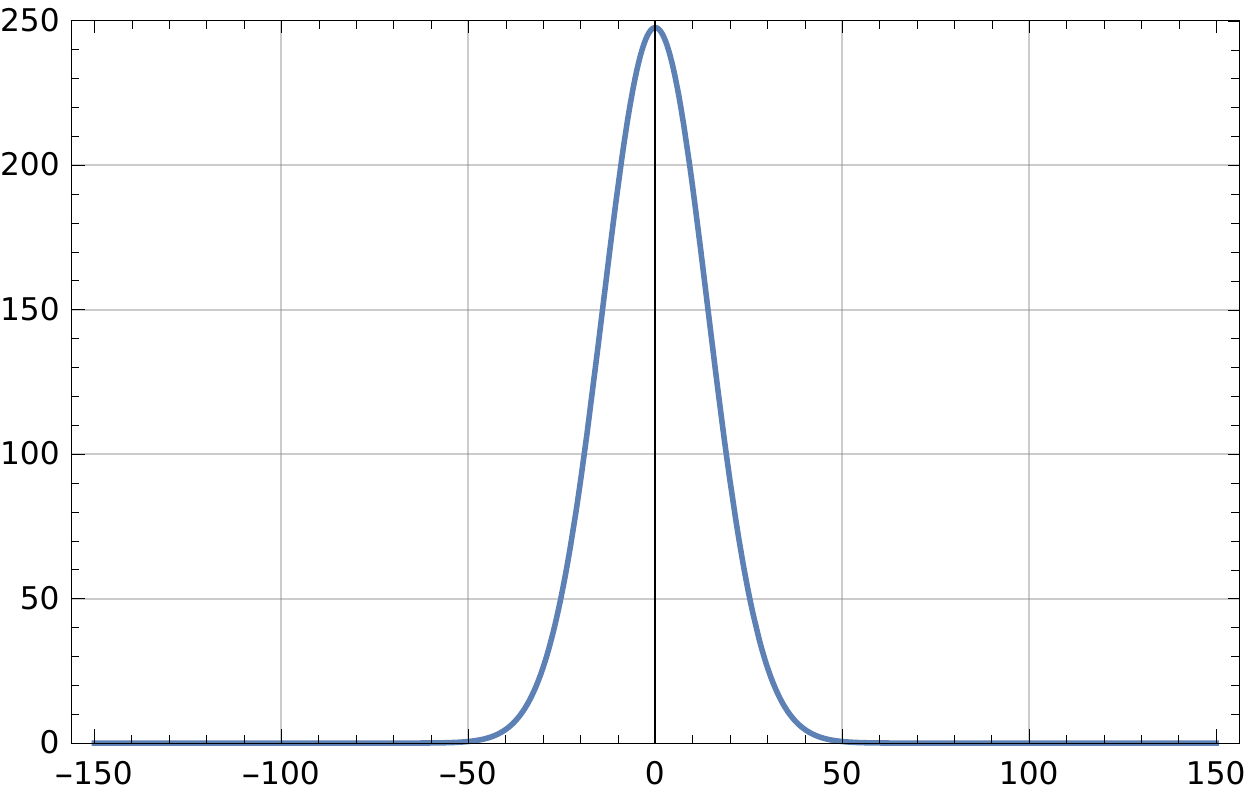}
\hfil
\caption{On the left is a plot of the function $1/\Gamma(\nu) K_\mu(t) t^\nu$ when $\nu=100$ and $\mu=20$.
On the right is a plot of $1/\Gamma(\nu)^2 G(t)$, where  $G$ is the function defined  in (\ref{phase:alf:convo}),
when $\nu=100$ and $\mu=20$.}
\label{phase:figure:1}
\end{figure}

When performing numerical computations involving the associated Legendre
functions, it is convenient to introduce a change of variables which eliminates the singular points of the
equation  (\ref{experiments:alf:ode}) at $\pm 1$.  If $y$ solves (\ref{experiments:alf:ode}), then
\begin{equation}
z(w) = y(\tanh(w))
\end{equation}
satisfies
\begin{equation}
z''(w) + \left(\nu (\nu+1) \sech^2 (w) - \mu^2 \right) z(w) = 0, \ \ -\infty < w < \infty.
\label{experiments:alf:ode2}
\end{equation}
We use $\alpha_{\nu,\mu}^{\mbox{\tiny alf}}$ to denote the phase function for (\ref{experiments:alf:ode2})
constructed via the algorithm of this paper with values of the derivative of the coefficient
supplied,  and 
$\widetilde{\alpha}_{\nu,\mu}^{\mbox{\tiny alf}}$ to denote the phase function for (\ref{experiments:alf:ode2})
constructed via the algorithm of this paper with derivatives of the coefficient
calculated using spectral differentiation.

In each experiment of this subsection, we first fixed a value of $\mu$
and  sampled $m=200$ equispaced points $x_1,\ldots,x_m$ in the interval $[0,6]$.
Then, for each $\lambda=10^{x_1},10^{x_2},\ldots,10^{x_m}$,
we constructed phase functions $\alpha_{|\mu|+\lambda,\mu}^{\mbox{\tiny alf}}$ and 
$\widetilde{\alpha}_{|\mu|+\lambda,\mu}^{\mbox{\tiny alf}}$
using the algorithm of this paper.  
In the case of  $\alpha_{|\mu|+\lambda,\mu}^{\mbox{\tiny alf}}$, the values of the derivatives
of the coefficient were supplied to the algorithm, while  spectral
differentiation was used to evaluate $q'(t)$ in the case of 
$\widetilde{\alpha}_{|\mu|+\lambda,\mu}^{\mbox{\tiny alf}}$.
To be clear, the degree $\nu$ of the associated Legendre function was taken to be $\nu = |\mu|+\lambda$ (it
is necessary for $\nu \geq |\mu|$).  The phase functions were given 
on intervals of the form $[0,b]$ with $b$ chosen by our algorithm
so that  the value of the derivative of $\alpha_{|\mu|+\lambda,\mu}^{\mbox{\tiny alf}}$ there was close
to the smallest representable IEEE double precision number.
For each chosen value of $\lambda$, we used both of the corresponding phase functions
to evaluate
\begin{equation}
f_{|\mu|+\lambda,\mu}^{\mbox{\tiny alf}}(w) = 
\widetilde{P}_{|\mu|+\lambda}^\mu(\tanh(w)) + 
\frac{2}{\pi} i \widetilde{Q}_{|\mu|+\lambda}^\mu(\tanh(w))
\label{experiments:alf:f}
\end{equation}
at $5,000$ equispaced points on the interval  $[0,b]$.  Figure~\ref{figure:alfplots1} reports the results.
Each row in that figure corresponds to one experiment.  The plot on the left-hand side
gives the maximum relative errors observed while evaluating (\ref{experiments:alf:f}) 
using each of the phase functions  $\alpha_{|\mu|+\lambda,\mu}^{\mbox{\tiny alf}}$ and 
$\widetilde{\alpha}_{|\mu|+\lambda,\mu}^{\mbox{\tiny alf}}$, as well as the maximum relative error
predicted by the condition number of evaluation of $f_{|\mu|+\lambda,\mu}^{\mbox{\tiny alf}}$,
as functions of $\lambda$.  The plot on the right gives the time (in milliseconds)
required to construct $\alpha_{|\mu|+\lambda,\mu}^{\mbox{\tiny alf}}$ as a function of $\lambda$.

\begin{remark}
 Formula~(\ref{phase:alf:int2}) implies that, after a suitable change of variables,
the function $\widetilde{P}_\nu^{-\mu}(t) + i \frac{2}{\pi}  \widetilde{Q}_\nu^{-\mu}(t)$
is an element of one of the Hardy spaces of functions analytic
in the upper half of the complex plane (see, for instance, \cite{Koosis} for a careful discussion
of such spaces).     
Likewise, Formula~(\ref{turning:beslaplace}) shows that $J_\nu(z)+ i Y_\nu(z)$
is an element of one of these Hardy spaces as well.
It is not a coincidence that these
two functions give rise to slowly varying phase functions.
 A relatively straightforward generalization
of the discussion in Subsection~\ref{section:experiments:alf} shows that when a second
order differential equation with slowly varying coefficients admits a solution
in one of the classical Hardy spaces of functions analytic in the upper half
of the complex plane (perhaps after a change of variables), 
it will have a slowly varying phase function.  
This argument applies to Bessel's differential equation, the spheroidal wave equation, 
and many other second order linear ordinary differential equations.
\end{remark}

\end{subsection}

\begin{subsection}{Equations with turning points of higher orders}
\label{section:experiments:high}

The experiments described in this subsection concerned the
differential equation
\begin{equation}
y''(t) + t^k y(t) =0, \ \ -10 < t < 10,
\label{experiments:high:ode}
\end{equation}
with $k$ a positive integer greater than $1$.    
In the event that $k$ is even, the functions $u_k^{\mbox{\tiny even}}$ and
$v_k^{\mbox{\tiny even}}$ defined in (\ref{turning:ueven})
and (\ref{turning:veven}) form a basis in the space of solutions
of (\ref{experiments:high:ode}).  It can be verified
using a argument similar to that of Subsection~\ref{section:turning:basis}
that
\begin{equation}
u_k^{\mbox{\tiny odd}}(t) =
 \begin{cases}
\sqrt{\frac{\pi t}{2+k}} J_{\frac{1}{2k}} \left(t^{1+\frac{k}{2}}{1+\frac{k}{2}}\right) & \mbox{if} \ \ t > 0\\
\sqrt{-\frac{-\pi t}{2+k}} I_{\frac{1}{2k}} \left((-t)^{1+\frac{k}{2}}{1+\frac{k}{2}}\right) & \mbox{if} \ \ t < 0\\
\end{cases}
\end{equation}
and
\begin{equation}
v_k^{\mbox{\tiny odd}}(t) =
 \begin{cases}
\sqrt{\frac{\pi t}{2+k}} Y_{\frac{1}{2k}} \left(\frac{t^{1+\frac{k}{2}}}{1+\frac{k}{2}}\right) & \mbox{if} \ \ t > 0\\
d_1\sqrt{-\frac{-\pi t}{2+k}} I_{\frac{1}{2k}} \left(\frac{(-t)^{1+\frac{k}{2}}}{1+\frac{k}{2}}\right) +
d_2\sqrt{-\frac{-\pi t}{2+k}} I_{-\frac{1}{2k}} \left(\frac{(-t)^{1+\frac{k}{2}}}{1+\frac{k}{2}}\right) & \mbox{if} \ \ t < 0,\\
\end{cases}
\end{equation}
where $I_\nu$ is the modified Bessel function of the first kind,
\begin{equation}
\begin{aligned}
d_1 = \sqrt{\pi} \cos\left(\frac{\pi(1+k)}{2+k}\right) \csc\left(\frac{\pi}{2+k}\right)\ \ \mbox{and}\ \ 
d_2 = \frac{\sqrt{\pi}\csc\left(\frac{\pi}{2+k}\right)}{\sqrt{2+k}},
\end{aligned}
\end{equation}
form a basis in the space of solutions of (\ref{experiments:high:ode}) when $k$ is odd.

In each experiment, a value of $k$ was fixed and the algorithm of this paper
was applied twice to evaluate
\begin{equation}
f_k^{\mbox{\tiny high}}(t) = 
\begin{cases}
u_{k}^{\mbox{\tiny even}}(t) + i v_{k}^{\mbox{\tiny even}}(t) & \mbox{if k is even}\\
u_{k}^{\mbox{\tiny odd}}(t) + i v_{k}^{\mbox{\tiny odd}}(t) & \mbox{if k is odd}
\end{cases}
\label{experiments:high:f}
\end{equation}
at $200$ equispaced points.  The first time the algorithm was executed,
the values of the derivative of the coefficient in (\ref{experiments:high:ode})
were supplied, and during the second application of the algorithm,
the derivative of the coefficient was evaluated using spectral differentiation.
The results appear in Figure~\ref{figure:highplots1}.
Each plot there corresponds to one value of $k$ and gives the relative errors
in the calculated values of $f_k^{\mbox{\tiny high}}(t)$ incurred by each variant of the algorithm, as well as the relative
errors predicted by the condition number of evaluation of $f_k^{\mbox{\tiny high}}$.

\end{subsection}

\begin{subsection}{An equation whose coefficient has two bumps}
\label{section:experiments:bumps}

In this experiment, we considered the boundary value problem
\begin{equation}
\left\{
\begin{aligned}
y''(t) + \nu^2 q(t) y(t) &= 0, \ \ -10 < t < 10,\\
y(0) &= 0 \\
y'(10) &= 1,
\end{aligned}
\right.
\label{experiments:bumps:problem}
\end{equation}
where
\begin{equation}
q(t) = \exp\left(-(t-5)^2\right) + \exp\left(-(t+5)^2\right) + \frac{\sin\left(\frac{t}{2}\right)^2}{1+t^2}.
\end{equation}
The function $q$, a graph of which appears on the left-hand side of 
Figure~\ref{experiments:bumps:plot}, has a single zero at the point $t=0$.
A plot of the solution of (\ref{experiments:bumps:problem}) when $\nu=50$
appears on the right-hand side of Figure~\ref{experiments:bumps:plot}.

\begin{figure}[h!]
\begin{center}

\hfil
\includegraphics[width=.44\textwidth]{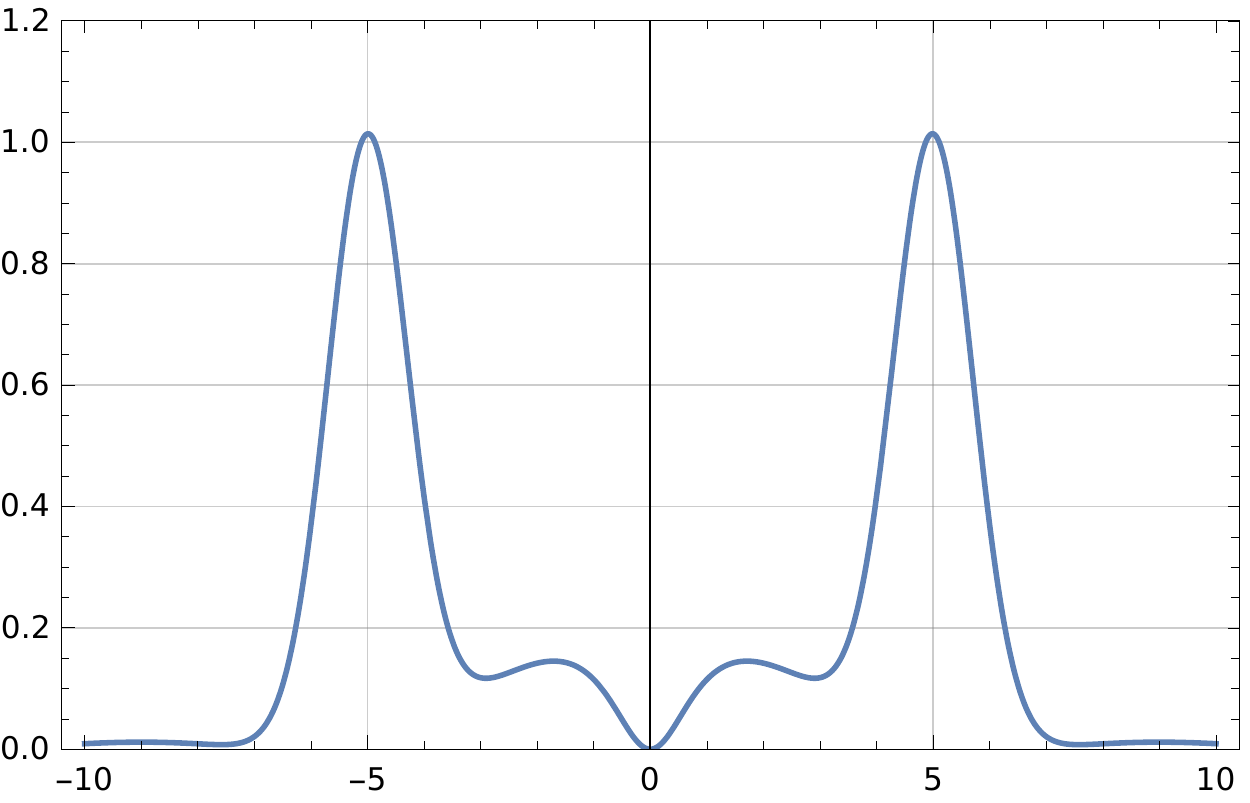}
\hfil
\includegraphics[width=.45\textwidth]{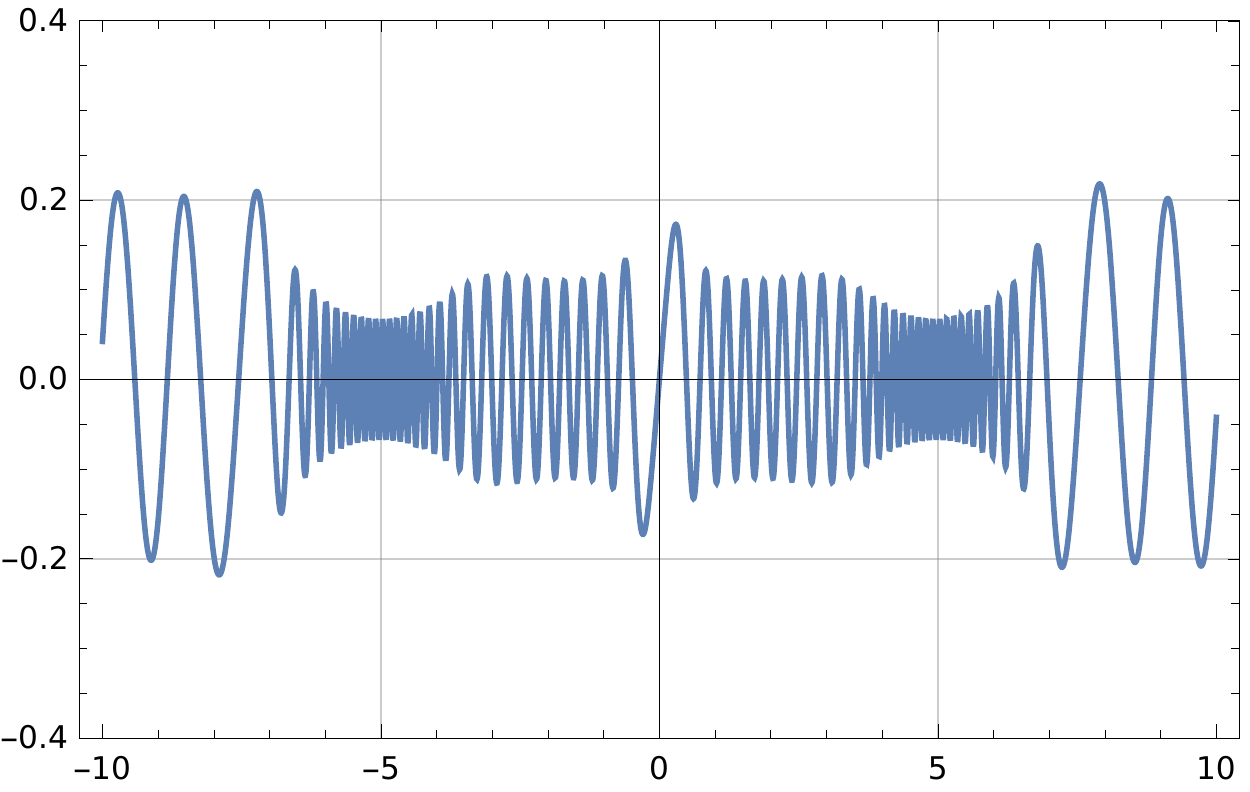}
\hfil

\end{center}
\caption{A graph of the coefficient $q(t)$ 
in the ordinary differential equation considered in Subsection~\ref{section:experiments:bumps}
appears on the left.  It has a single zero at the point $t=0$.
On the right is  a graph of the solution of the problem (\ref{experiments:bumps:problem})
when $\nu=50$.
  In the experiments of the paper, we consider values of $\nu$ as large
as $10^6$ --- we plot the solution in the case of a relatively small value of $\nu$ to
give the reader an indication of the behavior of the solutions and not as an illustration
of the performance of our algorithm.
}
\label{experiments:bumps:plot}
\end{figure}

In this experiment, we sampled $m=200$ points $x_1,\ldots,x_m$ in the interval $[0,6]$.  Then, for each
$\nu=10^{x_1},\ 10^{x_2},\ldots,10^{x_m},$
we solved (\ref{experiments:bumps:problem})
using the algorithm of this paper with the values of $q'$ specified.
We also solved  (\ref{experiments:bumps:problem}) by running the  standard solver
described in  Subsection~\ref{section:algorithm:odesolver} using extended precision arithmetic.
We then calculated the absolute difference in the obtained solutions at $5,000$ equispaced points
in the interval $(-10,10)$.    The plot on the left-hand side of 
Figure~\ref{figure:bumpsplots1} gives  the largest observed maximum absolute difference  as a function of $\nu$.  
The plot on the right-hand side gives the times required to solve (\ref{experiments:bumps:problem}) using
the phase method as a function of $\nu$. 


\end{subsection}

\begin{subsection}{An equation with three turning points}
\label{section:experiments:two}

In this experiment, we considered the problem
\begin{equation}
\left\{
\begin{aligned}
y''(t) + \nu^2 q(t) y(t) &= 0, \ \ -10 < t < 10,\\
y(0)  &= 1\\
y'(0) &= 0,
\end{aligned}
\right.
\label{experiments:two:problem}
\end{equation}
where
\begin{equation}
q(t) =  \exp\left(-(t+5)^2\right) -(t-5) \exp\left(-(t-5)^2\right) -6\exp(-25).
\end{equation}
The coefficient $q$ has three zeros on the interval $(-10,10)$:
\begin{equation}
\begin{aligned}
c_1 &\approx -9.817493179110059,\\
c_2 &=0\ \ \mbox{and}\\
c_3 &\approx 4.999999999916672.
\end{aligned}
\end{equation}
A graph of the function $q$ appears in  Figure~\ref{experiments:two:plot},
as does  a plot of the solution of (\ref{experiments:two:problem}) over
the interval $(-10,5.5)$ when $\nu=50$.  We truncated the plot of the solution
because its magnitude grows quickly in the interval $(c_3,10)$.

\begin{figure}[h!]
\begin{center}

\hfil
\includegraphics[width=.45\textwidth]{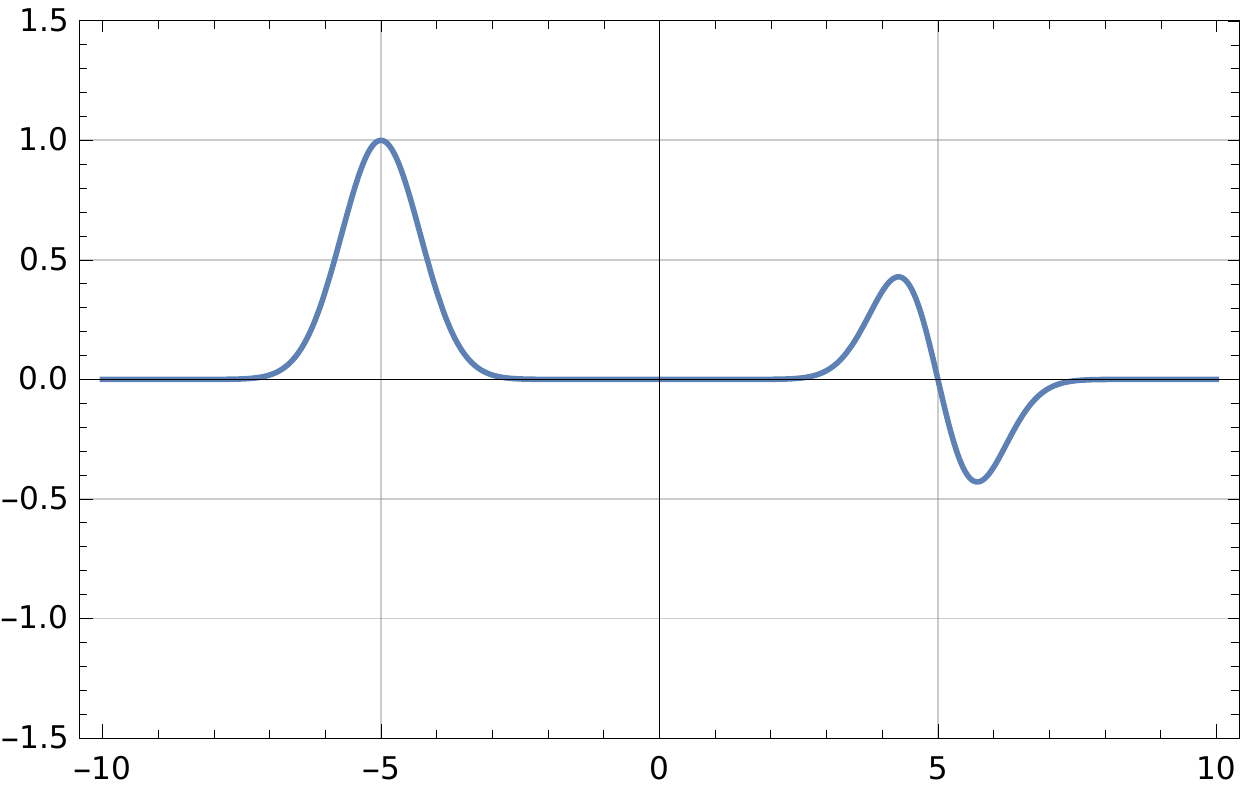}
\hfil
\includegraphics[width=.44\textwidth]{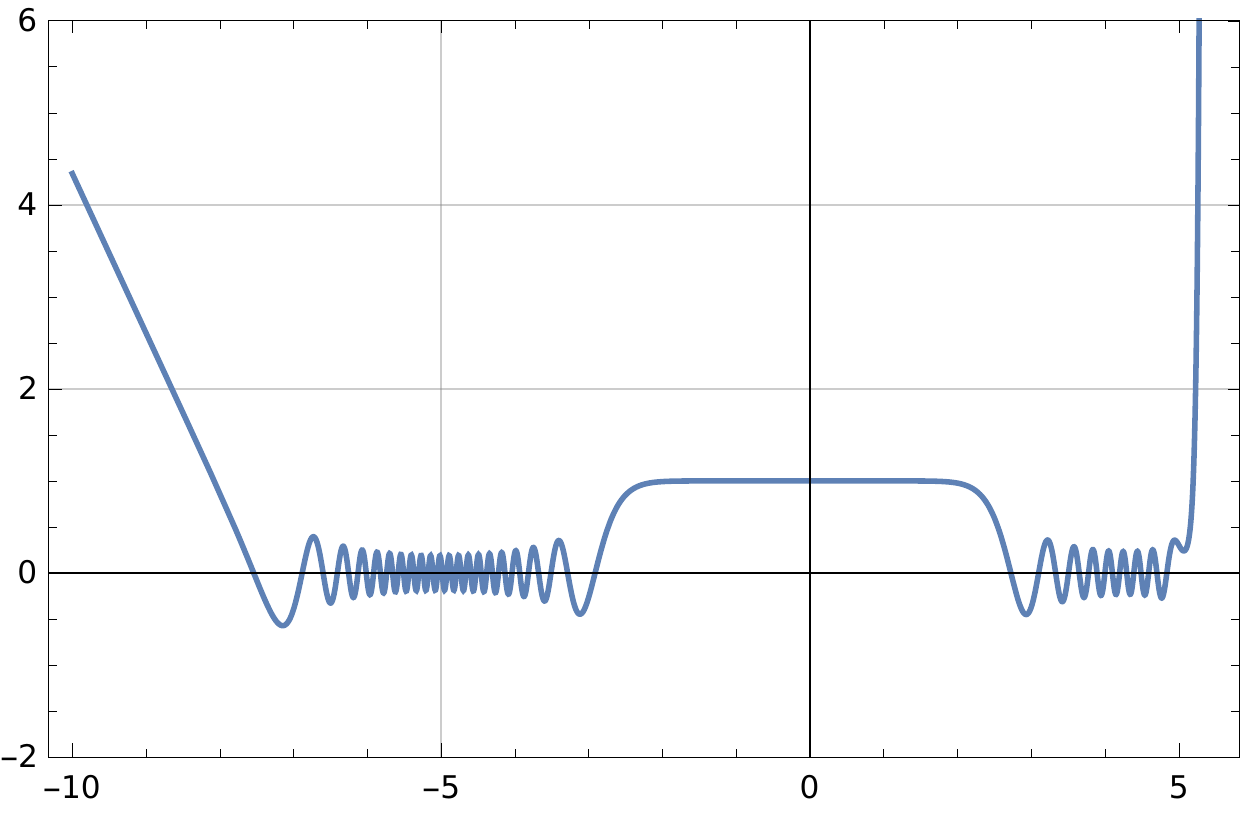}
\hfil

\end{center}
\caption{A graph of the coefficient $q(t)$ 
in the ordinary differential equation considered in Subsection~\ref{section:experiments:two}
appears on the left.
On the right is  a graph of the solution of the problem (\ref{experiments:two:problem})
over the interval $(-10,5.5)$ when $\nu=50$.  The plot was truncated because the magnitude
of the solution grows quickly in the interval $(c_3,10)$.
  In the experiments of the paper, we consider values of $\nu$ as large
as $10^5$ --- we plot the solution in the case of a relatively small value of $\nu$ 
to give the reader an indication of the behavior of the solutions and not as an illustration
of the performance of our algorithm.
}
\label{experiments:two:plot}
\end{figure}

In this experiment, we first sampled $m=200$ points $x_1,\ldots,x_m$ in the interval $[0,5]$.  Then, for each
$\nu=10^{x_1},\ 10^{x_2},\ldots,10^{x_m}$, we  solved (\ref{experiments:two:problem})
using a phase function method.
More explicitly, we constructed two phase functions, one defined on the interval $(-10,0)$ and the second defined
on the interval $(0,10)$.  The windowing algorithm was used on an interval
around the point $-5$ to determine the correct initial values for the phase
function on $(-10,0)$ and the windowing algorithm was applied on
an interval around the point $3$ to determine
the initial values for the phase function given on $(0,10)$.
The values of $q'$ were specified.
We next solved the problem (\ref{experiments:two:problem}) 
using the standard solver
described in  Subsection~\ref{section:algorithm:odesolver} running in extended
precision arithmetic.
Then, we evaluated the solution $y_1$ obtained with the phase method
and the solution $y_2$ constructed by the standard solver at  $r=5,000$ equispaced
points $z_1,z_2,\ldots,z_r$ on the interval $(-10,10)$ and measured the quantity
\begin{equation}
\xi_\nu = \max_{1\leq j \leq 5000} \frac{\left|y_1(z_j)-y_2(z_j)\right|}{1+\left|y_2(z_j)\right|}.
\label{experiments:two:error}
\end{equation}
Neither the absolute nor relative differences between the two solutions were appropriate for this problem because
the solution of (\ref{experiments:two:problem}) oscillates on part of the interval
and it is of large magnitude in another region.
Figure~\ref{figure:twoplots1} gives the results.   On the left-hand side
is a plot  of  $\xi_\nu$ as a function of $\nu$, while  
the plot on the right-hand side gives the times required to solve (\ref{experiments:two:problem}) 
via the phase method as a function of $\nu$.

\end{subsection}

\begin{subsection}{An equation with many turning points}
\label{section:experiments:many}

In this experiment, we considered the problem
\begin{equation}
\left\{
\begin{aligned}
y''(t) + \nu^2 q(t) y(t) &= 0, \ \ -11 < t < 11,\\
y(0)  &= 1\\
y'(0) &= 1,
\end{aligned}
\right.
\label{experiments:many:problem}
\end{equation}
where
\begin{equation}
q(t) = 1 + \cos\left(\pi t\right).
\end{equation}
The coefficient $q$ has $12$ zeroes in the interval $(-11,11)$,
at the points 
\begin{equation*}
-11,-9,-7,\ldots,-1,1,\ldots,7,9,11.
\end{equation*}
A graph of the function $q$ appears in  Figure~\ref{experiments:many:plot},
as does  a plot solution of (\ref{experiments:many:problem})  when $\nu=20$.
\begin{figure}[h!]
\begin{center}

\hfil
\includegraphics[width=.45\textwidth]{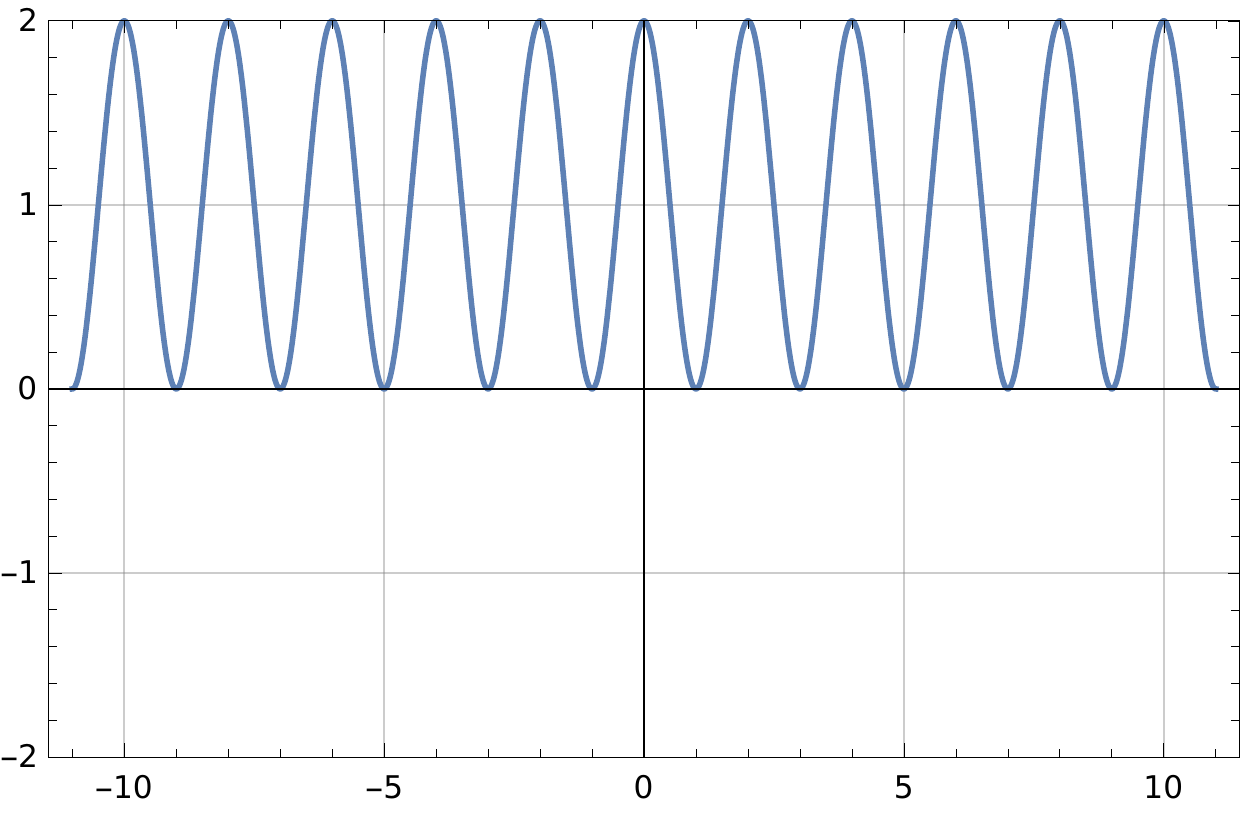}
\hfil
\includegraphics[width=.45\textwidth]{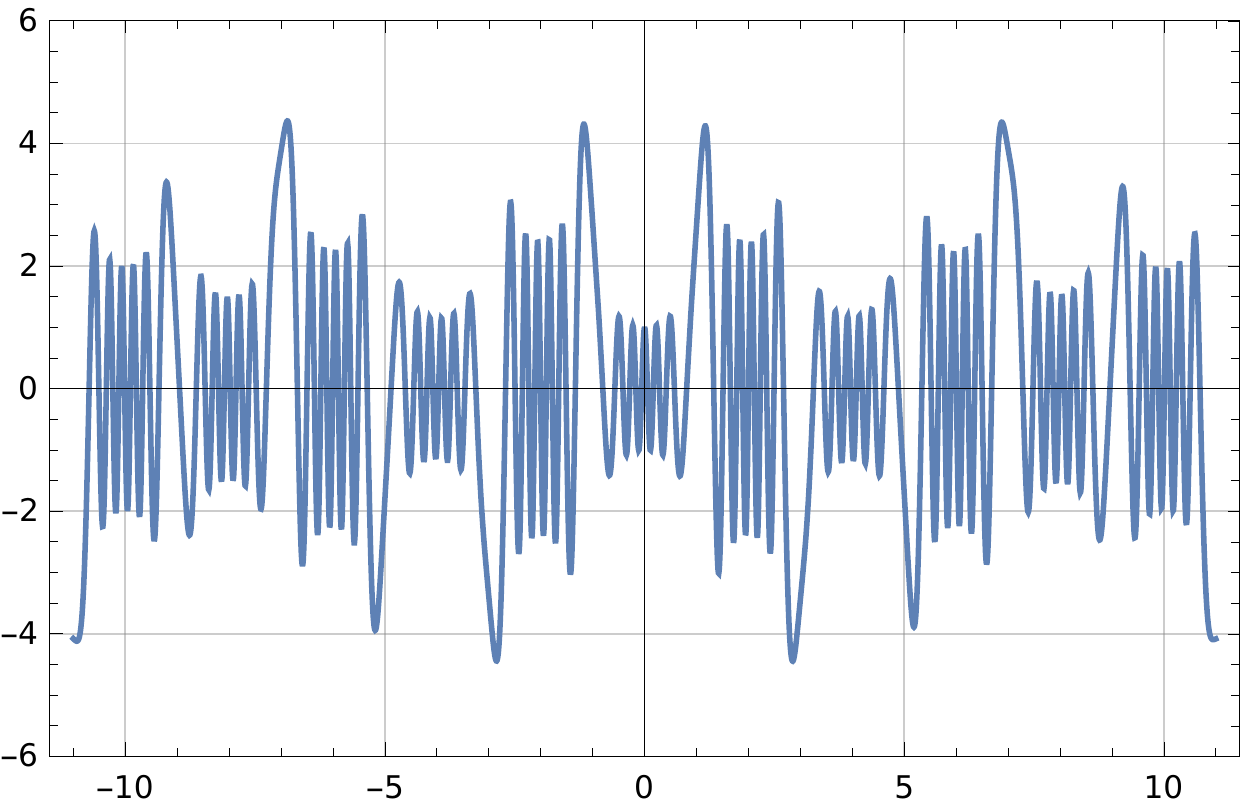}
\hfil

\end{center}
\caption{A graph of the coefficient $q(t)$ appearing
in the ordinary differential equation considered in Subsection~\ref{section:experiments:many}
appears on the left.
On the right is  a graph of the solution of the problem (\ref{experiments:many:problem})
when $\nu=20$.  In the experiments of the paper, we consider values of $\nu$ as large
as $10^5$ --- we plot the solution in the case of a relatively small value of $\nu$ to
give  the reader an indication of the behavior of the solutions and not as an illustration
of the performance of our algorithm.
}
\label{experiments:many:plot}
\end{figure}

We sampled $m=200$ points $x_1,\ldots,x_m$ in the interval $[0,5]$.  Then, for each
$\nu=10^{x_1},\ 10^{x_2},\ldots,10^{x_m}$, we solved (\ref{experiments:two:problem})
using a phase function method and with the standard solver
described in  Subsection~\ref{section:algorithm:odesolver} running in extended precision
arithmetic.
To solve (\ref{experiments:many:problem}) using phase functions, we constructed
eleven phase functions, one on each of the intervals
\begin{equation}
\left(-1+2k,1+2k\right), \ \ k=-5,-4,\ldots,-1,0,1,\ldots,4,5.
\end{equation}
The values of $q'$ were specified.  

Figure~\ref{figure:manyplots1} gives the results.   The plot on the left-hand side
gives the maximum absolute difference in the solutions obtained by the two solvers
observed while evaluating them at $5,000$ equispaced points in the interval $(-11,11)$ as a function of $\nu$.  
The plot on the right-hand side gives the times required to solve (\ref{experiments:many:problem}) using
the phase method as a  function of $\nu$.

\end{subsection}

\end{section}

\begin{section}{Conclusions}
\label{section:conclusions}

We have described a variant of the algorithm of \cite{BremerKummer}
for the numerical solution of a large
class of second order linear ordinary differential equations, including
many with turning points.
Unlike standard solvers, its running time is largely independent of the magnitude 
of the coefficients appearing in the equation, and, unlike asymptotic methods,
it consistently achieves accuracy on par with that indicated by the condition
number of the problem.     

One of the key observations of \cite{BremerKummer} is that
it is often preferable to calculate phase functions numerically
by simply solving the Riccati equation
rather than constructing asymptotic approximations of them \`a la
WKB methods.  A principal observation of this paper is that,
unlike standard asymptotic approximations of solutions
of the Riccati equation which become singular near turning points,
the phase functions themselves are perfectly well-behaved
near turning points, and this can be exploited to avoid 
using different mechanisms to represent solutions
in different regimes.  It appears from numerical
experiments that, in the vicinity of turning points, the cost of representing phase functions
via Chebyshev expansions and the like is comparable to that of representing the solutions
themselves.  And, perhaps  surprisingly, solutions can be represented to high
relative accuracy deep into the nonoscillatory regime via phase functions.

We have also presented experiments showing that the algorithm of this paper can be 
used to evaluate many families of special functions in time independent
of their parameters.  It should be noted, however,
that using  \emph{precomputed} piecewise polynomial expansions
of phase functions for the second order linear ordinary
differential equations satisfied by various families
of special functions is much faster.  Because these phase
functions are slowly varying, the expansions are surprisingly
compact.  We refer the interested reader to \cite{BremerALegendre}, where
such an approach is used to evaluate  the associated Legendre functions.

Phase functions methods are highly useful for performing several
other  computations involving special functions.  For instance, they can be used to calculate their zeros
and various generalized Gaussian quadrature rules \cite{BremerZeros},
and the author will soon report on a method for rapidly
applying special function transforms using phase functions.

In a future work, the author will describe an extension of the algorithm
of this paper which combines standard solvers for second order differential equations
with phase function methods.  The algorithm will operate by subdividing the solution domain
into intervals and constructing a local basis of solutions on each interval.
When the coefficients are of large magnitude, the local basis functions
will be represented via phase functions, and the basis functions  will be represented
directly, via piecewise Chebyshev expansions, on intervals in which the coefficients are of small magnitude.
One of the principal advantages of such an approach is that the solution interval
can be subdivided in a more or less arbitrary fashion and each computation can be performed
locally, thus allowing for large-scale parallelization.

\end{section}

\begin{section}{Acknowledgments}
The author would like to thank Alex Barnett and Fruzsina Agocs of the Flatiron Institute
and Kirill Serkh of the University of Toronto
for many thoughtful discussions and suggestions.
The author was supported in part by NSERC Discovery grant  RGPIN-2021-02613.
\end{section}

\bibliographystyle{acm}
\bibliography{odesolve.bib}

\vfil\eject
\vfil

\begin{figure}[h!!!!!!!!!!!!!!!!!!!!!!!!!!!!]
\centering
\hfil
\includegraphics[width=.43\textwidth]{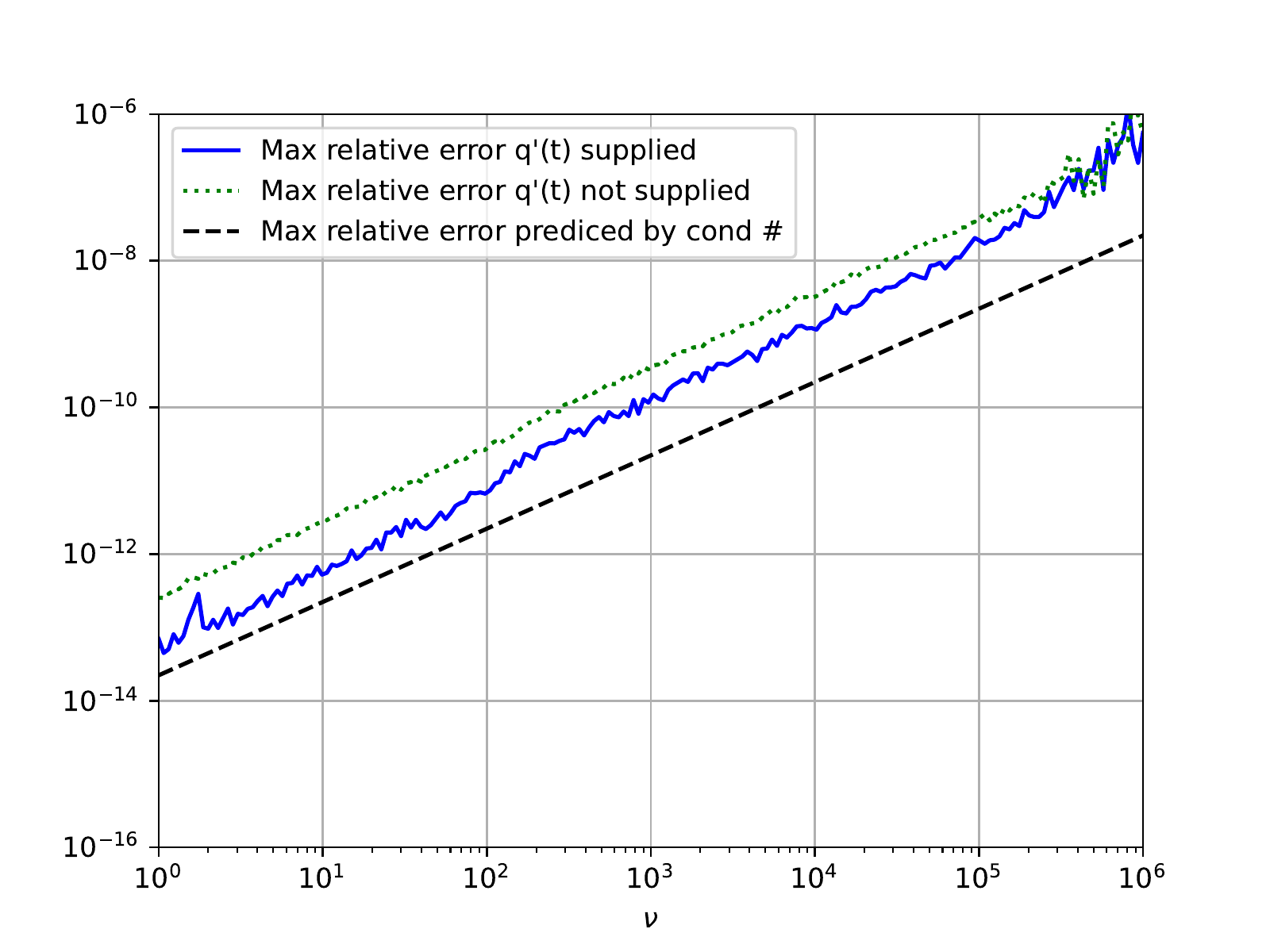}
\hfil
\includegraphics[width=.43\textwidth]{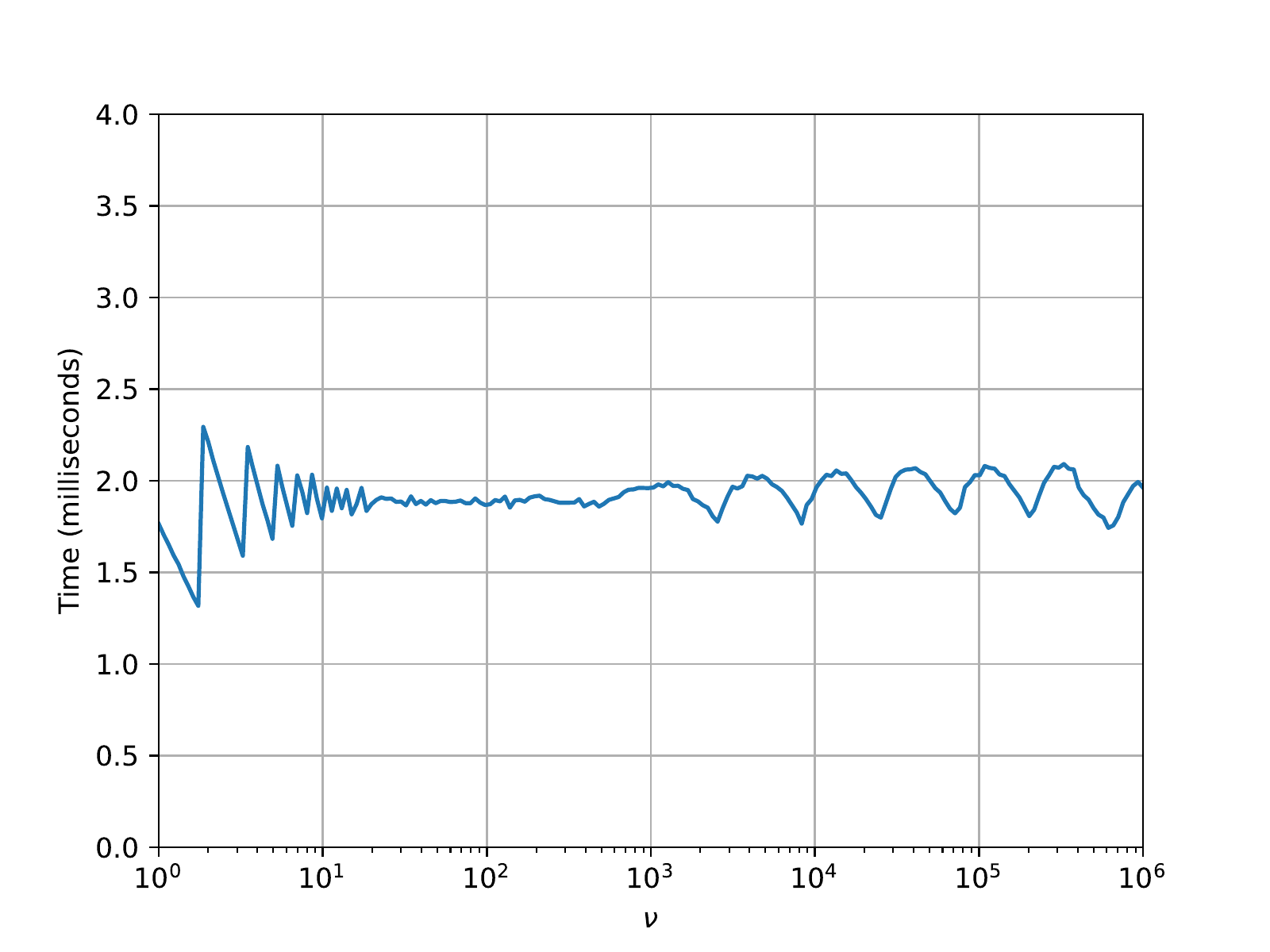}
\hfil

\caption{
The results of the first experiment of Section~\ref{section:experiments:bessel},
which concerned the Bessel functions.
The plot on the left gives the maximum relative errors which were observed
in the course of evaluating the Bessel function 
$f_\nu^{\mbox{\tiny bes}} = J_\nu(t) + i Y_\nu(t)$ at $5,000$ equispaced on the interval $(0,100\nu)$ 
using the phase functions $\alpha^{\mbox{\tiny bes}}_\nu$ 
and $\widetilde{\alpha}^{\mbox{\tiny bes}}_\nu$, 
as well as  the maximum relative error predicted by the condition number of evaluation of $f_\nu^{\mbox{\tiny bes}}$,
as functions of $\nu$.
The plot on the right gives the time (in milliseconds) which was required to construct the phase
function $\alpha^{\mbox{\tiny bes}}_\nu$  as a function of $\nu$.
}
\label{figure:besselplots1}
\end{figure}

\vfil\eject

\begin{figure}[h!!!]

\hfil
\includegraphics[width=.45\textwidth]{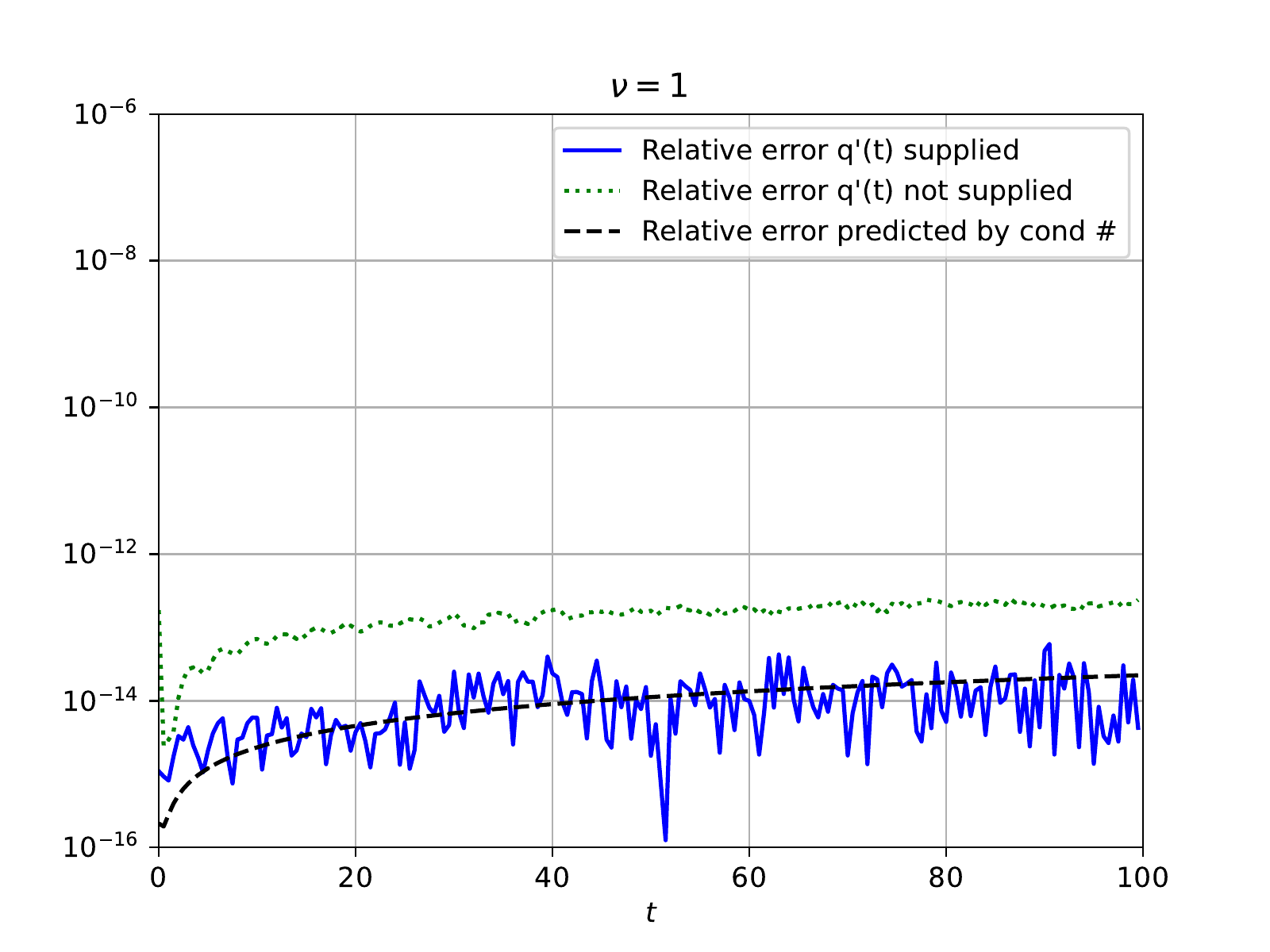}
\hfil
\includegraphics[width=.45\textwidth]{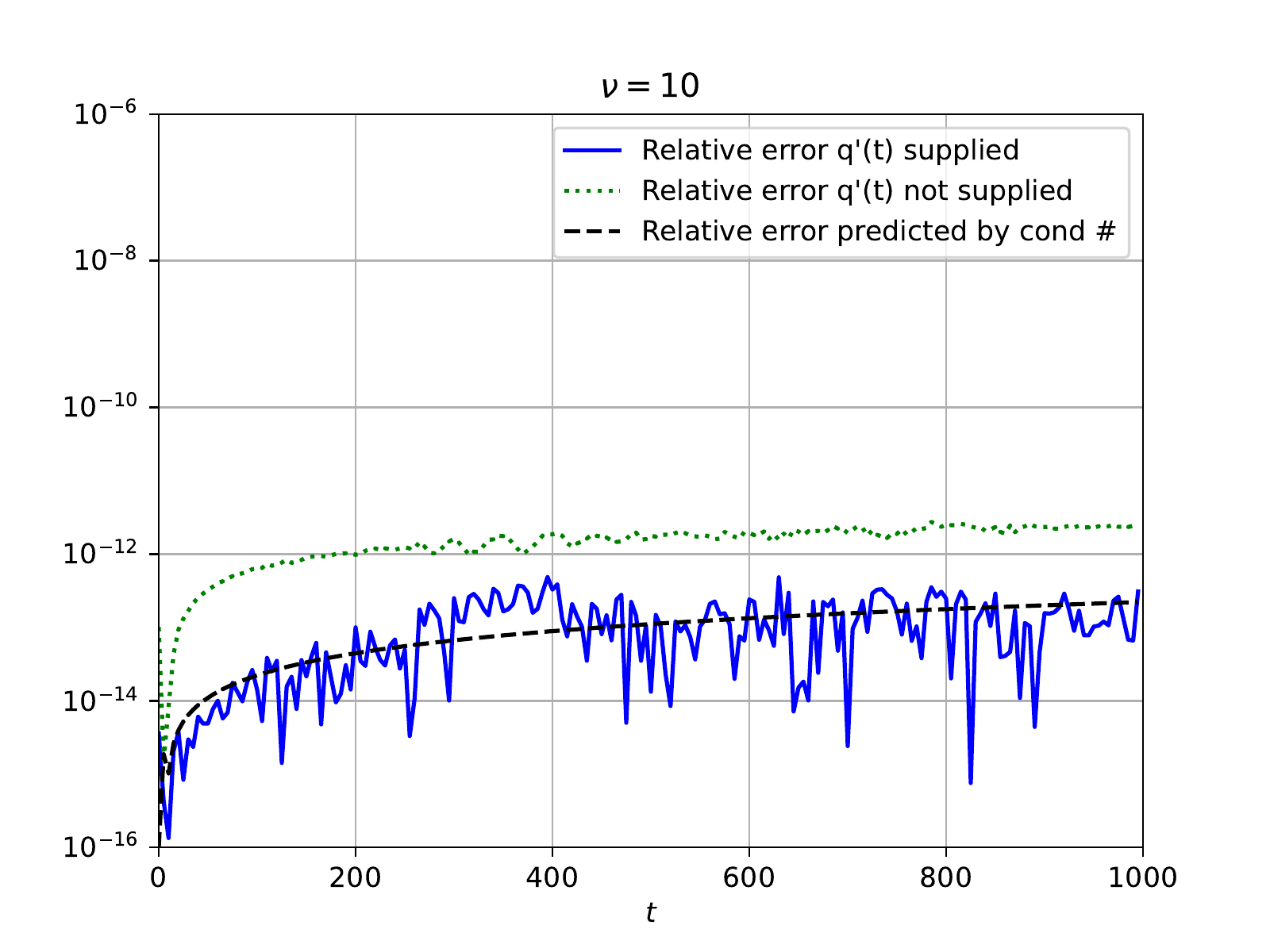}
\hfil
\vskip 3em

\hfil
\includegraphics[width=.45\textwidth]{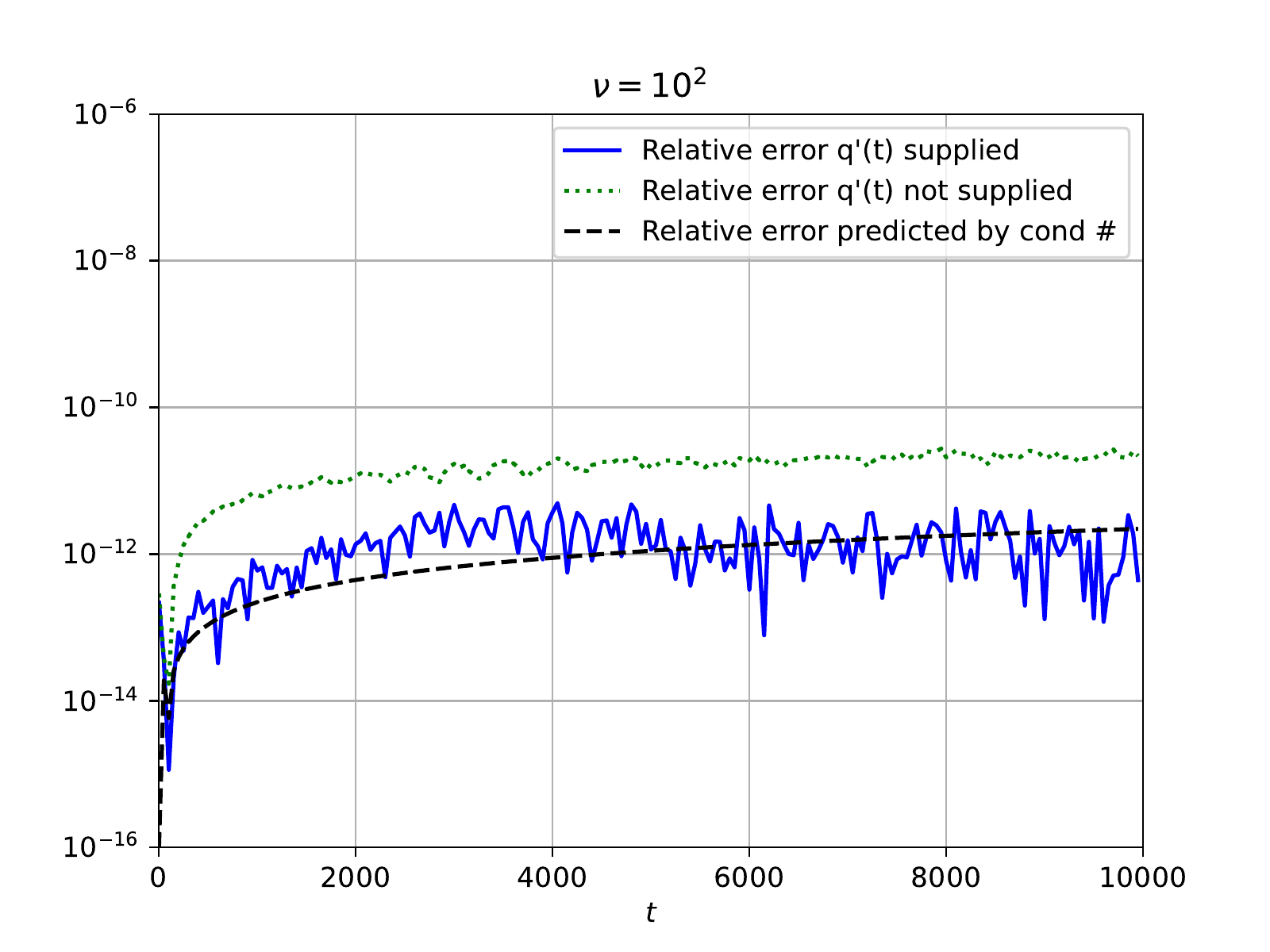}
\hfil
\includegraphics[width=.45\textwidth]{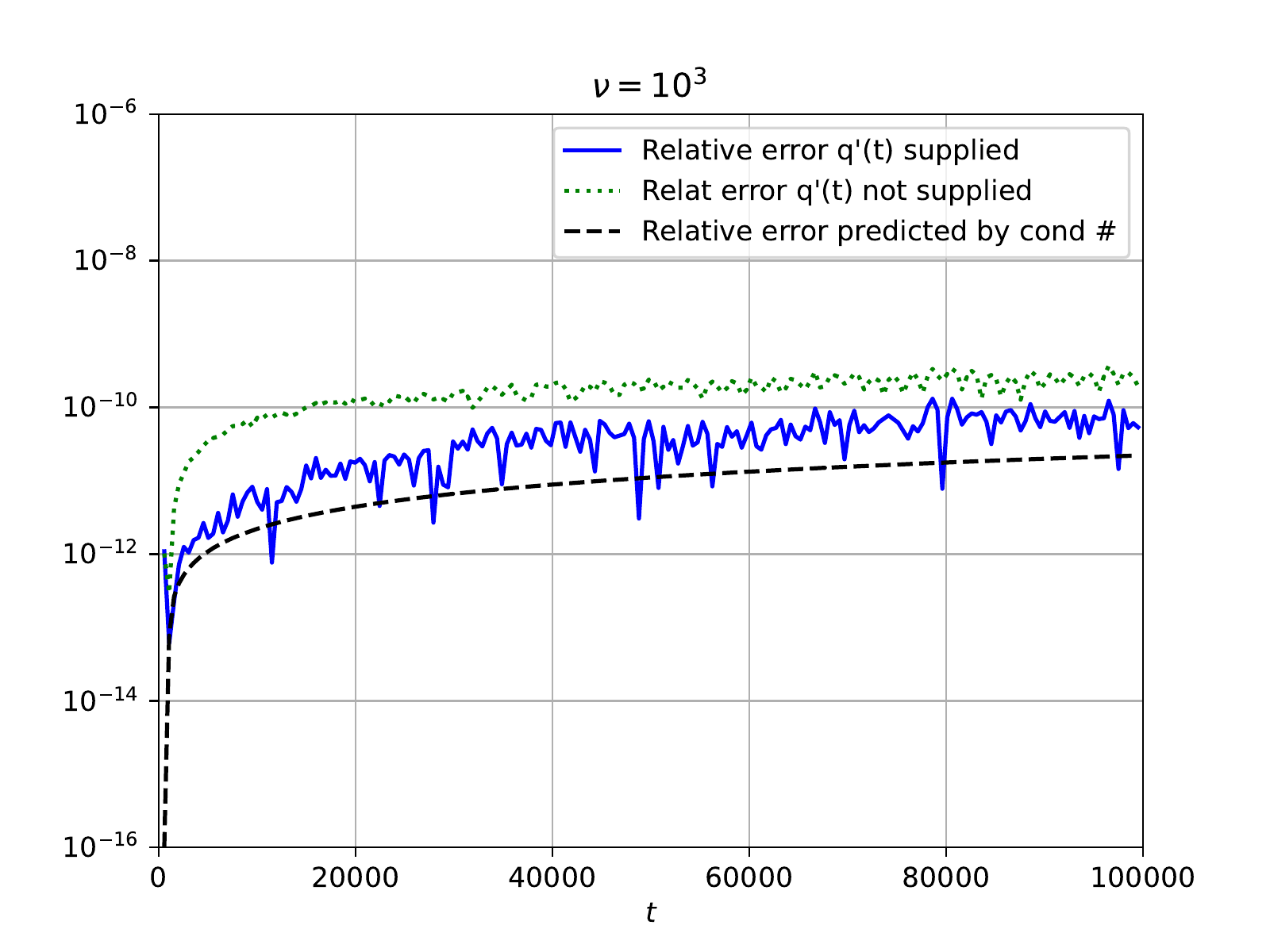}
\hfil
\vskip 3em

\hfil
\includegraphics[width=.45\textwidth]{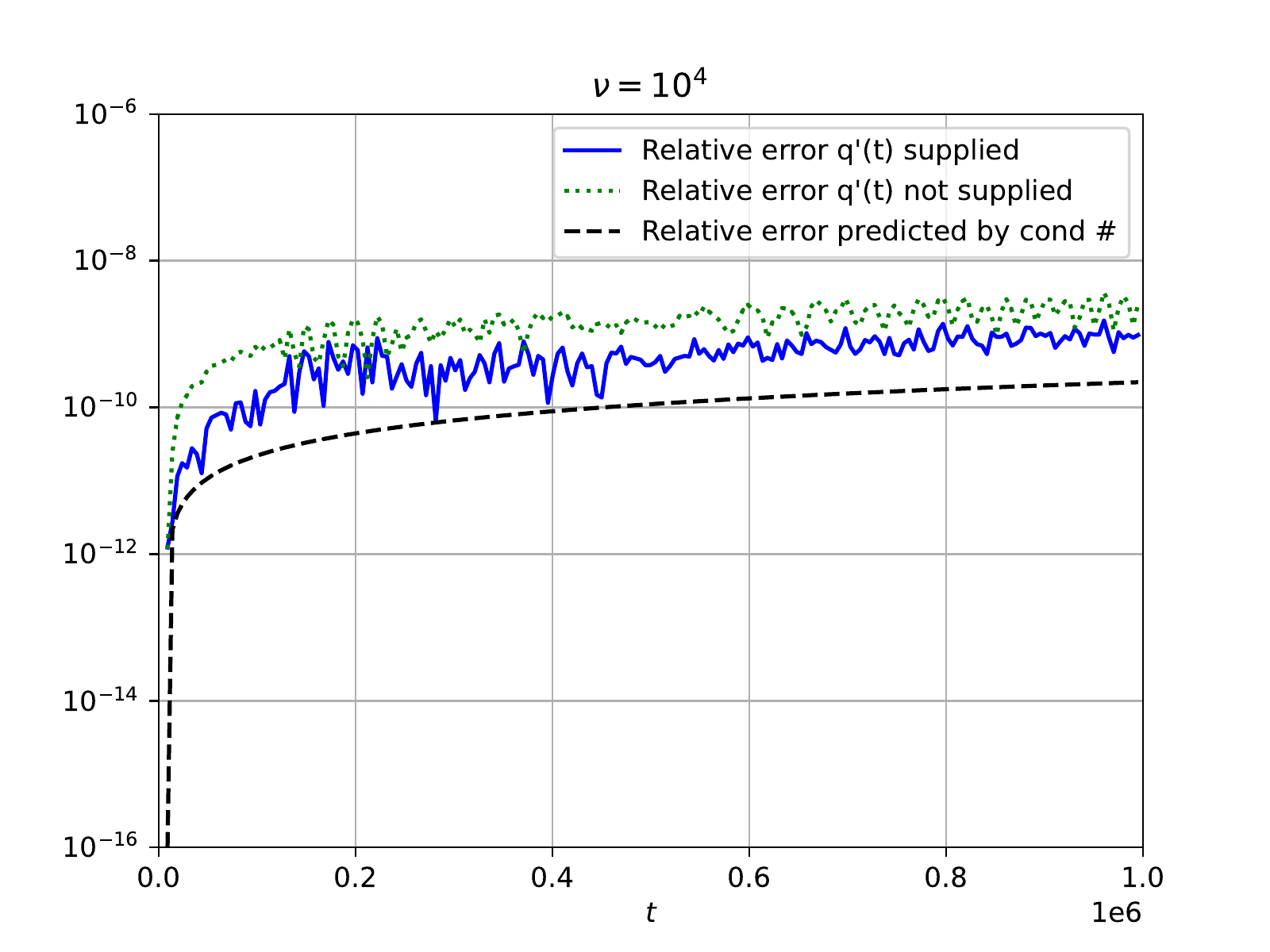}
\hfil
\includegraphics[width=.45\textwidth]{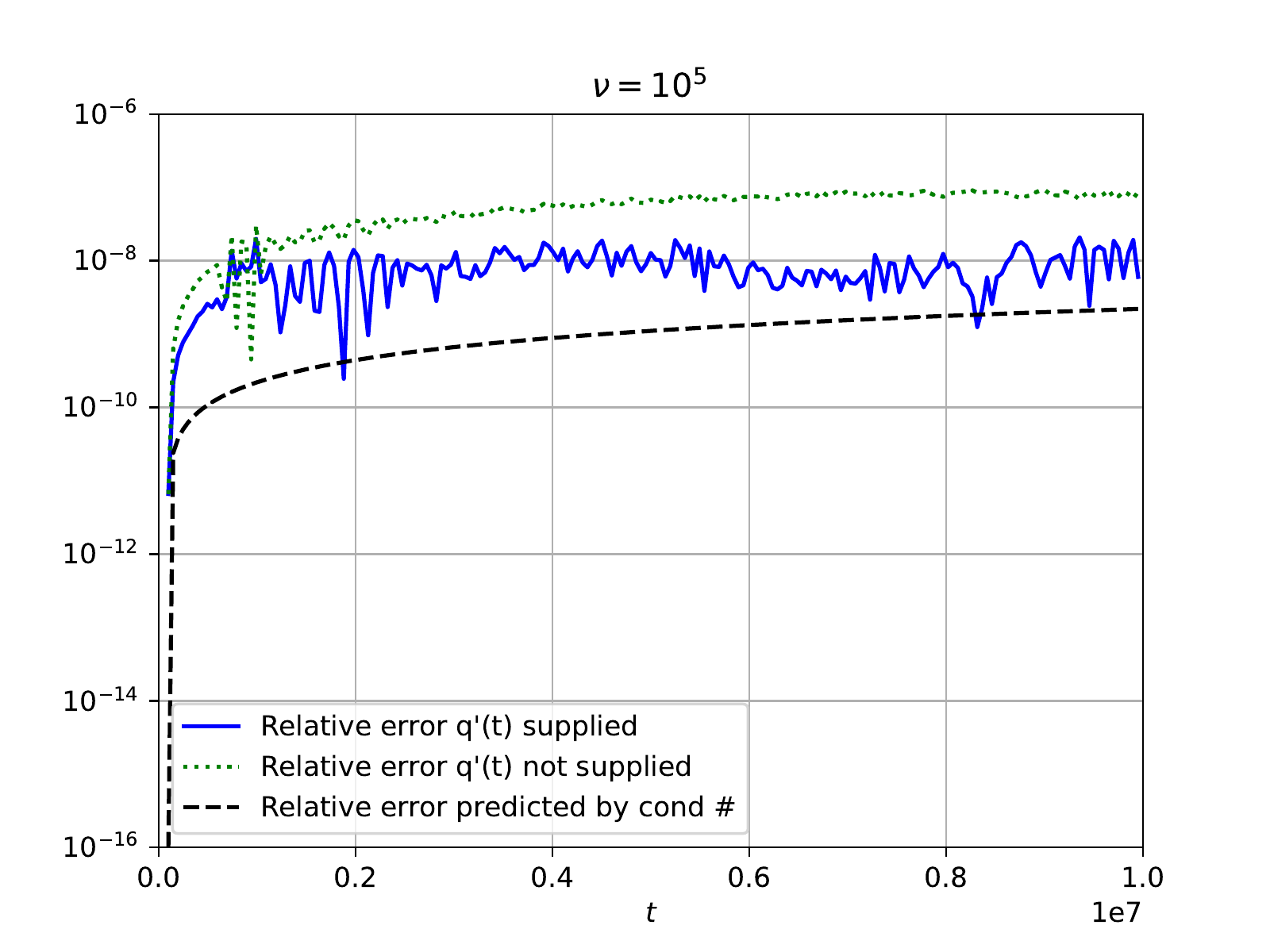}
\hfil

\caption{The results of the second set of experiments  discussed in 
Subsection~\ref{section:experiments:bessel}, which concerned the Bessel functions.
Each plot gives the relative errors in the values of  
$f_\nu^{\mbox{\tiny bes}}(t) = J_\nu(t) + i Y_\nu(t)$  calculated using
the phase functions $\alpha_\nu^{\mbox{\tiny bes}}$ and
$\widetilde{\alpha}_\nu^{\mbox{\tiny bes}}$  for a particular fixed value of $\nu$ as a function of $t$,
and compares them with the relative error predicted by the condition number of  evaluation of $f_\nu^{\mbox{\tiny bes}}$.
}
\label{figure:besselplots2}
\end{figure}

\vfil\eject

\begin{figure}[h!!!!!!!!!!!!!!!!!!!!!!!!!!!!!!!!!!]

\hfil
\includegraphics[width=.40\textwidth]{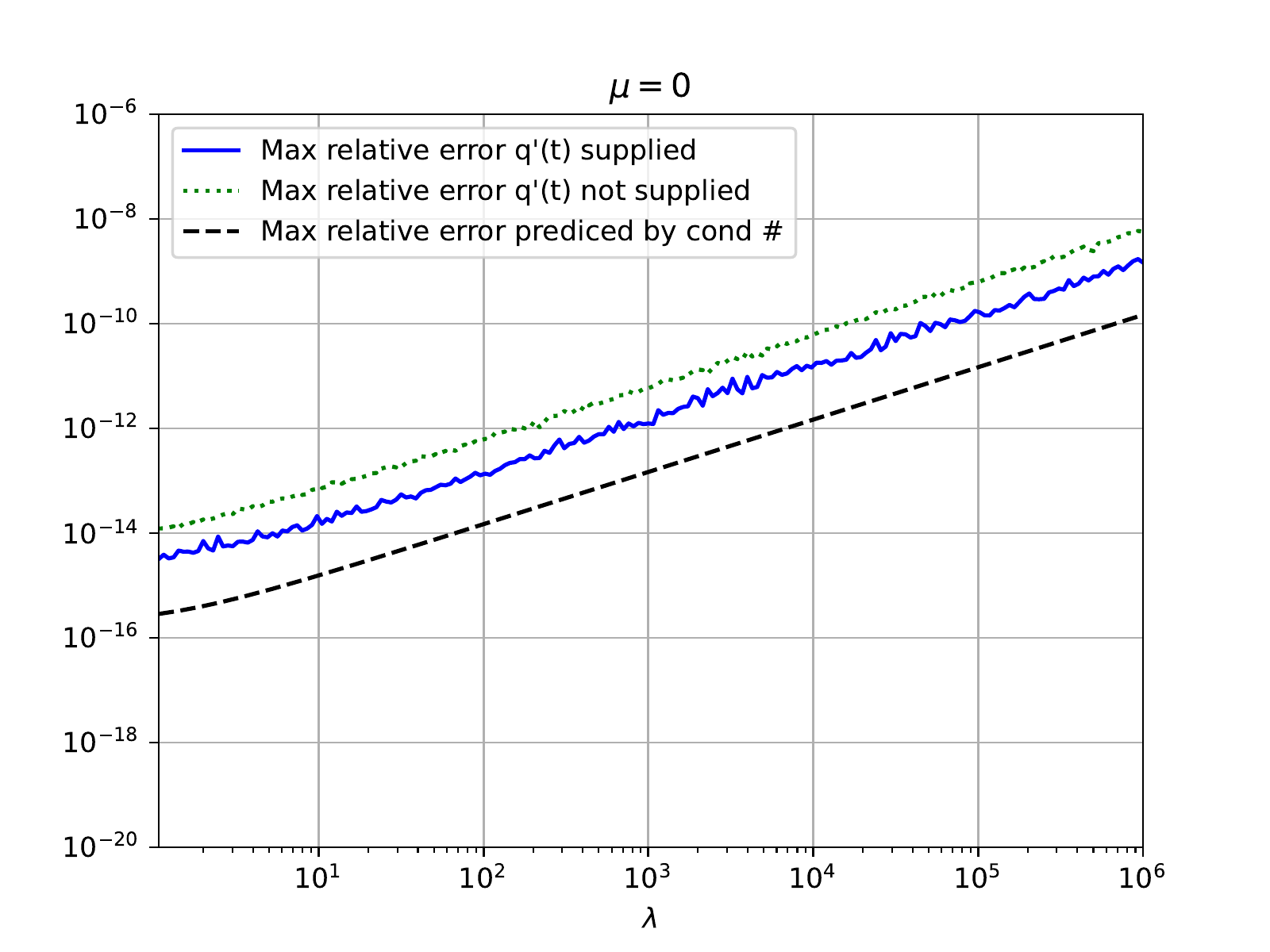}
\hfil
\includegraphics[width=.40\textwidth]{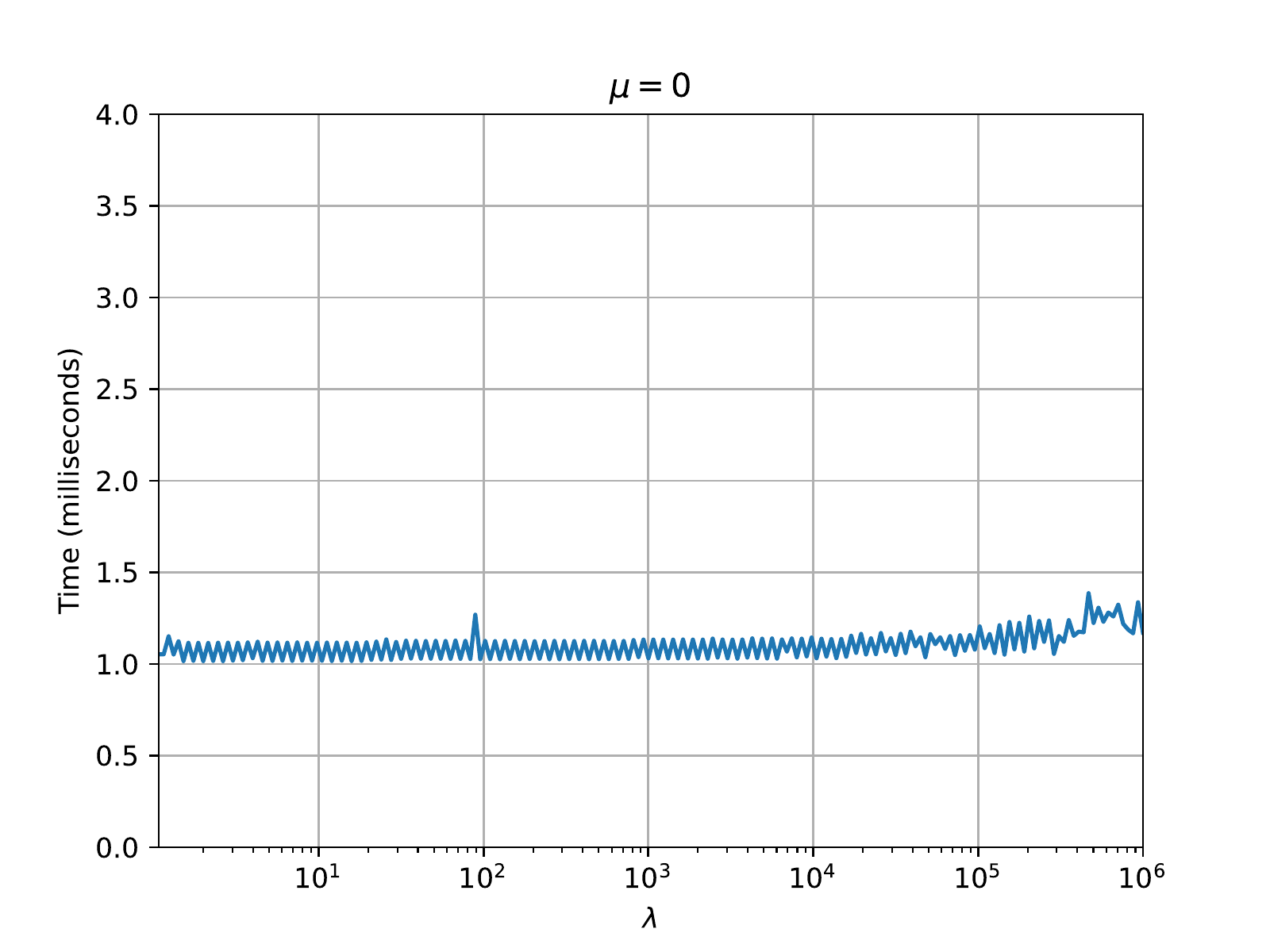}
\hfil

\hfil
\includegraphics[width=.40\textwidth]{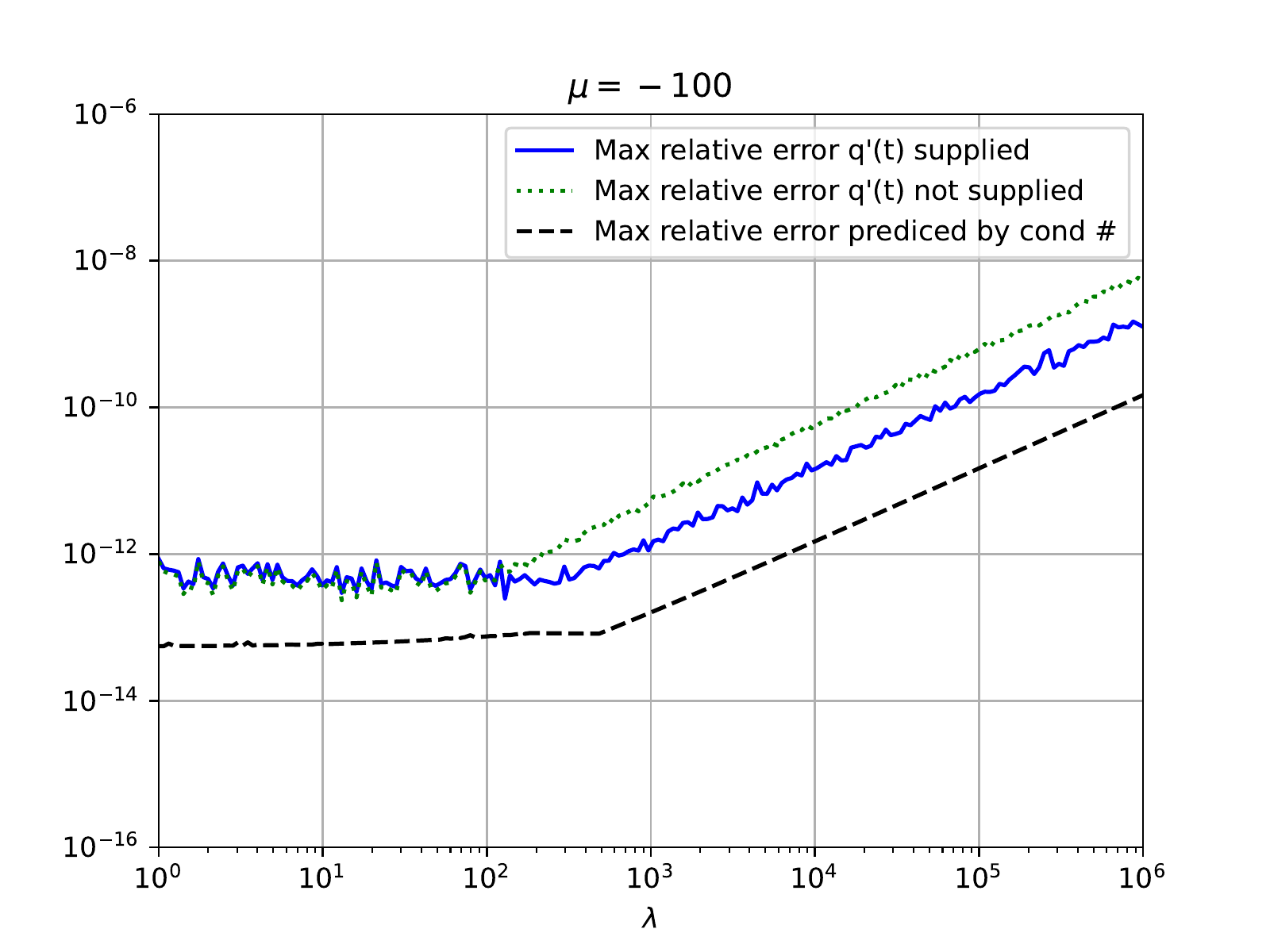}
\hfil
\includegraphics[width=.40\textwidth]{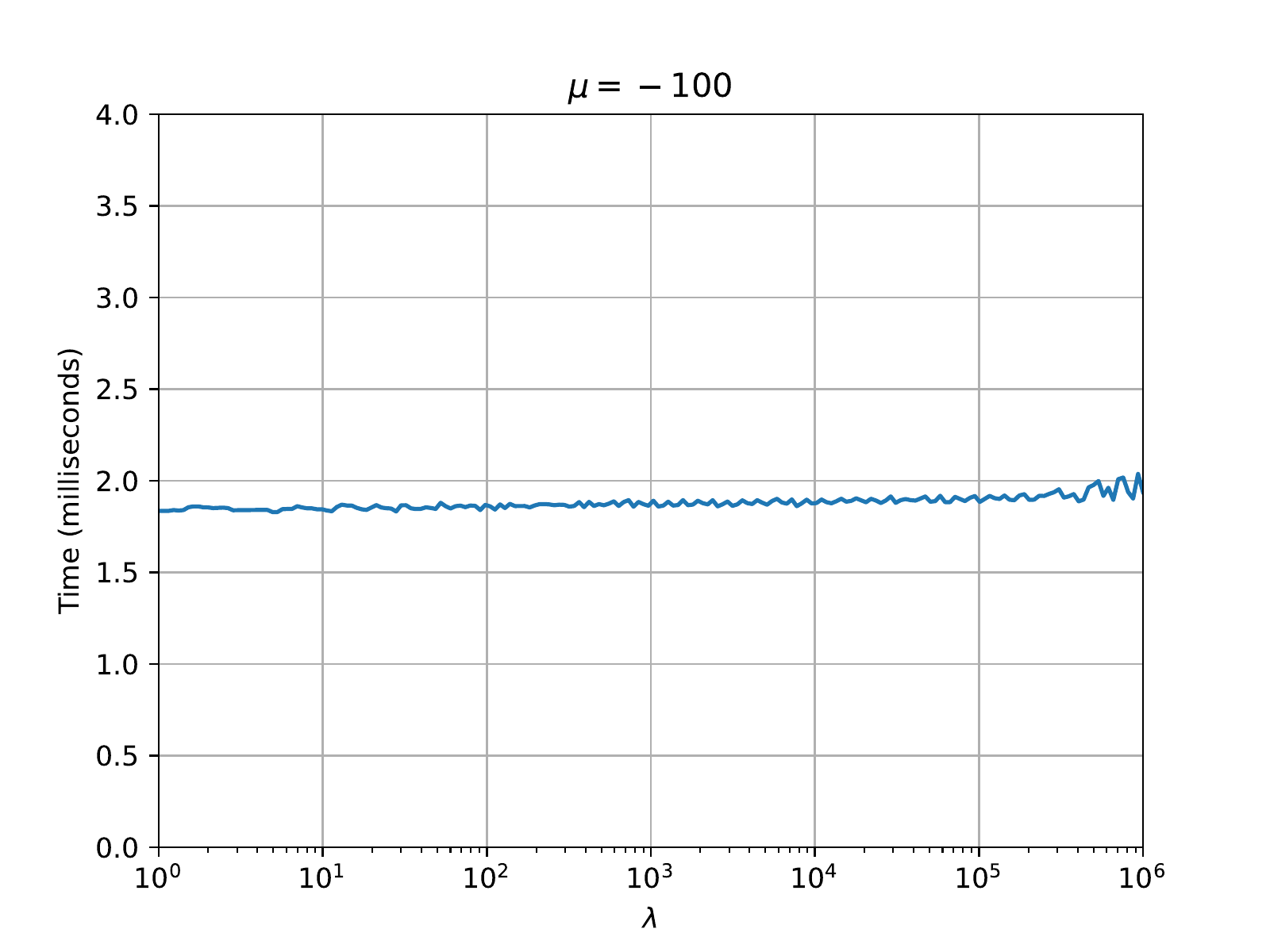}
\hfil

\hfil
\includegraphics[width=.40\textwidth]{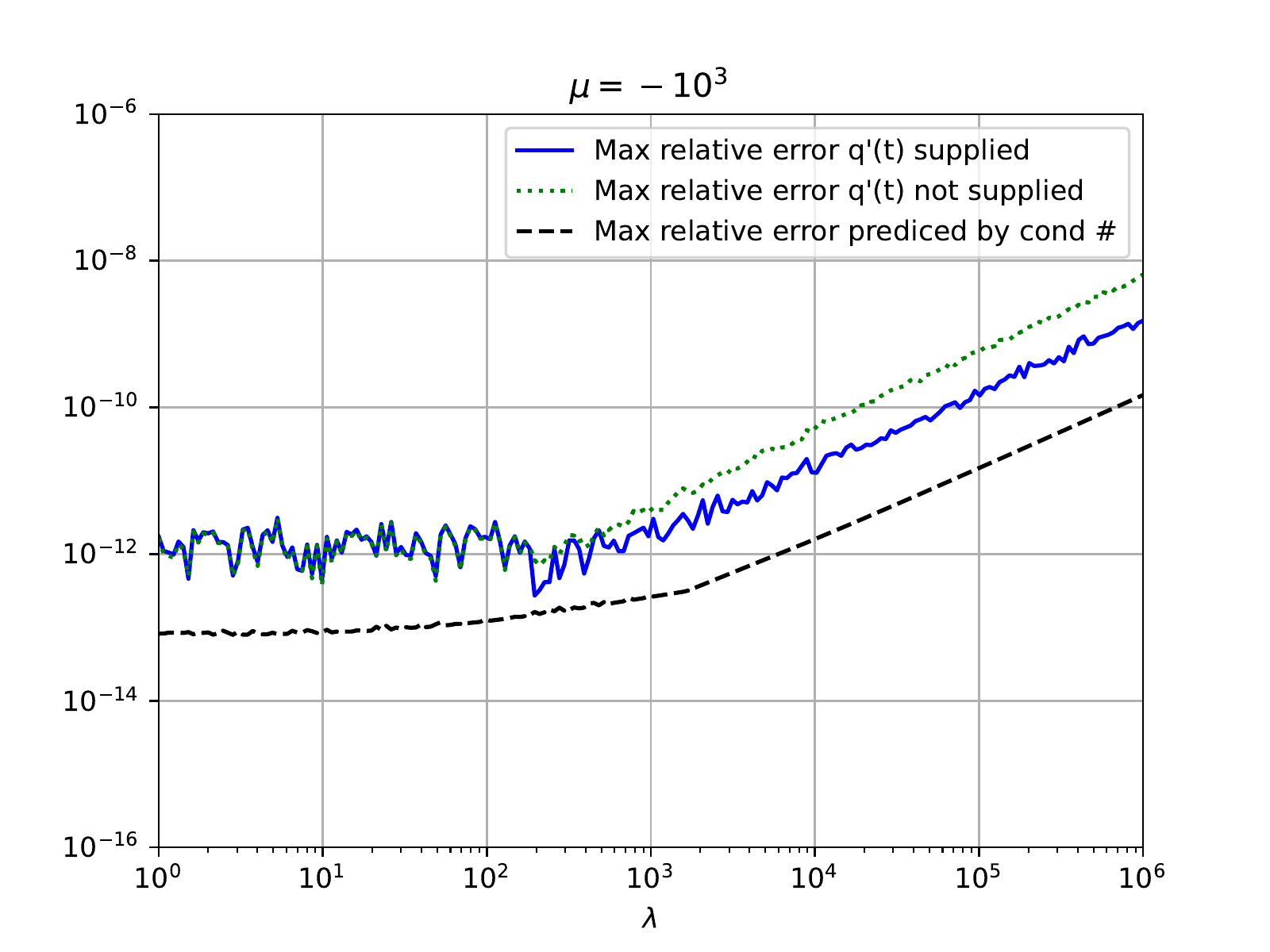}
\hfil
\includegraphics[width=.40\textwidth]{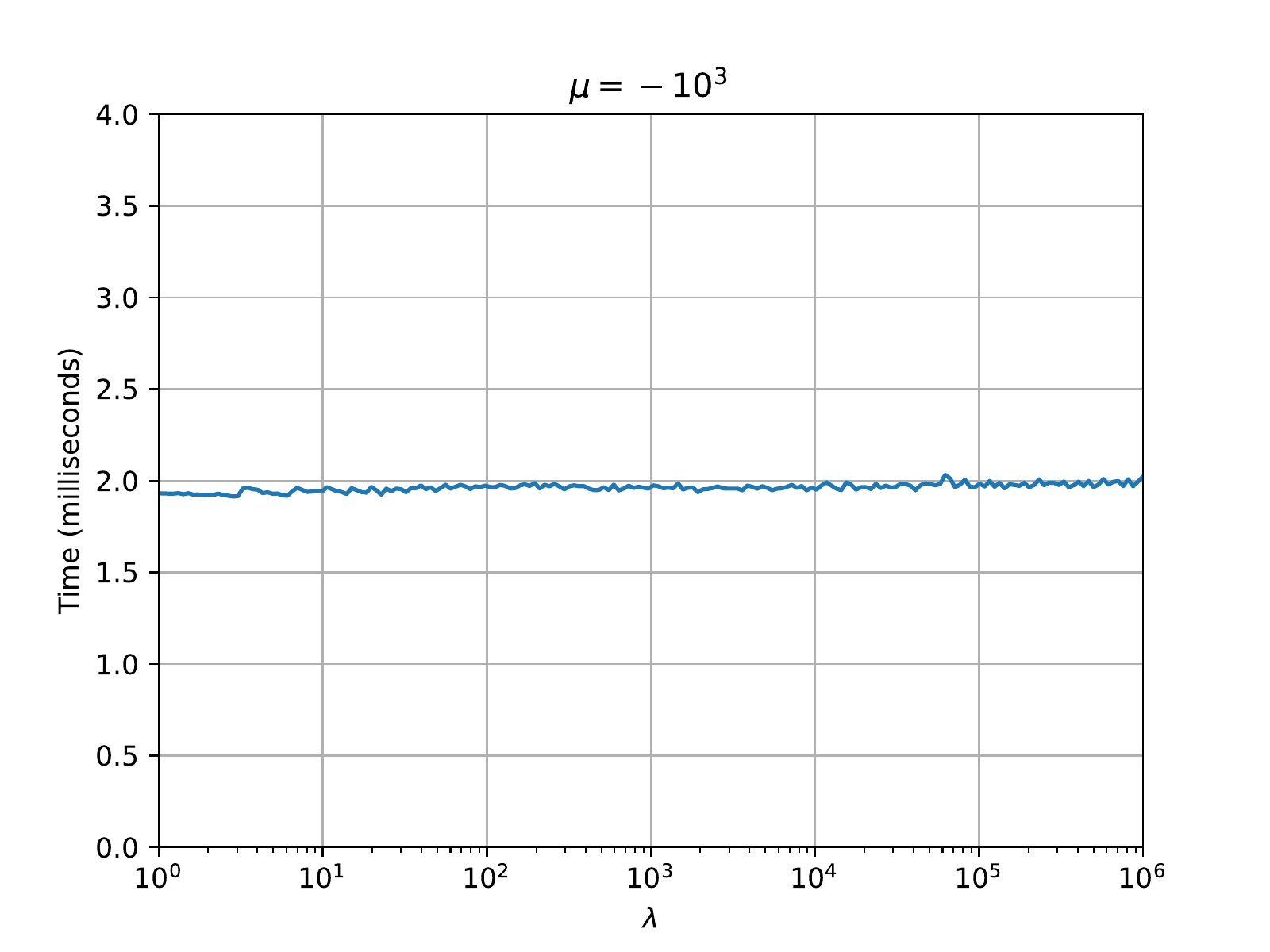}
\hfil

\hfil
\includegraphics[width=.40\textwidth]{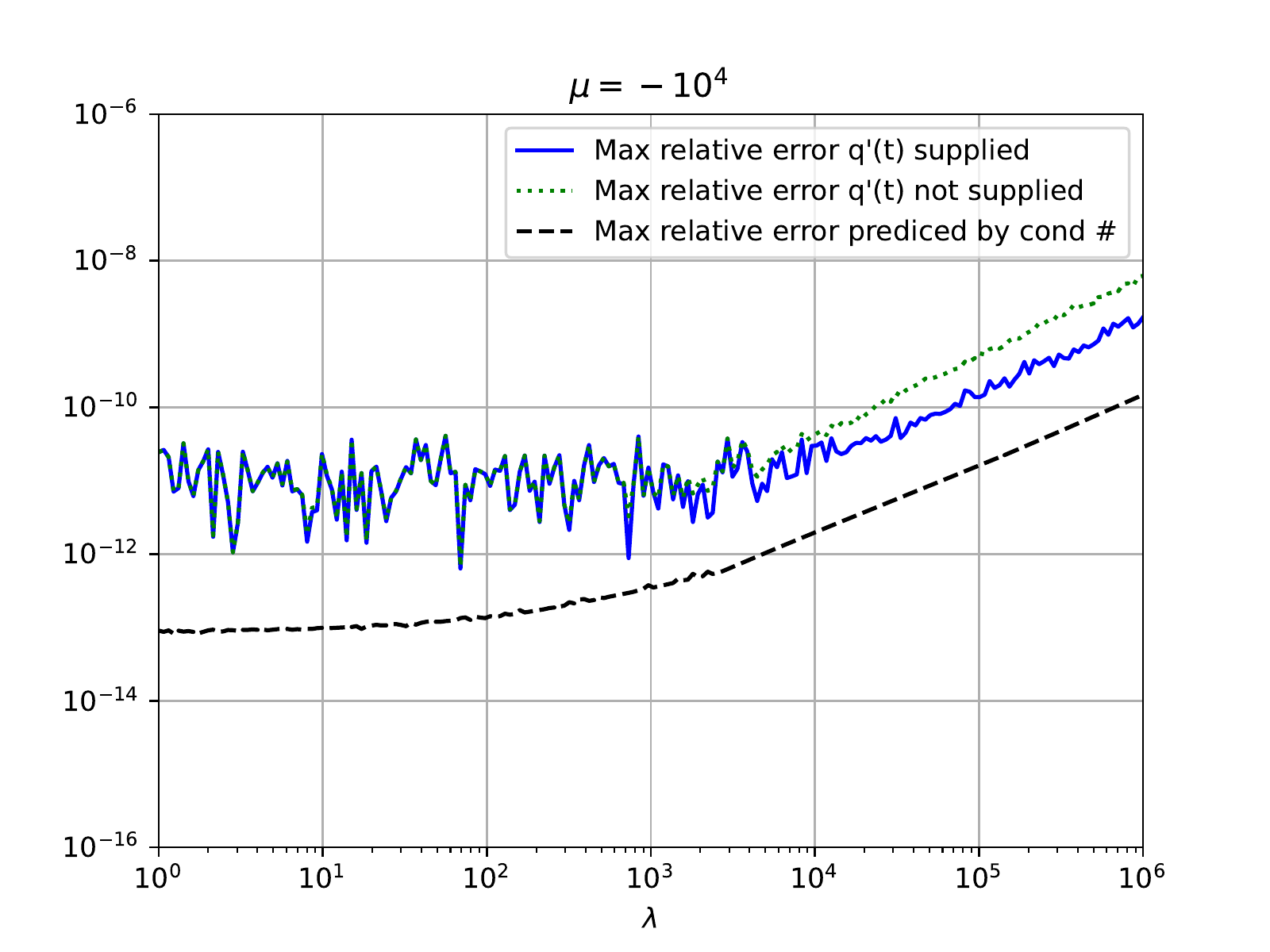}
\hfil
\includegraphics[width=.40\textwidth]{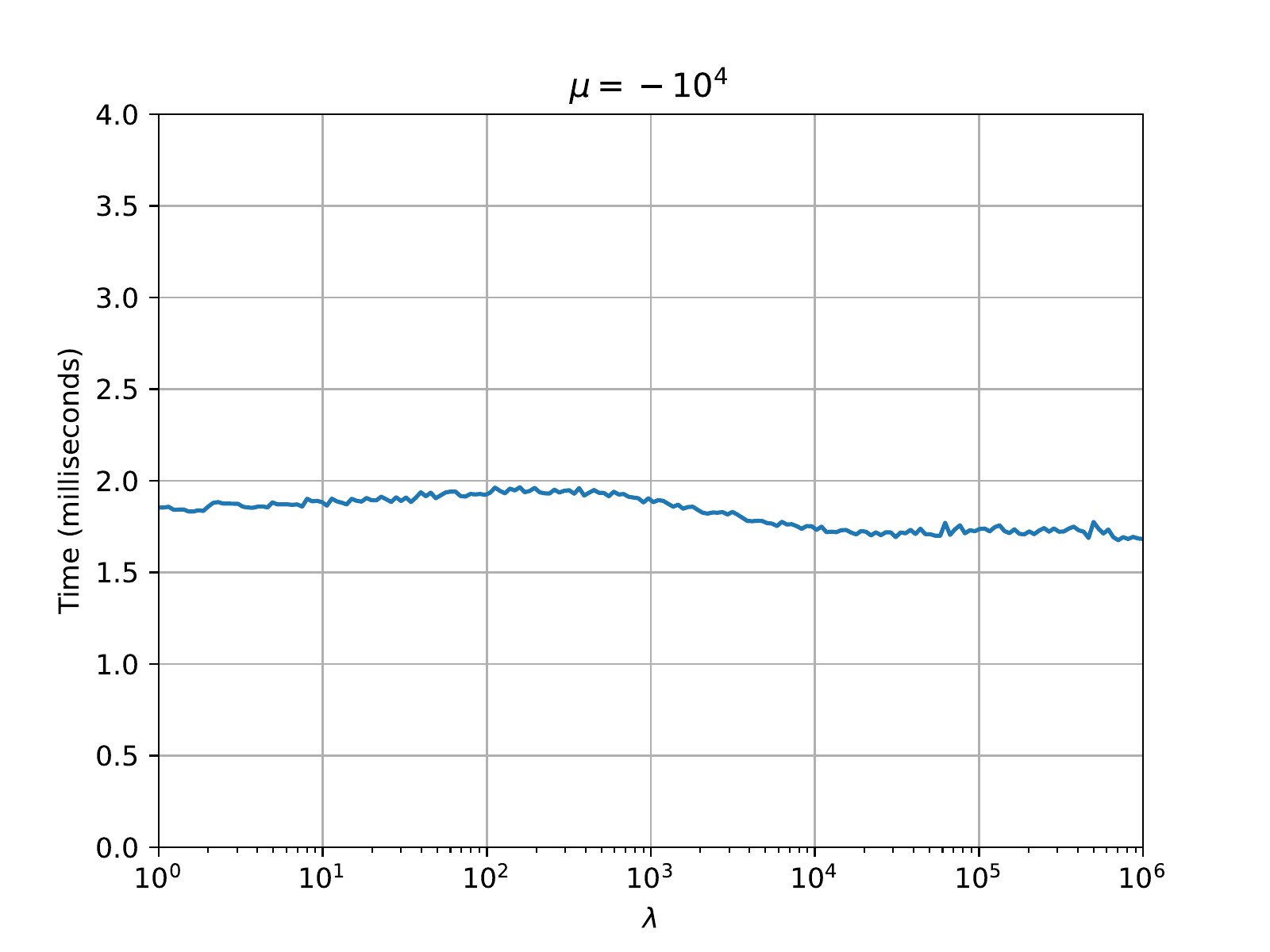}
\hfil

\caption{\small The results of the experiments discussed in  Subsection~\ref{section:experiments:alf},
which concerned the associated Legendre functions.
The left-hand side of each row contains a plot giving the maximum relative errors observed while  calculating the values of 
$f^{\mbox{\tiny alf}}_{\lambda,\mu}(w) = \widetilde{P}_{|\mu|+\lambda}^\mu(\tanh(w)) + i \frac{2}{\pi} \widetilde{Q}_{|\mu|+\lambda}^\mu(\tanh(w))$
at $5,000$ points using each of the phase functions
$\alpha_{|\mu|+\lambda,\mu}^{\mbox{\tiny alf}}$ and $\widetilde{\alpha}_{|\mu|+\lambda,\mu}^{\mbox{\tiny alf}}$, as well as the 
the maximum relative error predicted by the condition number of evaluation of  
$f^{\mbox{\tiny alf}}_{|\mu|+\lambda,\mu}(w)$,
as functions of $\lambda$ for a fixed value of $\mu$.
The right-hand side of each row contains a plot giving the time (in milliseconds) required
to construct the phase function representing  $f^{\mbox{\tiny alf}}_{|\mu|+\lambda,\mu}(w)$ as a function of $\lambda$.
}
\label{figure:alfplots1}
\end{figure}

\vfil\eject


\begin{figure}[h!!!!!!!!!!!!!!!!!!!!!!!!]

\hfil
\includegraphics[width=.40\textwidth]{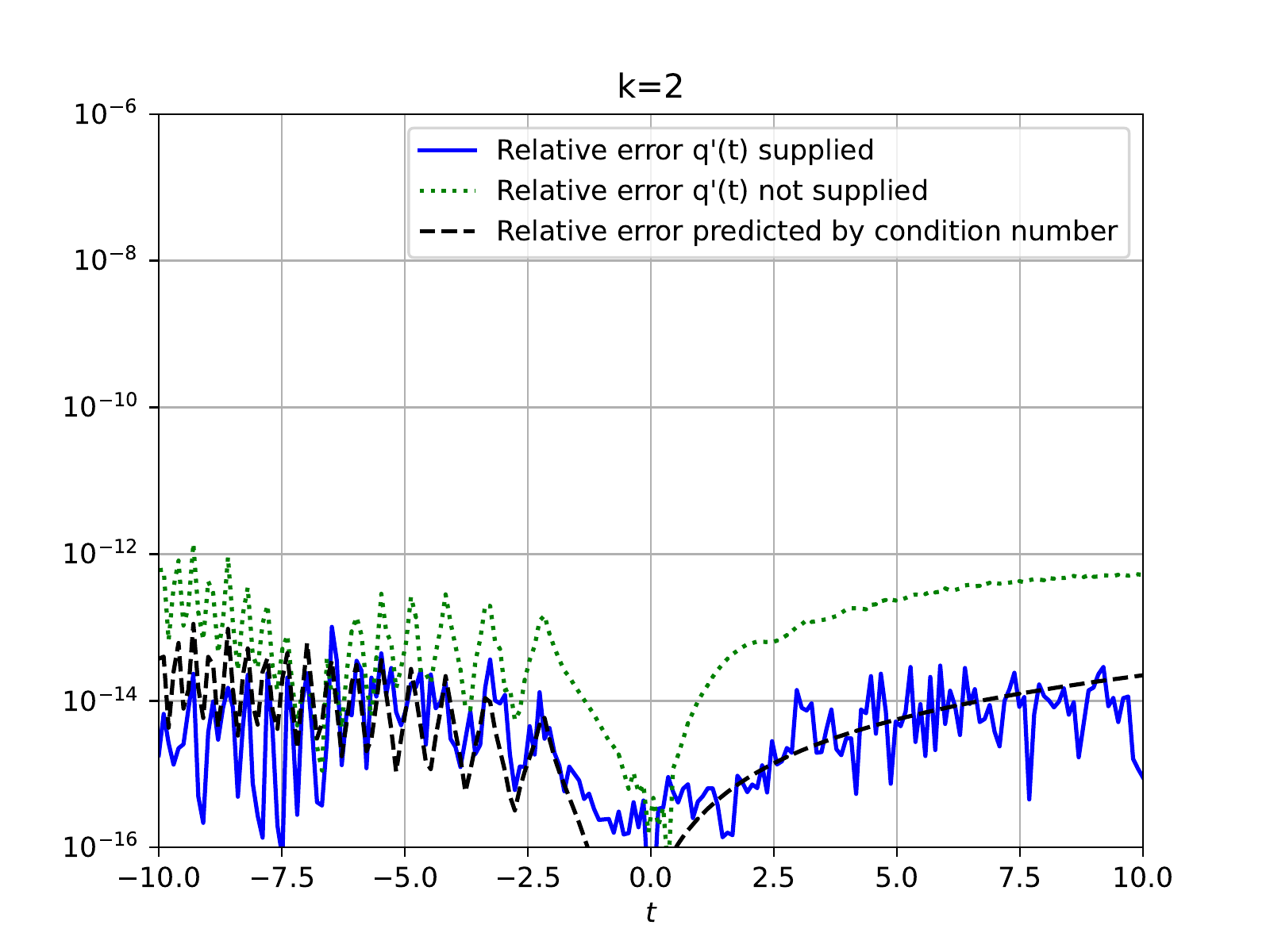}
\hfil
\includegraphics[width=.40\textwidth]{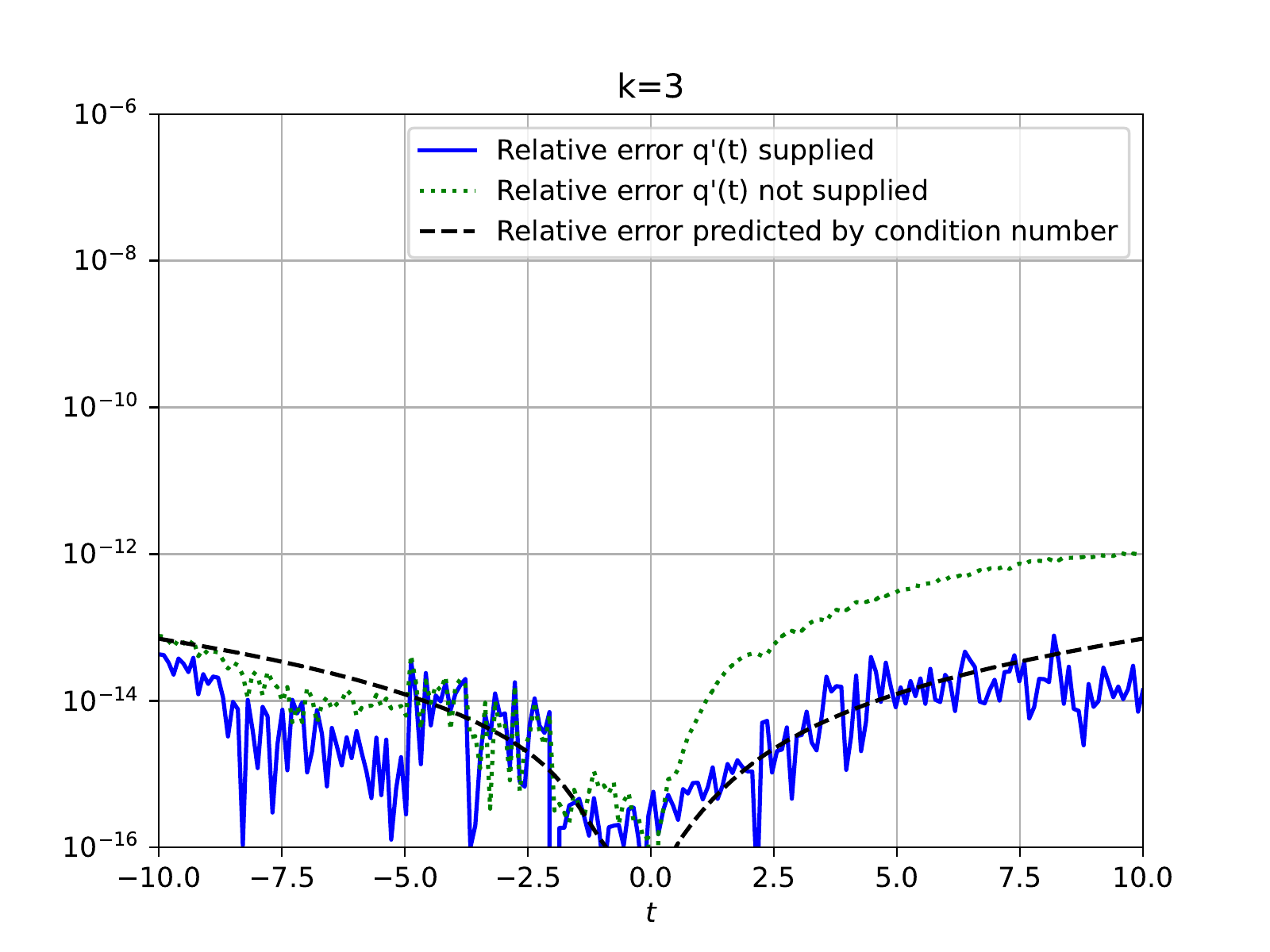}
\hfil

\hfil
\includegraphics[width=.40\textwidth]{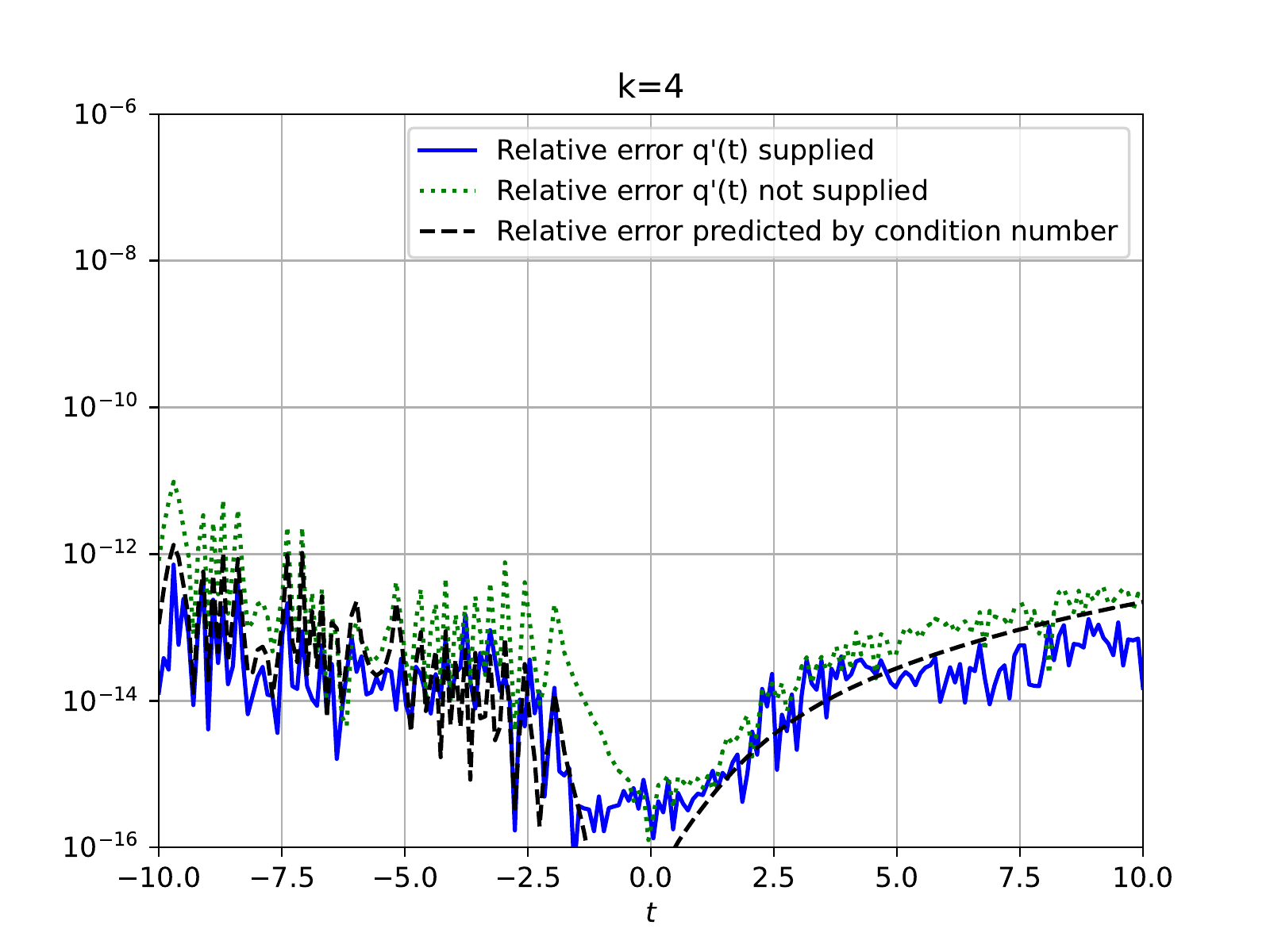}
\hfil
\includegraphics[width=.40\textwidth]{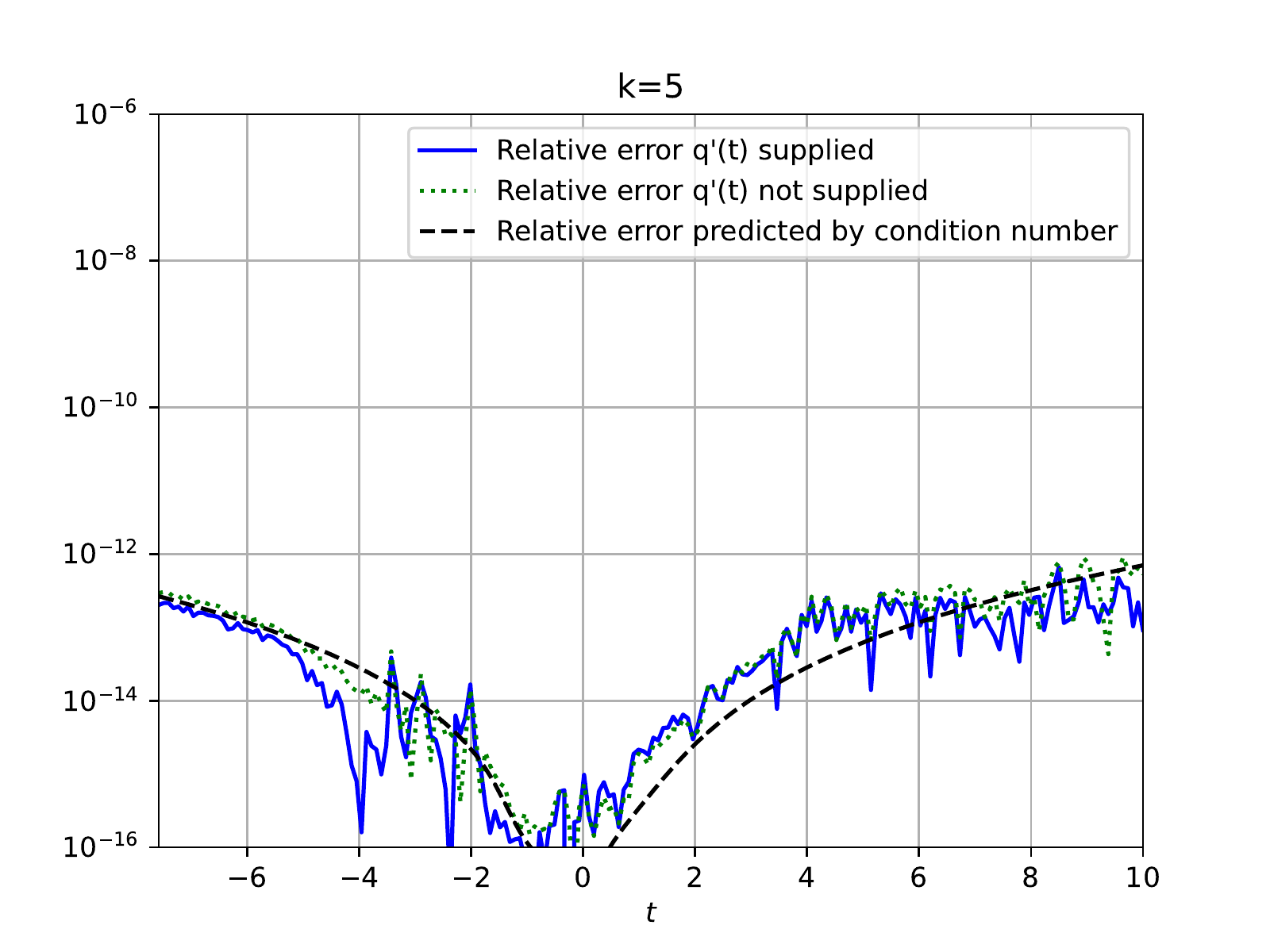}
\hfil

\hfil
\includegraphics[width=.40\textwidth]{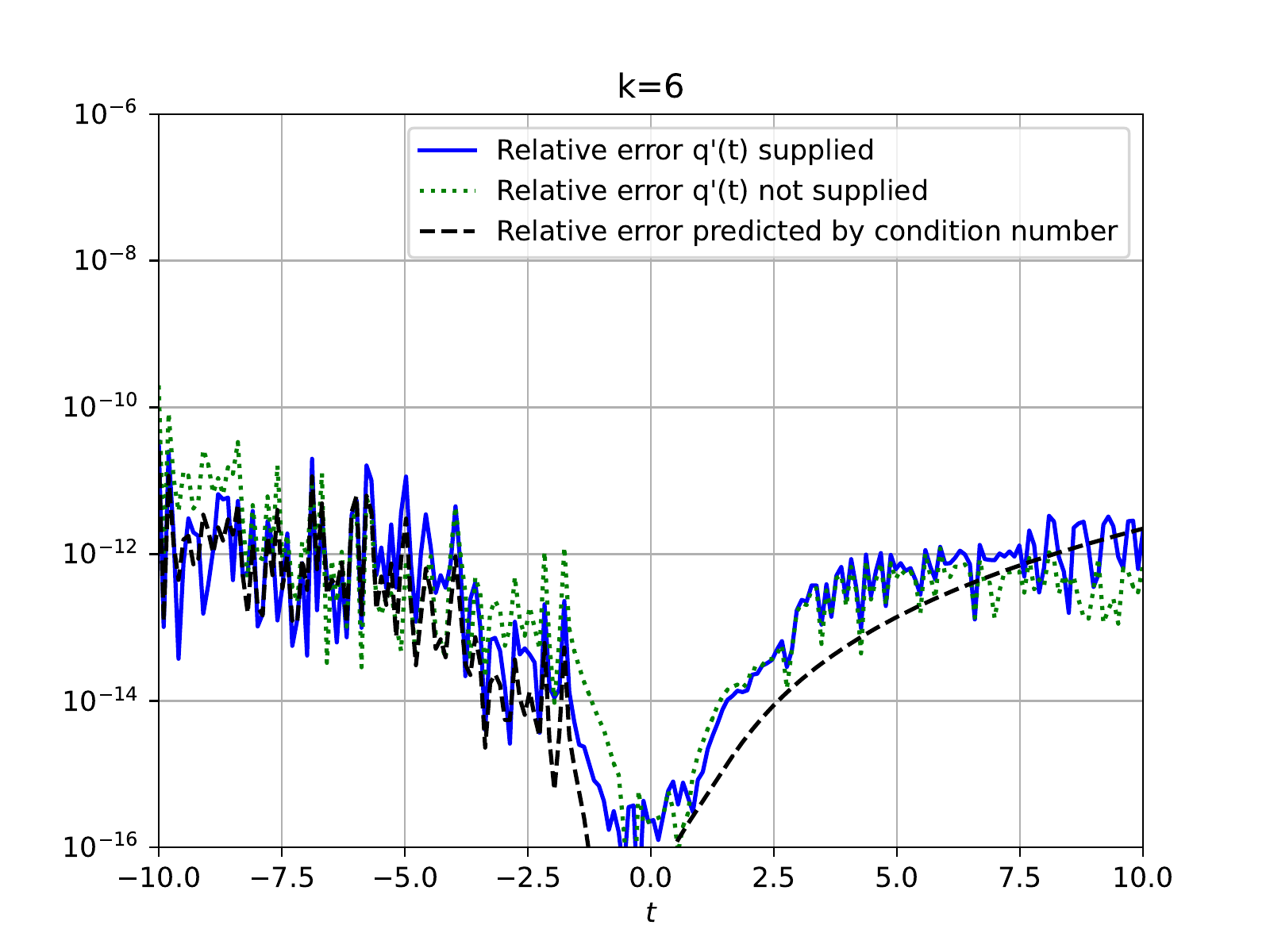}
\hfil
\includegraphics[width=.40\textwidth]{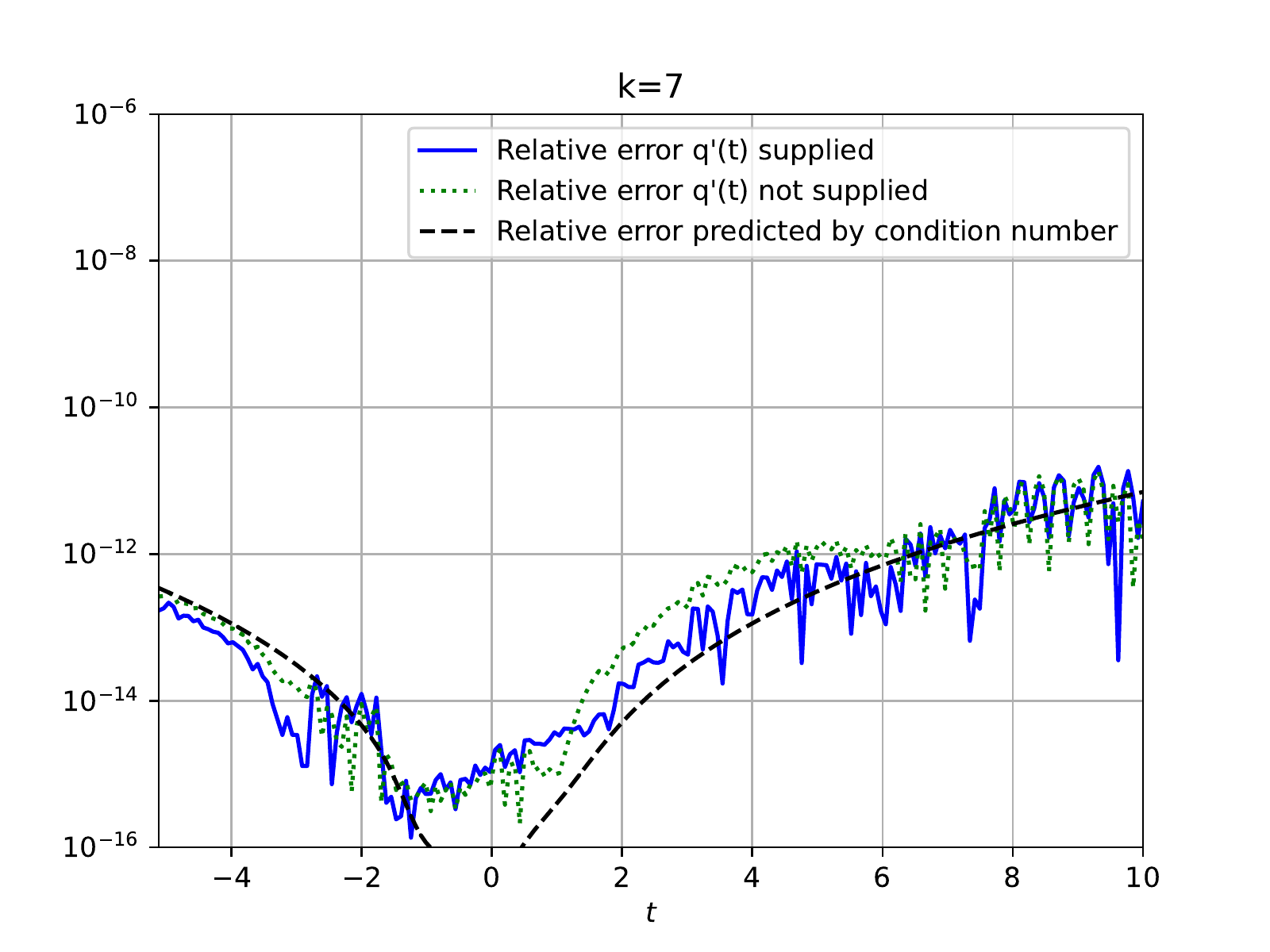}
\hfil

\hfil
\includegraphics[width=.40\textwidth]{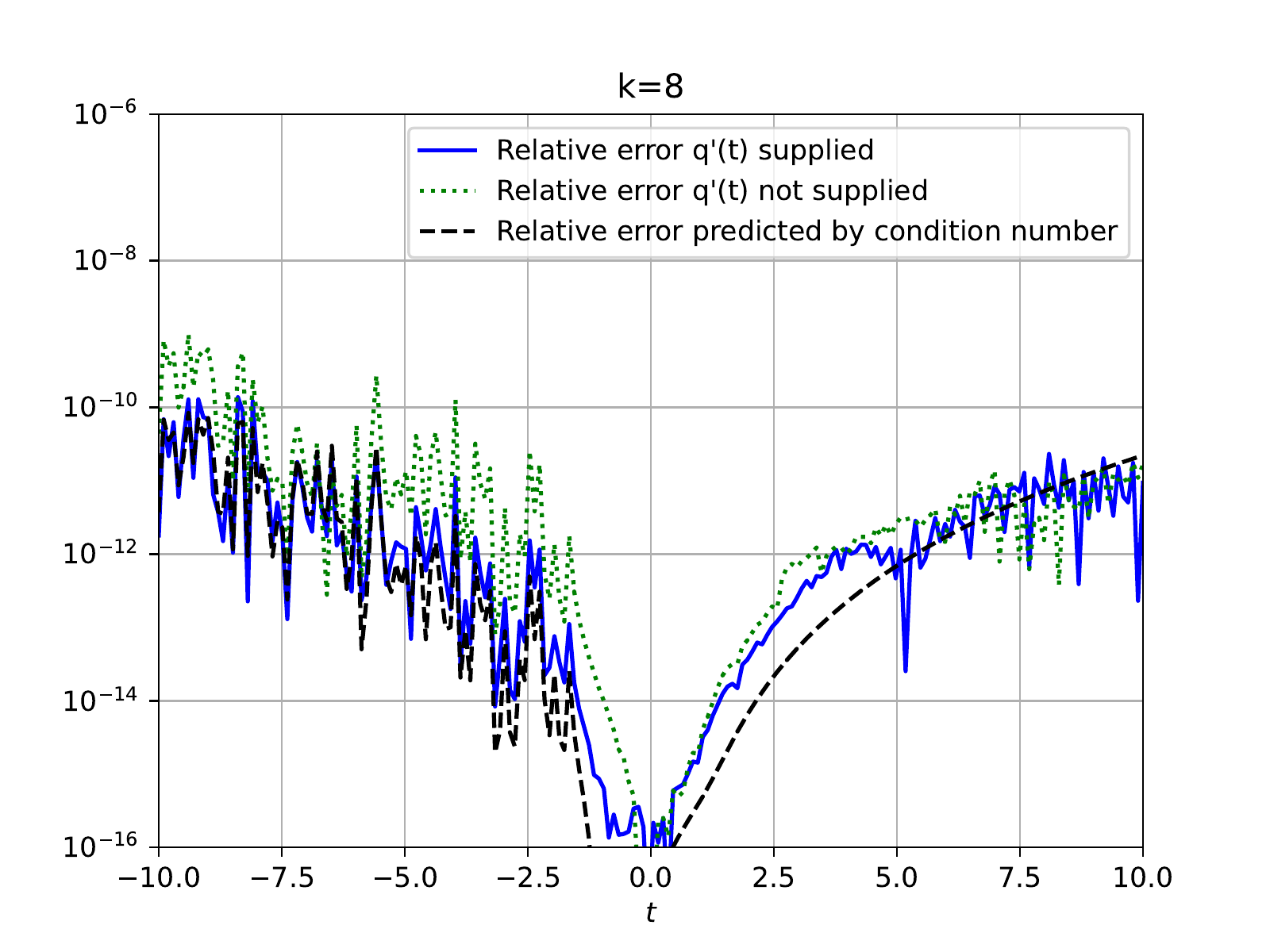}
\hfil
\includegraphics[width=.40\textwidth]{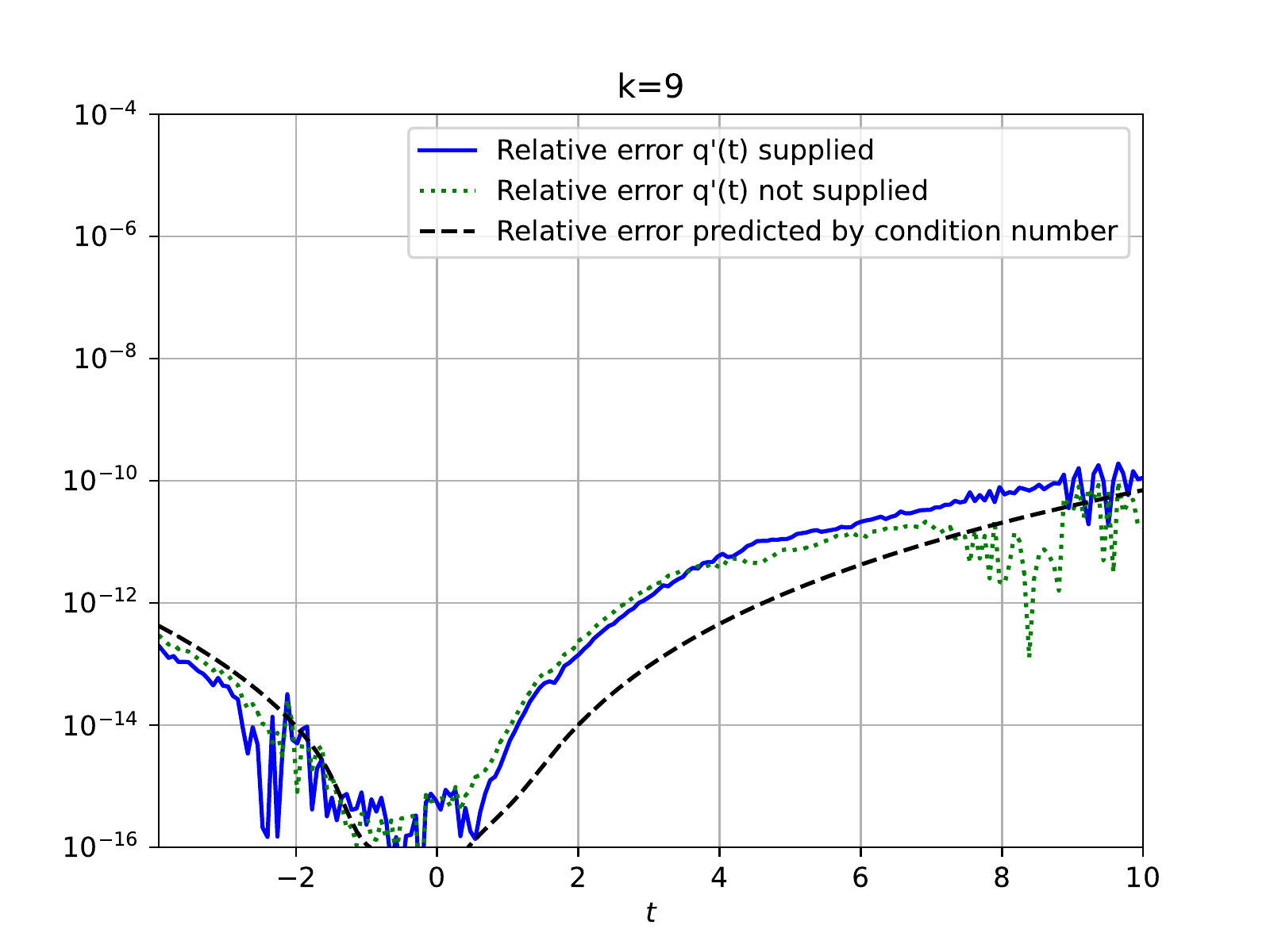}
\hfil

\caption{\small  The results of the experiments discussed in  Subsection~\ref{section:experiments:high},
which concerned an equation with a turning point of order $k$.
Each plot corresponds to one value of $k$ and gives the relative errors incurred when the function  $f_k^{\mbox{\tiny high}}$ 
defined via (\ref{experiments:high:f}) was evaluated at $200$ equispaced points using the algorithm
of this paper with the derivatives of the coefficients supplied, the 
relative errors  incurred when $f_k^{\mbox{\tiny high}}$  was evaluated at $200$
equispaced points using the algorithm of this paper with the derivatives of the coefficients calculated
via spectral differentiation, and the relative errors predicted by
the condition number of evaluation of $f_k^{\mbox{\tiny high}}$ .
}
\label{figure:highplots1}

\end{figure}

\vfil\eject

\begin{figure}[h!!!!!!]
\hfil
\includegraphics[width=.42\textwidth]{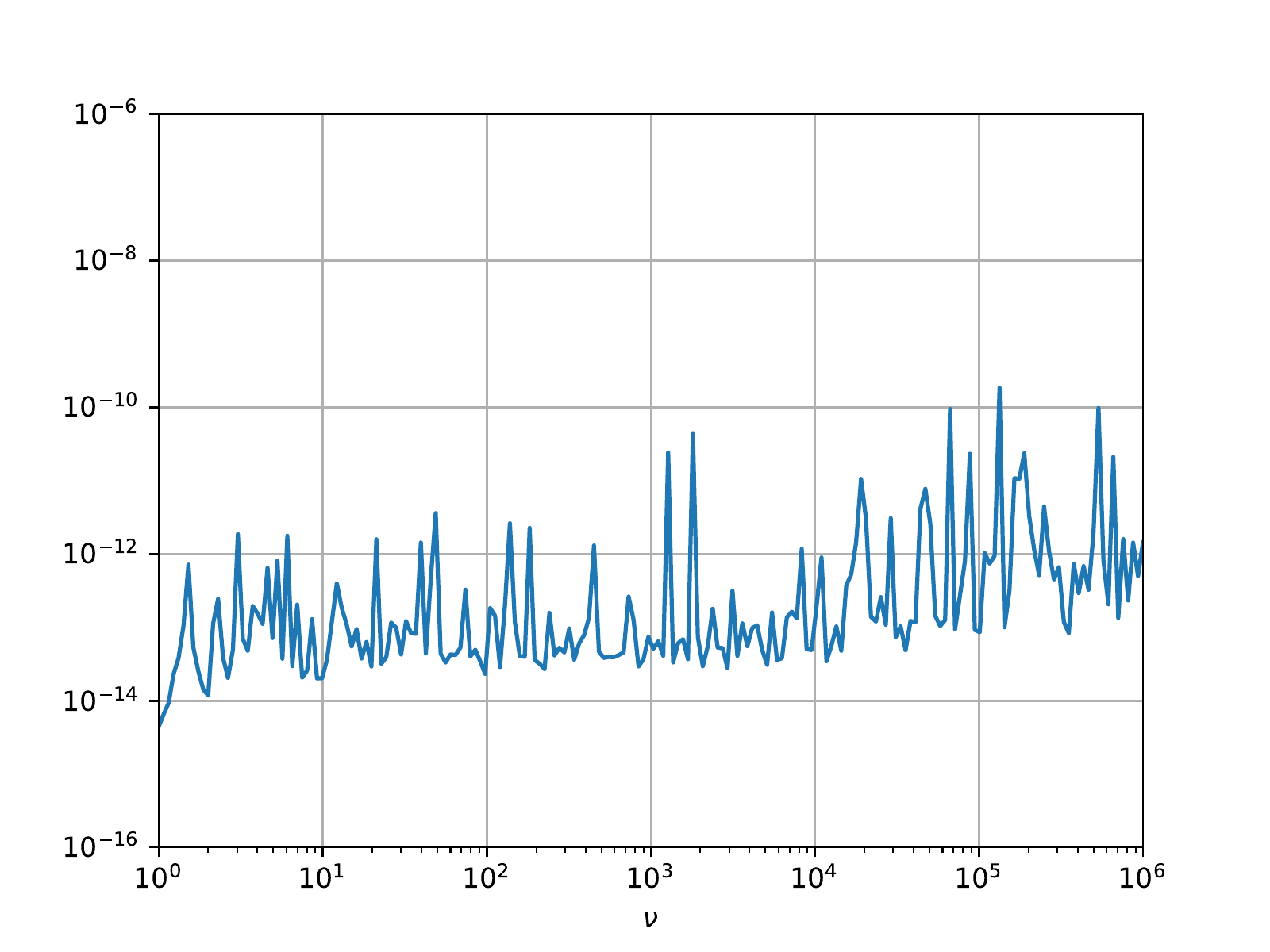}
\hfil
\includegraphics[width=.42\textwidth]{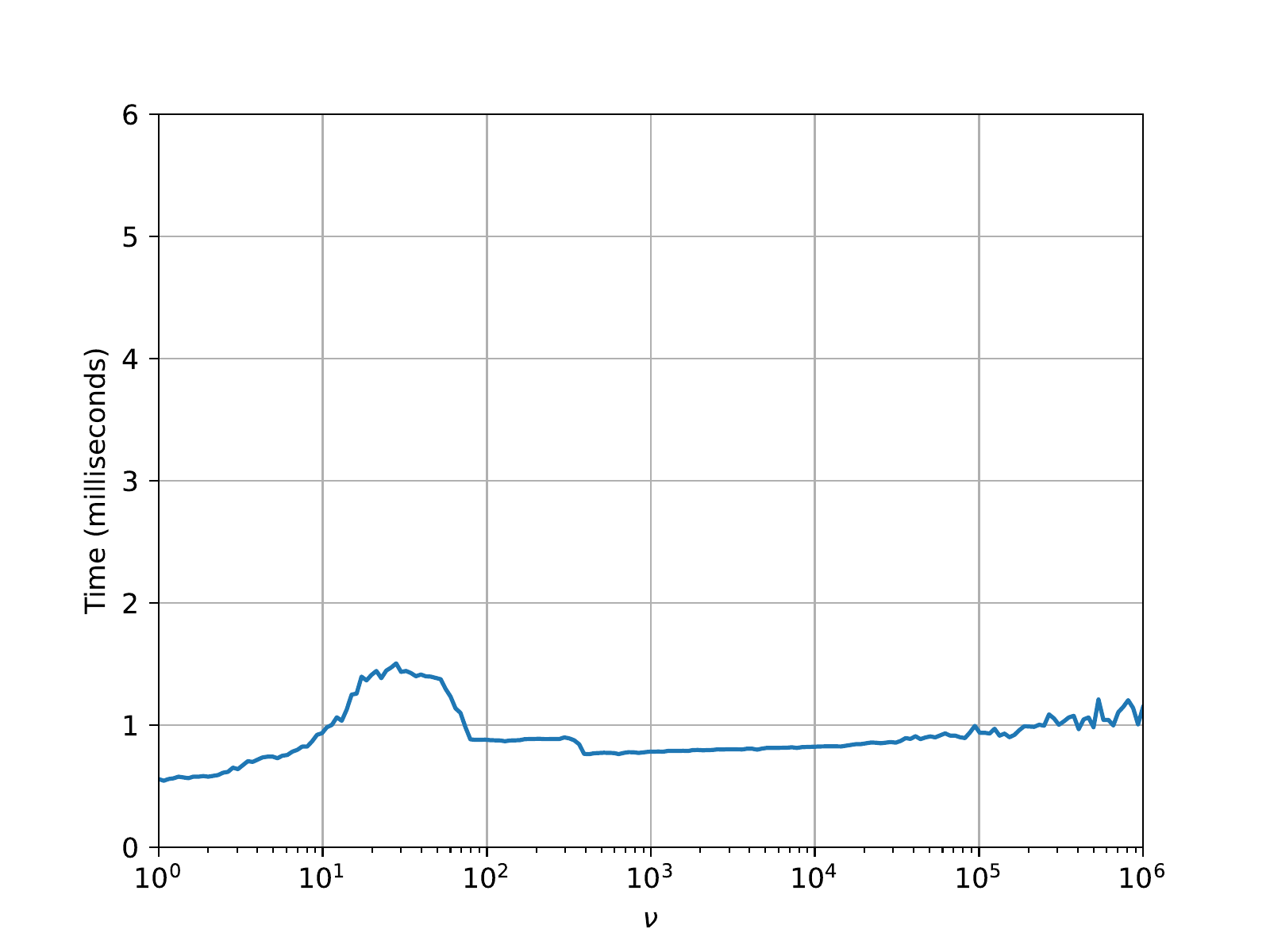}
\hfil

\caption{The results of the experiments discussed in  Subsection~\ref{section:experiments:bumps},
which concerned an equation whose coefficient has two bumps.
On the left we give the maximum absolute difference between the solution
obtained with the phase method and that obtained using a standard solver running in extended precision
 as a function
of $\nu$.  On the right is a plot giving the time taken to solve the problem (\ref{experiments:bumps:problem})
using the phase method as a function of $\nu$.}
\label{figure:bumpsplots1}
\end{figure}

\begin{figure}[!h]

\hfil
\includegraphics[width=.42\textwidth]{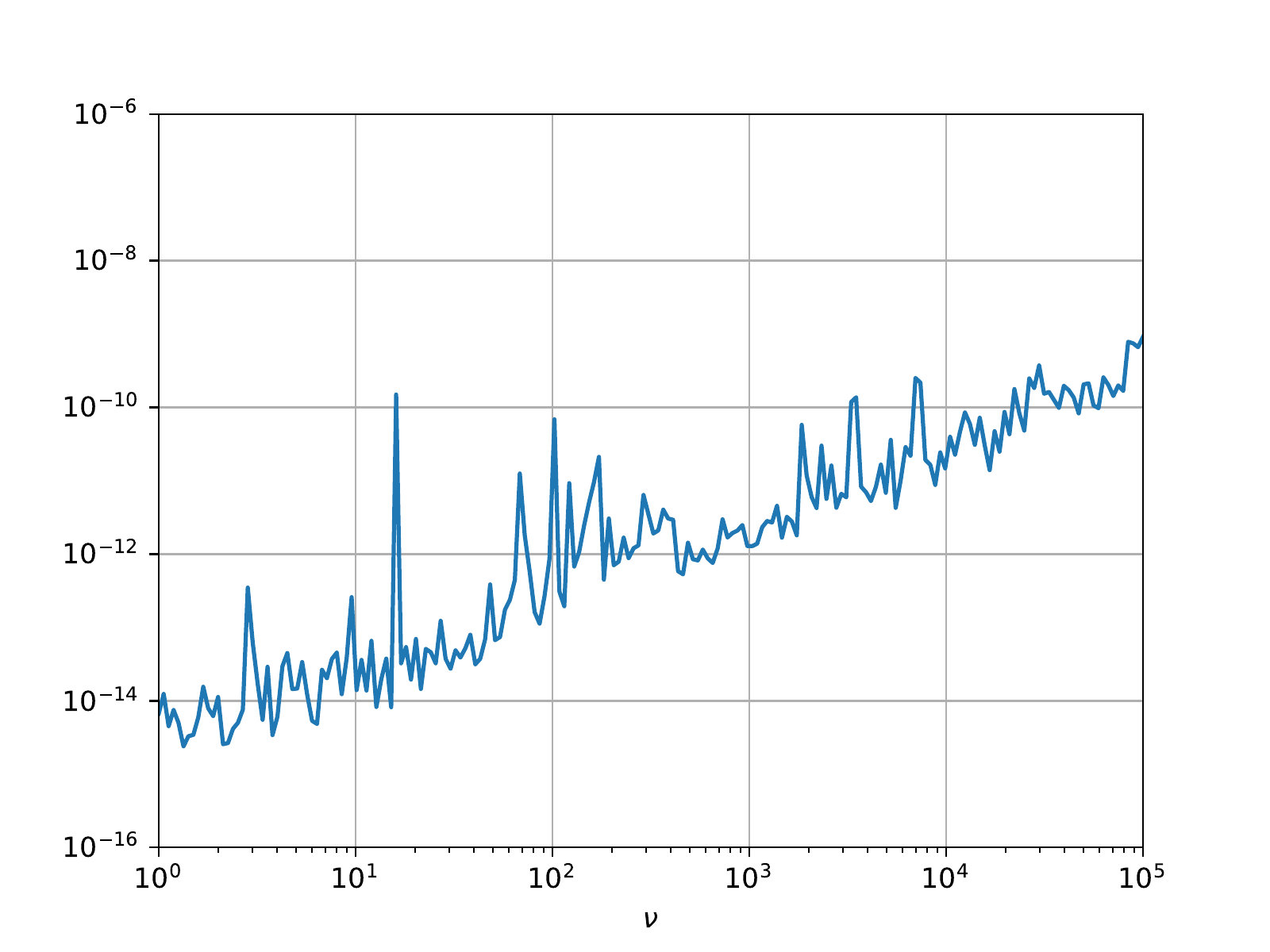}
\hfil
\includegraphics[width=.42\textwidth]{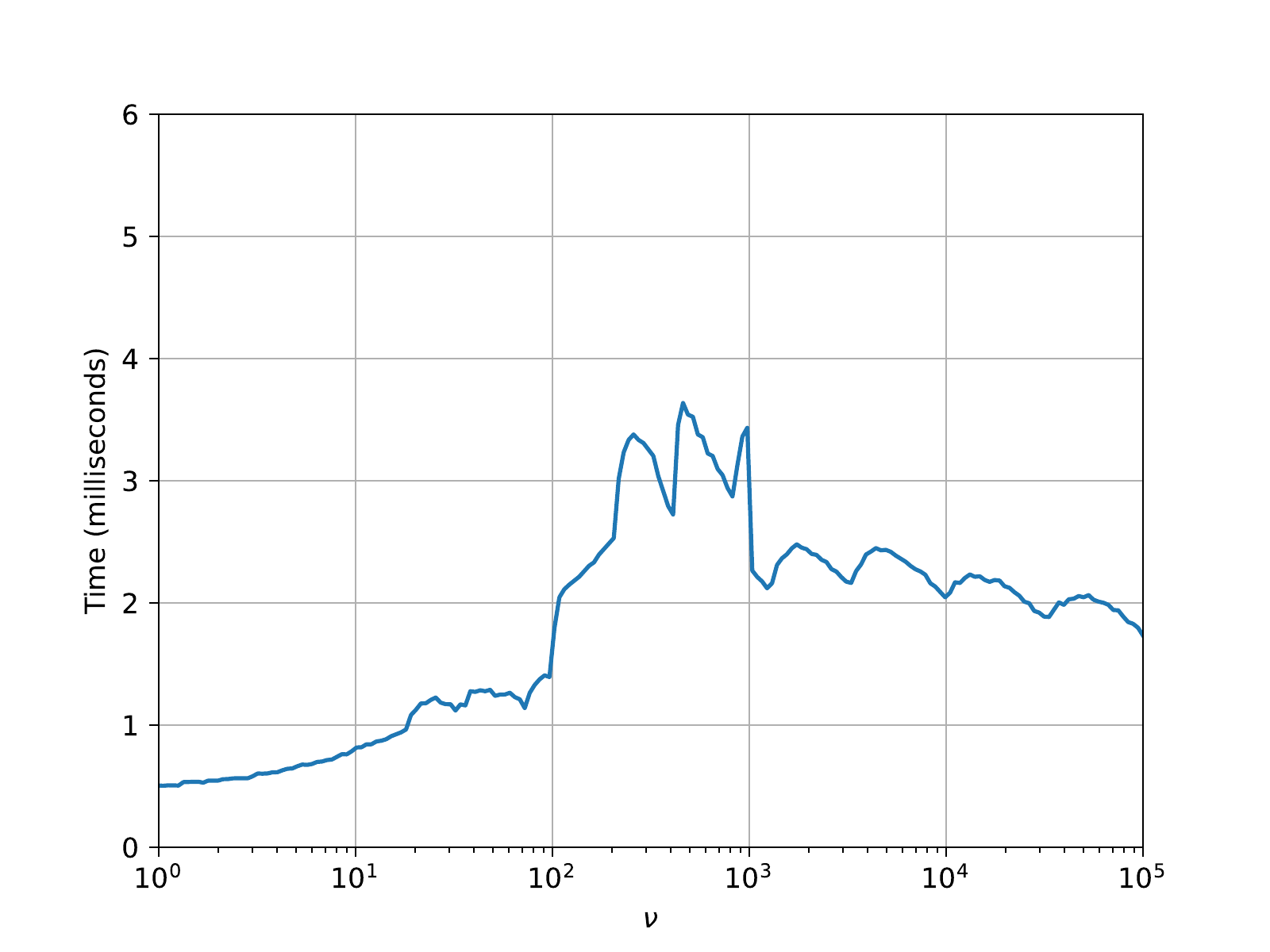}
\hfil

\caption{The results of the experiments discussed in  Subsection~\ref{section:experiments:two},
which concerned an equation whose coefficient has three turning points.
On the left is a plot of a measure of the maximum difference between the solution obtained
via a phase function method and that obtained using a standard solver
running in extended precision
 as a function of $\nu$.  On the right is a plot giving the time taken to solve the problem (\ref{experiments:bumps:problem})
using the phase method as a function of $\nu$.}
\label{figure:twoplots1}
\end{figure}

\begin{figure}[!h]

\hfil
\includegraphics[width=.42\textwidth]{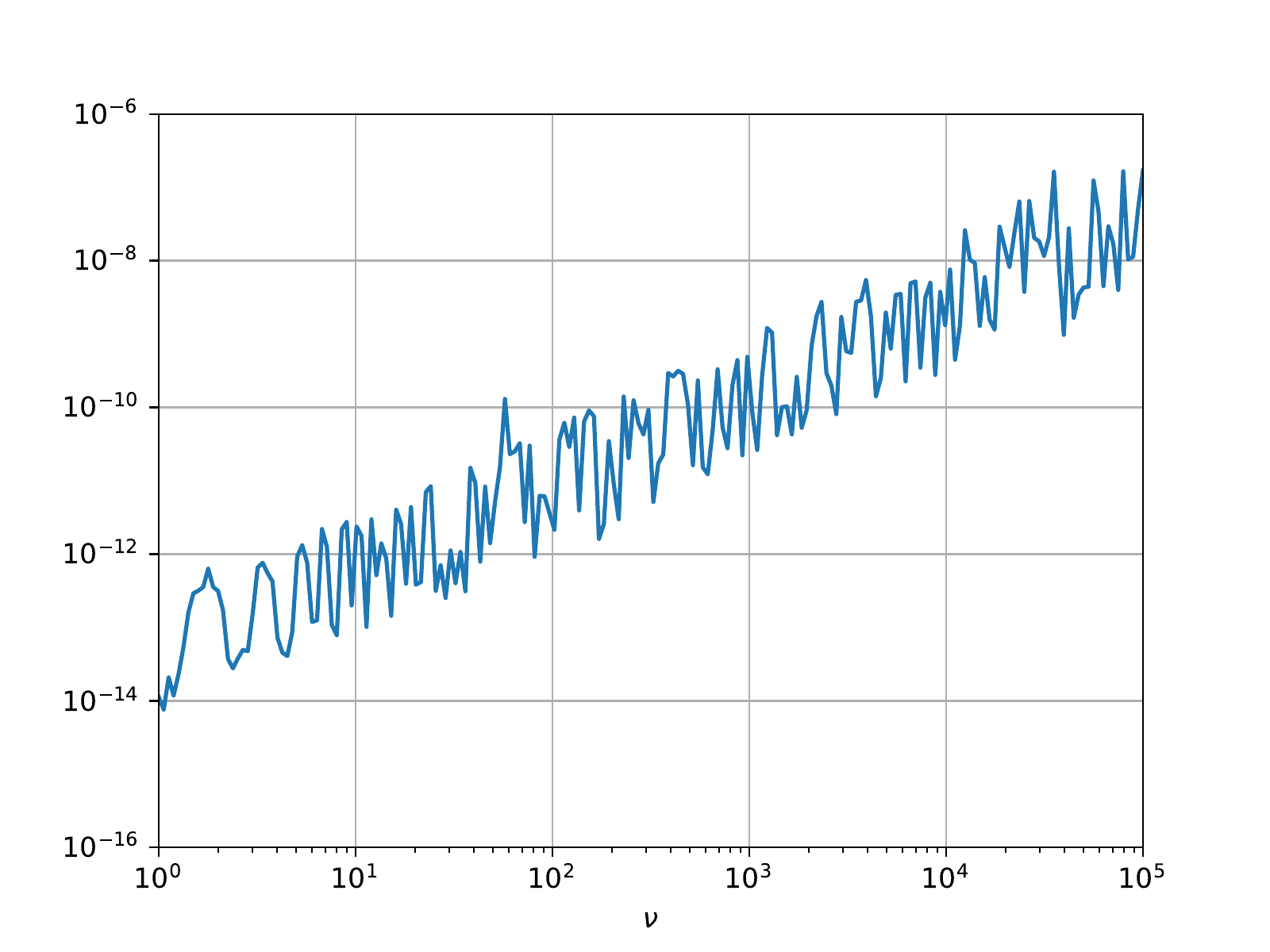}
\hfil
\includegraphics[width=.42\textwidth]{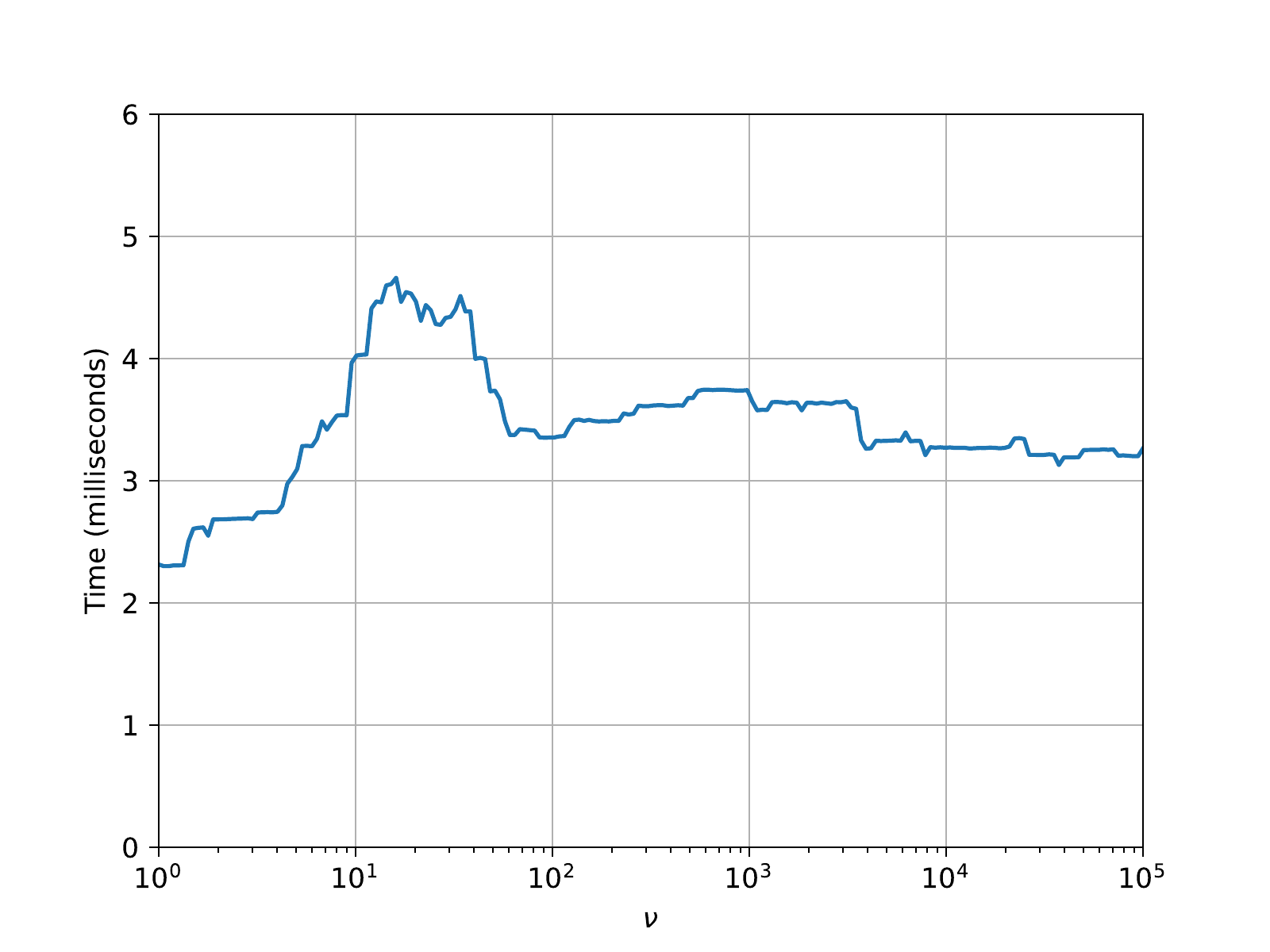}
\hfil

\caption{The results of the experiments discussed in  Subsection~\ref{section:experiments:many},
which concerned an equation whose coefficient has $12$ turning points.
On the left is a plot of the maximum absolute difference between the solution
obtained with the phase method and that obtained via a standard solver running in extended precision
 as a function of $\nu$.  On the right is a plot giving the time taken to solve the problem (\ref{experiments:bumps:problem})
using the phase method as a function of $\nu$.}
\label{figure:manyplots1}
\end{figure}

\end{document}